\documentclass[letterpaper,oneside]{amsart} 
\pdfoutput=1 
\usepackage{todonotes}
\setuptodonotes{inline}
\usepackage{color,soul}
\usepackage{amsrefs}
\usepackage[pagewise]{lineno}
\usepackage{xcolor}
\usepackage[utf8]{inputenc} 
\usepackage{comment}
\usepackage{mathtools} 
\usepackage{colonequals}
\usepackage{amsmath}
\usepackage{amssymb}
\usepackage{tikz}
\usepackage{mathrsfs}
\usepackage{enumitem}
\usepackage[foot]{amsaddr}
\usepackage{mathtools}

\usepackage{mleftright}

\usepackage[english]{babel} 
\usepackage{tikz}
\usepackage{pdfpages}
\usepackage{thmtools} 
\usepackage{amssymb}

\usepackage[all]{xy} 
\usepackage{hyperref}

\newcommand{\R}{\mathbb{R}}
\newcommand{\Z}{\mathbb{Z}}

\newcommand{\N}{\mathbb{N}}

\newcommand{\catname}[1]{\mathbf{#1}}

\newcommand{\Q}{\mathbb{Q}}

\title{Semi-coarse Spaces, Homotopy and Homology}
\author{Antonio Rieser and Jonathan Treviño-Marroquín}

\thanks{This material is based in part upon work supported by the US National
Science Foundation under Grant No. DMS-1928930 while the authors participated
in a program supported by the Mathematical Sciences Research Institute. The
program was held in the summer of 2022 in partnership with the the Universidad
Nacional Aut{\'o}noma de M{\'e}xico. This work was also supported by the
CONACYT Investigadoras y Investigadores por M{\'e}xico Project \#1076, the CONACYT Ciencia de Fronteras grant CF-2019-217392, and by
the grant N62909-19-1-2134 from the US Office of Naval Research Global and the
Southern Office of Aerospace Research and Development of the US Air Force
Office of Scientific Research. The second author was also supported by the CONACYT postgraduate studies scholarship number 839062. }

\theoremstyle{plain}
\newtheorem{teorema}{Theorem}[subsection]
\newtheorem{proposicion}[teorema]{Proposition}
\newtheorem{corolario}[teorema]{Corollary}
\newtheorem{lema}[teorema]{Lemma}

\theoremstyle{definition}
\newtheorem{definicion}[teorema]{Definition}
\newtheorem{definition}[teorema]{Definition}
\newtheorem{observacion}[teorema]{Remark}
\newtheorem{ejemplo}{Example}[teorema]

\newtheorem{innercustomthm}{Theorem}
\newenvironment{customthm}[1]
  {\renewcommand\theinnercustomthm{#1}\innercustomthm}
  {\endinnercustomthm}

\newtheorem{innercustomprop}{Proposition}
\newenvironment{customprop}[1]
  {\renewcommand\theinnercustomprop{#1}\innercustomprop}
  {\endinnercustomprop}
  
\newtheorem{innercustomlem}{Lemma}
\newenvironment{customlem}[1]
  {\renewcommand\theinnercustomlem{#1}\innercustomlem}
  {\endinnercustomlem}

\declaretheoremstyle[
qed=\qedsymbol
]{mystyle}

\begin{document}
\begin{abstract}

We begin the study the algebraic topology of semi-coarse spaces, 
which are generalizations of coarse spaces that enable one to endow non-trivial 
`coarse-like' structures to compact metric spaces,
something which is impossible in coarse geometry. We first study
homotopy in this context, and we then construct homology groups which are
invariant under semi-coarse homotopy equivalence. We further show that any
undirected graph $G=(V,E)$ induces a semi-coarse structure on its set of vertices $V_G$, and
that the respective semi-coarse homology is isomorphic to the
Vietoris-Rips homology. This, in turn, leads to a homotopy invariance theorem for
the Vietoris-Rips homology of undirected graphs. 
\end{abstract}

\maketitle

\tableofcontents

\section{Introduction}
\pagenumbering{arabic}\label{cap.introduccion}

Coarse geometry \cite{Roe_2003} is often referred to as `geometry in the large', or `large
scale geometry', for the reason that, on a metric space endowed with the natural coarse
structure, all bounded phenomena are trivial from the coarse point of view.
The study of coarse geometry has led to a number of important advances, in particular in group theory,
where it has enabled many questions about groups to be addressed geometrically
via the analysis of their Cayley graphs (see, for instance, \cite{Drutu_Kapovich_2018} for a
recent book-length discussion of this approach). Nonetheless,
the large-scale nature of coarse geometry forces all finite-diameter spaces to
be coarsely-equivalent to a point. It is not difficult, however, to imagine
scenarios where one might like to `coarsen' a space up to a certain scale, but
where one does not want to erase all bounded phenomena at every scale. Many
problems of topological data analysis, for instance, in which one would like to
deduce the topological invariants of a space from algebraic invariants built
from a finite subset of the space, present themselves naturally as problems of
`medium-scale' coarsening. Several possibilities for a formal context for
`coarsening' a space up to a scale were developed by the first author in
\cite{Rieser_2021} and \cite{Rieser_arXiv_VR_2020}, where he studied the problem from the point of 
view of \v Cech closure spaces and semi-uniform spaces, respectively. The connection to coarse 
geometry in these earlier works, however, remained unclear, and this question forms the 
primary motivation for the work presented here.

Whereas the closure structures and semi-uniform structures studied in \cite{Rieser_2021} 
and \cite{Rieser_arXiv_VR_2020} are generalizations
of topological structures and uniform structures, respectively, in this article, we focus instead 
on a generalization of coarse structures,
modifying the axioms of coarse geometry \cite{Roe_2003} to allow `coarsenings'
only up to a preferred size. Doing so links the homotopy of spaces coarsened up to a fixed scale to the study of coarse spaces, and, in particular, the development of homotopy on semi-coarse spaces that we have begun here has been used by the second author in \cite{Trevino_2024_FundGroupoid_arXiv} to construct new homotopical invariants on coarse spaces themselves, something which is unclear how to accomplish through either closure or semi-uniform spaces directly. Technically,
this `coarsening up to a scale' is achieved by eliminating the product axiom in the definition of the
coarse structure, which, indeed, has a similar flavor to the elimination of
the idempotence axiom in the passage from Kurotowski (topological) closure
structures to more general \v Cech closure structures, or to the elimination
of the product axiom in the passage from uniform spaces to semi-uniform spaces.  While removing the
product axiom effectively destroys the notion of coarse equivalence, we observe in the following
that coarse equivalence can be profitably generalized to a useful notion of homotopy
equivalence in this setting.  The resulting
structures and spaces are called \emph{semi-coarse structures} and \emph{semi-coarse spaces}, respectively, and were first introduced in \cite{Zava_2019}.
Many examples of semi-coarse spaces exist. The examples of most interest to topological data analysis may be built from pseudo-metric spaces
and a preferred scale $r>0$, but many others exist as well, in
particular those constructed from semi-uniform spaces. Furthermore, when specialized to graphs on lattices, the semi-coarse homotopy introduced here is also similar to the digital homotopy studied in \cites{Staecker_2021,Boxer_Staecker_2018,Boxer_Staecker_2017,Boxer_Staecker_2016,Lupton_etal_2023,Lupton_Scoville_2022,Lupton_etal_2022a,Lupton_etal_2022,Lupton_etal_2021}. 

The outline of this article is as follows. In Section \ref{sec:Semi-Coarse
Spaces}, we give the basic definitions, point-set properties, and principal
examples of semi-coarse spaces. We will also show in that section how to
obtain a coarse space from a semi-coarse space through a certain limiting
process. In Section \ref{sec:Homotopy}, we begin the study of homotopy
invariants in the semi-coarse category. Unfortunately, the interval
does not appear to have a natural non-trivial semi-coarse structure,
significantly complicating the construction. We are able circumvent this
shortcoming by adapting the homotopy construction from \cite{Babson_etal_2006}, 
using finite-length subsequences of $\Z$ in place of the
interval to construct cylinders. While technically delicate, this nonetheless 
allows us to define homotopy groups, prove a long-exact sequence of pairs, and
compute the fundamental group of a `semi-coarse circle' with four points. 
Finally, in Section \ref{sec:Homology}, we construct homology groups for
semi-coarse spaces, inspired by the Vietoris-Rips construction now commonly
used in topological data analysis. We give a construction for any semi-coarse
space, demonstrate the invariance of the homology groups with respect to
the homotopy introduced in Section \ref{sec:Homotopy}, and show that, for a countable
semi-coarse space, this homology is exactly the Vietoris-Rips homology of an 
associated graph. 

Many of the results in this paper were first presented in the Master's thesis of the second author
\cite{Trevino_2019}, written under the supervision of the first author, where \emph{semi-coarse} spaces 
and structures were called \emph{pseudo-coarse}.

\section{Semi-Coarse Spaces}
\label{sec:Semi-Coarse Spaces}

In this section, we define semi-coarse spaces, bornologous functions, and we give examples of 
semi-coarse structures constructed from a metric space and a scale parameter $r>0$. We then 
define semi-coarse
quotients, disjoint unions, and products, and we show how to build an induced
coarse structure from a semi-coarse structure. 

\subsection{Fundamental Concepts and Examples}
\label{SecConcepFund}

We begin by setting some notation which we will use throughout the article.

\begin{definition}
Let $X$ be a set. We denote by $\mathcal{P}(X)$ the collection of all
subsets of $X$, and 

\begin{enumerate}[wide]
\item $X^n$ will denote $\overbrace{X\times X \times \cdots \times X}^{n-\mbox{times}}$.
\item $\Delta_X:=\{(x,x)\in X\times X \mid x\in X\}$ will be called \emph{the diagonal of $X$}.
\item For $V \in \mathcal{P}(X\times X)$, we define \[V^{-1}:=\{(y,x)\in X\times
    X: (x,y)\in V\},\] which we call \emph{the inverse of $V$}.
\item For $V,W\in\mathcal{P}(X\times X)$, we define \[ V\circ W:=\{(x,y)\in
            X\times X \mid
    \exists z\in X,(x,z)\in V\mbox{ and }(z,y)\in W\},\] 
    which will be called \emph{the set product of $V$ and $W$}.

    Given an integer $n\geq 2$, we write $V^{\circ n}$ for 
    $\overbrace{V\circ\ldots\circ V}^{n-\mbox{times}}$.
\item Let $X$ and $Y$ be sets, $f:X \to Y$ a set function, and $V\in\mathcal{P}(X\times X)$.
    Then \[(f\times f)(V):=\{(f(x),f(x'))\in Y\times Y \mid (x,x')\in V\},\]  and we call 
    $(f \times f)(V)$ \emph{the image of $V$ under $f\times f$}.
\end{enumerate}
\end{definition}

We now define semi-coarse spaces, our principal object of interest. 

\begin{definicion}[Semi-coarse space]
\label{def:Semi-coarse} \index{Semi-coarse}
Let $X$ be a set, and let $\mathcal{V}\subset\mathcal{P}(X\times X)$
be a collection of subsets of $X \times X$ which satisfies 
\begin{enumerate}[label=(sc\arabic*)]
    \item \label{sc1}$\Delta_X \in \mathcal{V}$,
    \item \label{sc2} If $B\in\mathcal{V}$ and $A\subset B$, then $A\in\mathcal{V}$,
    \item \label{sc3} If $A,B\in \mathcal{V}$, then $A\cup B\in\mathcal{V}$,
    \item \label{sc4} If $A\in\mathcal{V}$, then $A^{-1}\in\mathcal{V}$.
\end{enumerate}
We call $\mathcal{V}$ a \emph{semi-coarse structure on $X$}, and we say
that the pair $(X,\mathcal{V})$ is a \emph{semi-coarse space}. 

If, in addition, $\mathcal{V}$ satisfies
\begin{enumerate}[label=(sc\arabic*),resume]
    \item\label{item:Coarse condition} If $A,B\in\mathcal{V}$, then $A\circ B\in\mathcal{V}$.
\end{enumerate}
then $\mathcal{V}$ will be called a \emph{coarse structure}, and
$(X,\mathcal{V})$ will be called a \emph{coarse space}, as in \cite{Roe_2003}. 

\end{definicion}

The elements of $\mathcal{V}$ will be called \emph{controlled sets}. Moreover,
if there exist $a,b$ such that $\{(a,b)\}\in\mathcal{V}$ we will say that \emph{$a$ and $b$ are adjacent}. 
If $\mathcal{V}$ and $\mathcal{V}'$ are semi-coarse structures on $X$ such that $\mathcal{V} \subset \mathcal{V}'$, then we say that \emph{$\mathcal{V}'$ is finer than $\mathcal{V}$} and
\emph{$\mathcal{V}$ is coarser than $\mathcal{V}'$}. Finally, when the structure $\mathcal{V}$ is 
unambiguous, we will sometimes refer to the semi-coarse space only by $X$.

The functions of interest between semi-coarse spaces will be those which preserve the semi-coarse 
structure.

\begin{definicion}
    We will say that $f:X\rightarrow Y$ is a \emph{$(\mathcal{V},\mathcal{W})$-bornologous function}, or simply \emph{bornologous}, if $f\times f$ maps each controlled set $V \in \mathcal{V}$ to a controlled set
    $(f \times f)(V) \in \mathcal{W}$.
\end{definicion}

Since the composition of set maps is associative, the composition of bornologous maps is bornologous, 
and the identity is bornologous for every semi-coarse space $(X,\mathcal{V})$, we have

\begin{teorema}
    Semi-coarse spaces and bornologous functions form a category.
\end{teorema}

We denote the category of semi-coarse spaces and bornologous functions by $\catname{SCoarse}$.

For our first collection of examples, we show how undirected graphs induce semi-coarse structures. 
We begin with a pair of lemmas.

\begin{lema}
    \label{lem:Qc from power set}
    Let $X$ be a set, and let $W \subset X \times X$ such that
    \begin{enumerate}
        \item $\Delta_X \subset W$,
        \item $W = W^{-1}$.
    \end{enumerate}
    Then $(X,\mathcal{P}(W))$ is a semi-coarse space.
\end{lema}
\begin{proof}
    We verify directly that the axioms \ref{sc1} - \ref{sc4} from Definition \ref{def:Semi-coarse} are
    satisfied.
    \begin{enumerate}[label=(sc\arabic*)]
        \item $\Delta_X \subset W$, so $\Delta_X \in \mathcal{P}(W)$.
        \item $\mathcal{P}(W)$ is closed under taking subsets by definition.
                \item Since any two sets $A,B \in \mathcal{P}(W)$ are subsets of $W$, 
                    the union $A \cup B \subset W$, and therefore $A \cup B \in \mathcal{P}(W)$.
                \item Suppose $A \in \mathcal{P}(W)$. Then $A \subset W$. However, by hypothesis on  
            $W$, if $(a,b) \in W$ then $(b,a) \in W$. Therefore $A^{-1} \subset W$, so
            $A^{-1} \in \mathcal{P}(W)$ as well.
    \end{enumerate}
    It follows from the above that $(X,\mathcal{P}(W))$ is a semi-coarse space, as desired.
\end{proof}

It will sometimes be convenient to use this lemma in the following alternative form.

\begin{lema}
    \label{lem:Qc from power set of union}
    Let $X$ be a set and suppose that $\mathcal{U} \subset \mathcal{P}(X \times X)$ is a collection
of subsets of $X \times X$ such that
    \begin{enumerate}
        \item $U \in \mathcal{U} \implies U^{-1} \in \mathcal{U}$, and
        \item $\Delta_X \in \mathcal{U}$.
    \end{enumerate}
    Let $W \coloneqq \bigcup_{U\in \mathcal{U}} U$.
    Then $(X,\mathcal{P}(W))$ is a semi-coarse space.
\end{lema}

\begin{proof}
    By construction, $W$ satisfies the hypotheses of Lemma \ref{lem:Qc from power set}. The conclusion
    follows.
\end{proof}

We now use these lemmas to give several important examples of semi-coarse spaces.

\begin{ejemplo}
\label{ex:Qc from graphs}
Let $G = (V,E)$ be an undirected graph (i.e. $(u,v) \in E \iff (v,u) \in E$). 
Define $\mathcal{V} \subset \mathcal{P}(V \times V)$ to be
\begin{equation*}
    \mathcal{V}_G \coloneqq \mathcal{P}(E \cup \Delta_V),
\end{equation*}
where $\Delta_V$ is the diagonal of $V \times V$. By Lemma \ref{lem:Qc from power set of union} above, 
the pair $(V,\mathcal{V}_G)$ is a semi-coarse space.
\end{ejemplo}

\begin{definicion}
	\label{def:SC from graphs}
    Given a graph $G=(V,E)$, we say that the semi-coarse space $(V,\mathcal{V}_G)$ constructed in 
    Example \ref{ex:Qc from graphs} is \emph{the semi-coarse space generated by the graph $G$}.
\end{definicion}

Another important class of examples may be constructed from metric spaces combined with a positive
scale parameter $r>0$.

\begin{ejemplo}
    \label{ex:Qc from metric}
    \begin{enumerate}[wide, align=left, leftmargin=0pt, labelindent=0pt,listparindent=0pt, itemindent=!,
        label=(\alph*)]
    \item \label{ex:Qc from metric less}Let $(X,d)$ be a metric space. Let $r>0$ be a positive real number and define
    \begin{equation*}
        U_r \coloneqq \{ (x,x') \in X \times X \mid d(x,x') \leq r \},
    \end{equation*}
and let $\mathcal{V}_r \coloneqq \mathcal{P}(U_r)$.
    Then $(X,\mathcal{V}_r)$ is a semi-coarse space by Lemma \ref{lem:Qc from power set}.
\item \label{ex:Qc from metric leq} Similarly, defining $U^<_r$ by 
        \begin{equation*}
            U^<_r \coloneqq \{ (x,x') \in X \times X \mid d(x,x') < r\},
        \end{equation*}
 and let $\mathcal{V}^<_r \coloneqq \mathcal{P}(U^<_r)$.
        Lemma \ref{lem:Qc from power set} gives that $(X,\mathcal{V}^<_r)$ is a semi-coarse space.
    \end{enumerate}
\end{ejemplo}

For the next example, which generalizes the ones in Example \ref{ex:Qc from metric}
above, we introduce semi-pseudometric spaces.

\begin{definicion}[Semi-Pseudometric; \cite{Cech_1966}, 18 A.1.]
\label{SemPseMet}
Let $X$ be a set and $d:X\times X\rightarrow \mathbb{R}$ be a function, we will say $d$ is a 
\emph{semi-pseudometric on $X$} if they satisfies the next conditions
\begin{enumerate}[label=(m\arabic*)]
    \item \label{item:m1} For each $x\in X$, $d(x,x)=0$.
\item For every $x,y\in X$, $d(x,y)=d(y,x)\geq 0$.
\end{enumerate}
A semi-pseudometric on $X$ is a \emph{pseudometric on $X$} if also
\begin{enumerate}[label=(m\arabic*),resume]
\item For each $x,y,z\in X$, $d(x,z)\leq d(x,y)+d(y,z)$, i.e. $d$ satisfies the triangle inequality.
\end{enumerate}
A semi-pseudometric will be called a \emph{semi-metric}, if it also satisfies
\begin{enumerate}[label=(m\arabic*),resume]
    \item \label{item:m4} $d(x,y)=0$ implies $x=y$.
\end{enumerate}
Finally, a semi-pseudometric $d$ will be called a \emph{metric} iff it satisfies \ref{item:m1}-
\ref{item:m4}, that is, iff $d$ is both a semi-metric and 
a pseudometric.

A \emph{semi-pseudometric space} is an ordered pair $(X,d)$ where $X$ is a set
and $d$ is a semi-pseudometric on $X$. Similarly, $(X,d)$ is a semi-metric, pseudometric,
or metric space when $d$ is a semi-metric, pseudometric, or metric, respectively. 

\end{definicion}

We will also need the following semi-pseudometrics constructed from a metric and a
pre-determined scale $r >0$.

\begin{definicion}
\label{def:Metric to semipseudometric}
Let $(X,d)$ be a metric space and let $r\geq 0$ be a non-negative real number. We will define the 
functions $d_r, d^<_r, d^\leq_r:X\times X\rightarrow [0,\infty)$ by
\begin{align*}
    d^{\leq}_r(x,y)=& \begin{cases} 0 & \text{if } d(x,y)\leq r,\\
        1& \text{if } d(x,y)> r,
    \end{cases}\\
        d^{<}_r(x,y) =& \begin{cases} 0 &\text{if } d(x,y) < r,\\
        1 & \text{if } d(x,y) \geq r.
    \end{cases}\\
            d_r(x,y) =& \max\{0,d(x,y)-r\}
\end{align*}
\end{definicion}

\begin{observacion}
    For $r\geq 0$, the function $d_r$ satisfies $d_r(x,x)=0$ and $d_r(x,y)=d_r(y,x)$,  
    by the symmetry of $d$. Therefore, $d_r$ is a semi-pseudometric on $X$ for any $r \geq 0$. Similarly,
    the functions $d^{\leq}_r$ and $d^<_r$ are semi-pseudometrics on $X$ for $r \leq 0$ and $r>0$,
    respectively. 
\end{observacion}
\begin{ejemplo}
    \label{ex:Qc from semipseudometric}
    Let $d:X \times X \to [0,\infty)$ be a semi-pseudometric,
    and define $U \coloneqq \{(x,x') \in X \times X \mid d(x,x') = 0\}$. Then $(X,\mathcal{P}(U))$ is
    a semi-coarse space by Lemma \ref{lem:Qc from power set}. The semi-pseudometrics $d^\leq_r$ and 
    $d^<_{r}$ from Definition \ref{def:Metric to semipseudometric} give the
    semi-coarse structures in Examples \ref{ex:Qc from metric}\ref{ex:Qc from metric less} and 
    \ref{ex:Qc from metric}\ref{ex:Qc from metric leq},
    respectively.
\end{ejemplo}

The following examples are independent of the Lemmas \ref{lem:Qc from power set} and \ref{lem:Qc from power set of union}, and give examples of semi-coarse spaces which are not coarse.

\begin{ejemplo}
    \begin{enumerate}[wide, align=left, leftmargin=0pt, labelindent=0pt,listparindent=0pt, itemindent=!,
        label=(\alph*)]

\item Let $(X,d)$ be a metric space which is at least countably infinite.
    Define $\mathcal{V} \subset \mathcal{P}(X \times X)$ by
    \begin{equation*}
        V \in \mathcal{V} \iff \{ |V| < \infty \mid \forall (x,x') \in V,\; d(x,x') < 1 \}.
    \end{equation*}
    Then $(X,\mathcal{V})$ is a semi-coarse space. If $X = \Q$ with the Euclidean metric,
    then $(\Q, \mathcal{V})$ is a semi-coarse space which is not coarse.
\item Let $(X,c)$ be an uncountable metric space. Define $\mathcal{V} \subset \mathcal{P}(X \times X)$
    by
    \begin{equation*}
        V \in \mathcal{V} \iff \{ |V| \leq \aleph_0 \mid \forall (x,x') \in V,\; d(x,x') < 1 \}.
    \end{equation*}
    Then $(X,\mathcal{V})$ is a semi-coarse space. If $X = \R$ with the Euclidean metric,
    then $(\R,\mathcal{V})$ is a semi-coarse space which is not coarse.
    \end{enumerate}
\end{ejemplo}

\subsection{Subspaces}
\label{SecSubEsp}
We will now proceed to build new semi-coarse spaces from existing ones, in particular,
constructing subspaces, products, and quotient spaces in the semi-coarse category. 
Our first step, in this section, will be to construct semi-coarse subspaces of a 
semi-coarse space $(X,\mathcal{V})$.

We first require the following lemma.

\begin{lema}
\label{InvUnionInter}
Let $X$ be a set, and let $\{A_\lambda\}_{\lambda\in\Lambda}$ a collection of elements of 
$\mathcal{P}(X\times X)$, where $\Lambda$ is an index set. Then
\begin{enumerate}[wide,label=(\roman*)]
    \item[(i)] $\bigcup\limits_{\lambda\in\Lambda} (A_\lambda)^{-1}= \left( \bigcup\limits_{\lambda\in\Lambda} A_\lambda \right)^{-1}$,
\item[(ii)] $\bigcap\limits_{\lambda\in\Lambda} (A_\lambda)^{-1}= \left( \bigcap\limits_{\lambda\in\Lambda} A_\lambda \right)^{-1}$.
\end{enumerate}
\end{lema}

\begin{proof}
Let $X$ be a set, let $\Lambda$ be a index set, and let $\{A_\lambda\}_{\lambda\in\Lambda}$ be a
collection of elements of $\mathcal{P}(X\times X)$.

\begin{enumerate}[wide,label=(\roman*)]
\item We see that 
    \begin{align*}
        (x,y)\in \bigcup\limits_{\lambda\in\Lambda} (A_\lambda)^{-1} \iff & \exists \lambda_0\in 
        \Lambda \text{ such that } (x,y)\in (A_{\lambda_0})^{-1}
        \iff (y,x)\in A_{\lambda_0}\\
        \iff & (y,x)\in\bigcup\limits_{\lambda\in\Lambda} A_\lambda
        \iff (x,y)\in\left(\bigcup\limits_{\lambda\in\Lambda} A_\lambda \right)^{-1}.
    \end{align*}
    Therefore $\bigcup\limits_{\lambda \in \Lambda} (A_\lambda)^{-1} =
    \left(\bigcup\limits_{\lambda \in \Lambda} A_\lambda \right)^{-1}$ 
\item Now note that
    \begin{align*}
        (x,y)\in \bigcap\limits_{\lambda\in\Lambda} (A_\lambda)^{-1} \iff &
        \forall \lambda\in \Lambda, (x,y)\in
    (A_{\lambda})^{-1}
    \iff  (y,x)\in A_{\lambda} \forall \lambda\in\Lambda\\
        \iff &(y,x)\in\bigcap\limits_{\lambda\in\Lambda} A_\lambda 
    \iff (x,y)\in\left(\bigcap\limits_{\lambda\in\Lambda} A_\lambda
    \right)^{-1}. 
    \end{align*}
    Therefore, $\bigcap\limits_{\lambda\in\Lambda} (A_\lambda)^{-1}= \left( \bigcap\limits_{\lambda\in\Lambda} A_\lambda \right)^{-1}$ as desired.\qedhere
    \end{enumerate}
\end{proof}

\begin{proposicion}
    \label{prop:Subspace}
Let $(X,\mathcal{V})$ be a semi-coarse space and suppose that $Y \subset X$. Define the
collection $\mathcal{V}_Y \in \mathcal{P}(Y\times Y)$ by
\begin{equation*}
    \mathcal{V}_Y \coloneqq \{ V \cap (Y\times Y) \mid V \in \mathcal{V}\}
\end{equation*}
Then the pair $(Y,\mathcal{V}_Y)$ is a semi-coarse space.
\end{proposicion}

\begin{proof}
Let $(X,\mathcal{V})$, $Y \subset X$, and $\mathcal{V}_Y$ be as in the statement of the 
proposition. We check that the axioms for a semi-coarse space are satisfied by 
$(Y,\mathcal{V}_Y)$.
\begin{enumerate}[wide,label=(sc\arabic*)]
\item We observe that $\Delta_Y=\Delta_X\cap (Y\times Y)$. Therefore $\Delta_Y\in \mathcal{V}_Y$.
\item Let $A\in\mathcal{V}_Y$ and suppose that $A'\subset A$. By definition, 
    there is a $B\in\mathcal{V}$ such that $A=B\cap (Y\times Y)$. We define $B':=A'\cap B$
    and we observe that $B'\in \mathcal{V}$ since $B' \subset B \in \mathcal{V}$.
    However, $B'\cap (Y\times Y)= A'\cap B\cap (Y\times Y) = A'\cap A = A'$, 
    so $A'\in \mathcal{V}_Y$.
\item If $A,B\in \mathcal{V}_Y$, then there are $A',B'\in\mathcal{V}$ such that
    $A=A'\cap (Y\times Y)$ and $B=B'\cap (Y\times Y)$. Therefore
    \begin{align*}
        A\cup B= & (A' \cap (Y \times Y))\cup (B' \cap Y \times Y)\\
        =& (A'\cup B') \cap (Y \times Y)
    \end{align*}
    Since $A' \cup B' \in \mathcal{V}$, it follows that $A\cup B \in \mathcal{V}_Y$.
\item If $A\in \mathcal{V}_Y$, then there is $A'\in\mathcal{V}_Y$ such that
    $A=A'\cap (Y\times Y)$, then by \autoref{InvUnionInter} we get 
    $A^{-1}=(A\cap(Y\times Y)^{-1}=A'^{-1}\cap (Y\times Y)^{-1}=A'^{-1}\cap (Y\times Y)$. 
    Therefore $A^{-1}\in \mathcal{V}_Y$.
\end{enumerate}
It follows that $(Y,\mathcal{V}_Y)$ is a semi-coarse space, as desired.
\end{proof}

\begin{definicion}[Semi-Coarse Subspace]
\label{DefSubPseCoar}
Let $(X,\mathcal{V})$ be a semi-coarse space and let $Y \subset X$. 
The ordered pair $(Y,\mathcal{V}_Y)$ from \autoref{prop:Subspace} will be called 
a \emph{semi-coarse subspace of $X$}. When 
the structure $\mathcal{V}_Y$ is clear from the context, we will simply refer to
the subspace $(Y,\mathcal{V}_Y)$ as $Y$.
\end{definicion}

The following proposition gives a useful criterion for checking whether a function is bornologous on a semi-coarse
space $(X,\mathcal{V})$.

\begin{proposicion}
    \label{prop:Bornologous function from subspaces}
  Let $(X, \mathcal{V})$ and $(Y, \mathcal{W})$ be semi-coarse spaces, and
  suppose that $(X_i, \mathcal{V}_i) \subset (X, \mathcal{V})$, $i \in \{ 1,
  \ldots, n \},$ are subspaces of $(X, \mathcal{V})$ such that $\cup_{i = 1}^n
  X_i = X$ and every set $V \in \mathcal{V}$ may be written in the form
  \begin{equation*}
      V = \bigcup_{i = n}^n V_i
\end{equation*}
where each $V_i \in \mathcal{V}_i$. Now suppose that $f : X
  \rightarrow Y$ is a map such that the restrictions $f |_{X_i} : (X_i,
  \mathcal{V}_i) \rightarrow (Y, \mathcal{W})$ are bornologous for all $i \in
  \{ 1, \ldots, n \}$. Then $f : (X, \mathcal{V}) \rightarrow (Y,
  \mathcal{W})$ is bornologous. 
\end{proposicion}

\begin{proof}
  Let $V \in \mathcal{V} .$ Then, by hypothesis, $V \in \mathcal{V }_i$ for
  some $i \in \{ 1, \ldots, n \}$, and since $f |_{X_i}$ is bornologous, we
  have that $(f \times f) (V) = (f |_{X_i} \times f |_{X_i}) (V) \in
  \mathcal{W}$. Since $V \in \mathcal{V}$ is arbitrary, it follows that $f$ is
  bornologous.
\end{proof}

\subsection{Product Semi-Coarse Spaces}
\label{SecProdPseuCoar}

In this section, we construct product semi-coarse structure on the product of sets. We
begin with the following definition and several preliminary results.

\begin{definicion}
Let $X$ and $Y$ be sets, and suppose that $V\in \mathcal{P}(X\times X)$ and 
$W\in \mathcal{P}(Y\times Y)$. We define
\begin{align*}
    V \boxtimes W := & \{ ((a,b),(c,d))\in (X\times Y) \times (X \times Y):
    (a,c)\in V, (b,d)\in W \}
\end{align*}
which we call the \emph{Cartesian cross product} of $V$ and $W$.
\end{definicion}

\begin{lema}
    \label{lem:Cross product associative}
    The Cartesian cross product is associative.
\end{lema}
\begin{proof}
    Let $X, Y,$ and $Z$ be sets, and suppose that $U \in \mathcal{P}(X\times X)$, 
    $V\in\mathcal{P}(Y\times Y)$, and $W\in\mathcal{P}(Z\times Z)$. We wish to show that
    $(V\boxtimes W)\boxtimes U \cong V \boxtimes (W\boxtimes U)$. By definition, we
    have that
    \begin{align*}
        (U \boxtimes V) \boxtimes W =\, &\{ (((a,c),e),((b,d),f)) \in ((X \times Y) \times Z)
            \times ((X \times Y) \times Z)\\
            & \mid \,
        (a,b) \in (X \times X), (c,d) \in (Y \times Y), \text{ and } (e,f) \in (Z \times Z) \},
        \text{ and}\\
            U \boxtimes (V \boxtimes W) =& \{ (((a,(c,e)),(b,(d,f))) \in (X \times 
                (Y \times Z)) \times (X \times (Y \times Z))
                \\ & \mid (a,b) \in (X \times X), (c,d) \in (Y \times Y), 
            \text{ and } (e,f) \in (Z \times Z) \}.
    \end{align*}
    Since the sets $(X \times Y) \times Z \cong X \times (Y \times Z)$ are isomorphic, 
    this proves the result.
\end{proof}

\begin{lema}
\label{boxunion}
Let $X$ and $Y$ be sets, $A\subset X\times X$, $B\subset Y\times Y$, 
$\{A_\lambda\}_{\lambda\in\Lambda}$ a collection of subsets of $X\times X$ 
indexed by $\Lambda$ and $\{B_\gamma\}_{\gamma\in\Gamma}$ a collection of subsets of 
$Y\times Y$ indexed by $\Gamma$. Then
\begin{enumerate}[wide,label=(\roman*)]
\item $(A\boxtimes B)^{-1}=A^{-1}\boxtimes B^{-1}$.
\item $\bigcup_{\lambda\in\Lambda} (A_\lambda \boxtimes B) = \left( \bigcup\limits_{\lambda\in\Lambda} A_\lambda \right)\boxtimes B$.
\item $\bigcup_{\gamma\in\Gamma} (A \boxtimes B_\gamma) = A \boxtimes \left( \bigcup\limits_{\gamma\in\Gamma} B_\gamma \right)$.
\end{enumerate}
\end{lema}

\begin{proof}
Let $X$, $Y$, $A\subset X\times X$, $B\subset Y\times Y$, $\{A_\lambda\}_{\lambda\in\Lambda}$,
and $\{B_\gamma\}_{\gamma\in\Gamma}$ be as in the statement of the proposition. 
\begin{enumerate}[wide,label=(\roman*)]
\item We observe that
\begin{align*}
    (A\boxtimes B)^{-1} = & \{ ((x',y'),(x,y)) : (x,x')\in A, (y,y')\in B\}\\
    = & \{ ((x',y'),(x,y)) : (x',x)\in A^{-1},(y', y)\in B^{-1}\}\\
                    = & A^{-1} \boxtimes B^{-1}.
\end{align*}
\item\label{item:Box product union 2} 
    If $((a,c),(b,d))\in \bigcup\limits_{\lambda\in\Lambda} (A_\lambda \boxtimes B)$,
    then there exists a $\lambda_0\in\Lambda$ such that $((a,b),(c,d))\in
    A_{\lambda_0}\boxtimes B$, so $(a,b)\in A_{\lambda_0}$ and $(c,d)\in
    B$. This implies that $(a,b)\in \bigcup\limits_{\lambda\in\Lambda} A_\lambda$ and
    $(b,d)\in B$, from which we conclude that $((a,c),(b,d))\in \left(
    \bigcup\limits_{\lambda\in\Lambda} A_\lambda \right)\boxtimes B$.

    Conversely, if $((a,c),(b,d))\in \left( \bigcup\limits_{\lambda\in\Lambda}
    A_\lambda \right)\boxtimes B$, then $(a,b)\in
    \bigcup\limits_{\lambda\in\Lambda} A_\lambda$ and $(c,d)\in B$. This implies that there
    is $\lambda_0\in\Lambda$ such that $(a,b)\in A_{\lambda_0}$, which gives that
    $((a,c),(b,d))\in A_{\lambda_0}\boxtimes B$. This implies that $((a,c),(b,d))\in
    \bigcup\limits_{\lambda\in\Lambda} (A_\lambda\boxtimes B)$, as desired. 

\item The proof of this point is
    analogous to the proof of part \ref{item:Box product union 2}
    above.\qedhere
\end{enumerate}
\end{proof}

In the next proposition, we construct the product of two semi-coarse spaces.

\begin{proposicion}
    \label{prop:Product QC structure}
Let $(X,\mathcal{V})$ and $(Y,\mathcal{W})$ be semi-coarse spaces, and let 
$\mathcal{V} \times \mathcal{W}$ be the collection of all subsets of finite unions of 
sets of the form $V \boxtimes W$, where $V \in \mathcal{V}$ and $W \in \mathcal{W}$, i.e. 
\begin{align*}
    \mathcal{V} \times \mathcal{W} \coloneqq \{ & U \in \mathcal{P}((X \times Y) 
        \times (X \times Y))\\
        & \mid \exists n \in \N \text{ such that }U \subset \cup_{i=1}^n V_i \boxtimes W_i, \\
        & \text{ where } V_i \in \mathcal{V}, W_i \in \mathcal{W} \;\forall i \in \{1,\dots,n\}\}.
\end{align*}
Then $(X\times Y, \mathcal{V}\times\mathcal{W})$ is a semi-coarse space.
\end{proposicion}

\begin{proof}
Let $(X,\mathcal{V})$ and $(Y,\mathcal{W})$ be semi-coarse spaces and let the collection 
$\mathcal{V}\times\mathcal{W}$ be as in the statement of the proposition. We check that 
$(X \times Y, \mathcal{V} \times \mathcal{W})$ satisfies the axioms for a semi-coarse
space.
\begin{enumerate}[wide,label=(sc\arabic*)]
    \item We observe that $\Delta_X\boxtimes\Delta_Y=\{ ((x,y),(x,y))\in (X\times Y)\times
        (X\times Y)\mid x\in X, y \in Y  \}=\Delta_{X\times Y}$.
\item Let $A\in \mathcal{V}\times\mathcal{W}$ and $B\subset A$. Then there are $n\in\mathbb{N}$,
    $\{V_1,...,V_n\}\subset \mathcal{V}$ and $\{W_1,...,W_n\}\subset\mathcal{W}$ such that
    $B\subset A\subset \bigcup\limits_{1\leq k\leq n} (V_k\boxtimes W_k)$, and therefore $B\in
    \mathcal{V}\times\mathcal{W}$.
\item If $A,B\in\mathcal{V}\times\mathcal{W}$. Then there are $m,n\in\mathbb{N}$, $\{V_1,\dots, V_m\}, 
    \{V'_1, \dots, V'_n\}\subset\mathcal{V}$ and $\{W_1,\dots, W_m\}, \{W'_1\dots,W'_n\}\subset\mathcal{W}$
    such that 
    \begin{align*}
        A\subset & \bigcup\limits_{1\leq k \leq m} (V_k\boxtimes W_k),\\
        B\subset & \bigcup\limits_{1\leq k \leq n} (V'_k\boxtimes W'_k).
    \end{align*}
It follows that $A\cup B \subset 
\left(\bigcup\limits_{1\leq j \leq m} (V_j\boxtimes W_j)\right) \bigcup \left(\bigcup\limits_{1\leq k 
    \leq n} (V'_k\boxtimes W'_k\right)$, which implies that $A\cup B\in \mathcal{V}\times 
    \mathcal{W}$.
\item If $A\in\mathcal{V}\times\mathcal{W}$, then there is an $n\in\mathbb{N}$, $\{V_1,...,V_n\}\subset 
    \mathcal{V}$, and $\{W_1,...,W_n\}\subset\mathcal{W}$ such that $A\subset \bigcup\limits_{1\leq k\leq n} 
    (V_k\boxtimes W_k)$, so $A^{-1}\subset \bigcup\limits_{1\leq k\leq n} (V_k^{-1}\boxtimes W_k^{-1})$ 
    by 
    \autoref{InvUnionInter} and \autoref{boxunion}, and therefore $A^{-1} \in \mathcal{V} \times 
    \mathcal{W}$.
\end{enumerate}
We conclude that $\mathcal{V}\times\mathcal{W}$ is a semi-coarse structure on the set $V\times W$.
\end{proof}

\begin{definicion}[Product Semi-Coarse Space]
\label{DefProdPseCoar}
Let $(X,\mathcal{V})$ and $(Y,\mathcal{W})$ be semi-coarse spaces, and let $\mathcal{V}\times \mathcal{W}$
be the semi-coarse structure on $X \times Y$ constructed in \autoref{prop:Product QC structure}
We call the ordered pair $(X\times Y,\mathcal{V}\times \mathcal{W})$ the 
\textbf{product semi-coarse space of $(X,\mathcal{V})$ and $(Y,\mathcal{W})$}.
\end{definicion}

It is important to observe that this product is associative, which is established in the next
proposition.

\begin{proposicion}
Let $(X,\mathcal{V})$, $(Y,\mathcal{W})$ and $(Z,\mathcal{Z})$ be semi-coarse spaces. Then $(\mathcal{V}\times\mathcal{W})\times\mathcal{Z}=\mathcal{V}\times (\mathcal{W}\times\mathcal{Z})$.
\end{proposicion}

\begin{proof}
Let $(X,\mathcal{V})$, $(Y,\mathcal{W})$ and $(Z,\mathcal{Z})$ be semi-coarse spaces. 
Consider $A\in (\mathcal{V}\times\mathcal{W})\times \mathcal{Z}$. Then there is a natural number 
$n$ and sets $\{D_k\}_{k=1}^n \subset \mathcal{V}\times \mathcal{W}$ and 
$\{E_k\}_{k=1}^n\subset\mathcal{Z}$ such that $A\subset \bigcup\limits_{k=1}^n (D_k\boxtimes E_k)$. 
It follows that, for each $k \in \{1,\dots,n\}$, there are natural numbers 
$\{m_k\}_{k=1}^n\subset \mathbb{N}$ and sets $\{B_{i}\}_{i\in 1}^{m_k}\subset \mathcal{V}$
and $\{C_{i}\}_{i\in 1}^{m_k}\subset\mathcal{W}$ such that $D_k\subset\bigcup\limits_{i=1}^{m_k}
(B_{i}\boxtimes C_{i})$. 
This gives that $A\subset \bigcup\limits_{k=1}^n\bigcup\limits_{i=1}^{m_k} (B_{i}\boxtimes C_{i}\boxtimes E_k)$,
by \autoref{lem:Cross product associative} and \autoref{boxunion}. Since this is a finite union, and,
moreover, since $(C_{i}\boxtimes E_k)\in \mathcal{W}\times\mathcal{Z}$ and $B_{i}\in \mathcal{V}$ for each $i\in\{1,\cdots,m_k\}$, $k\in\{1,\cdots,n\}$, we have that $A\in \mathcal{V}\times (\mathcal{W}\times \mathcal{Z})$. 

The proof that $(\mathcal{V}\times\mathcal{W})\times\mathcal{Z}\supset \mathcal{V}\times 
(\mathcal{W}\times\mathcal{Z})$ is analogous. Thus, the two structures are identical.
\end{proof}

The following alternative characterization of the product structure will be useful for characterizing
bornologous functions.

\begin{proposicion}
\label{LemaProdPseCoar}
Let $(X,\mathcal{V})$ and $(Y,\mathcal{W})$ be semi-coarse spaces. The product structure
$\mathcal{V}\times \mathcal{W}$ is the collection of subsets of sets of the form 
$V\boxtimes W$, where $V\in\mathcal{V}$ and $W\in\mathcal{W}$.
\end{proposicion}

\begin{proof}
Let $(X,\mathcal{V})$ and $(Y,\mathcal{W})$ be semi-coarse spaces, and let $\mathcal{V}\boxtimes \mathcal{W}$ 
denote the collection of subsets of sets of the form $V\boxtimes W$, where $V\in\mathcal{V}$, $W\in\mathcal{W}$.

By definition, it is clear that $\mathcal{V}\boxtimes\mathcal{W}\subset \mathcal{V}\times \mathcal{W}$. 
On the other hand, if $A\in\mathcal{V}\times \mathcal{W}$, then there is a natural number $n \in \N$ and sets 
$\{V_1,\cdots,V_n\}\in\mathcal{V}$ and $\{W_1,\cdots,W_n\}\in\mathcal{W}$ such that 
$A=\bigcup\limits_{k=1}^n (V_k\boxtimes W_k)$. However, we also have that
\begin{align*}
    \bigcup\limits_{k=1}^n V_k\in\mathcal{V}, &\ \bigcup\limits_{k=1}^n W_k \in\mathcal{W}, \text{ and }
    \\ A \subset \bigcup\limits_{k=1}^n (V_k\boxtimes  W_k &)\subset \left( \bigcup\limits_{k=1}^n V_k \right) \boxtimes 
\left( \bigcup\limits_{k=1}^n W_k \right).
\end{align*}
Therefore $A\in \mathcal{V}\boxtimes\mathcal{W}$, and the result follows.
\end{proof}

The following corollaries are now an immediate consequence of the above proposition.

\begin{corolario}
\label{CoroProdPseCoar}
Let $(X,\mathcal{V})$, $(Y,\mathcal{W})$, and $(Z,\mathcal{A})$ be semi-coarse spaces. 
Then $f:(X\times Y,\mathcal{V}\times\mathcal{W})\rightarrow (Z,\mathcal{A})$ is a bornologous function 
iff $f(B\boxtimes C)\in \mathcal{A}$ for each $B\in\mathcal{V}$ and $C\in\mathcal{W}$.
\end{corolario}

\begin{corolario}
\label{cor:Roofed product bornologous}
Let $(X,\mathcal{V}), (Y,\mathcal{W}),$ and $(Z,\mathcal{A})$ be semi-coarse spaces, and suppose
that $(Y,\mathcal{W})$ is roofed (\autoref{DefRoof}). Then $f:(X \times Y,\mathcal{V} \times \mathcal{W}) \to 
(Z,\mathcal{A})$ is a bornologous function iff $(f \times f)(V \times \mathfrak{R}(\mathcal{W}))
\in \mathcal{A}$ for all $V \in \mathcal{V}$.
\end{corolario}

\begin{proof}
    If $(f \times f)(V \times \mathfrak{R}(\mathcal{W})) \in \mathcal{A}$ for all $V \in \mathcal{V}$, 
    then for any $W \in \mathcal{W}$, $W \subset \mathfrak{R}(\mathcal{V})$, so
    $(f \times f)(V \times W) \subset (f \times f)(W \times \mathfrak{R}(\mathcal{V}) \in \mathcal{A}$.
    Therefore, $(f \times f)(V \times W) \in \mathcal{A}$ and $f$ is bornologous by 
    \autoref{LemaProdPseCoar}.

    Conversely, suppose that $f$ is bornologous. Then $(f \times f)(V \times \mathfrak{F}(\mathcal{W}))
    \in \mathcal{A}$.
\end{proof}

\subsection{Quotient Spaces}
\label{SecEspCoc}
In this section, we study quotient spaces of semi-coarse spaces, i.e. the semi-coarse structures constructed
on a set $Y$ given a semi-coarse space $(X,\mathcal{V})$ and a surjective 
map $g:X\rightarrow Y$. We begin with the following definition.

\begin{definicion}
	\label{def:Quotient structure}
	Let $(X,\mathcal{V})$ be a semi-coarse space, let $Y$ be a set, let $g:X\rightarrow Y$ be a surjective function, and define
	\begin{equation*}
		\mathcal{V}_g:=\{(g\times g)(V): V\in\mathcal{V}\}.
	\end{equation*}
	We call $\mathcal{V}_g$ the \emph{semi-coarse quotient structure 
		on $Y$ induced by $g$}, and we call $(Y,\mathcal{V}_g)$ the \emph{semi-coarse quotient (space) of $X$ induced by $g$}.
\end{definicion}

\begin{observacion}
		The set $\mathcal{V}_g$ constructed above will occasionally be referred to as 
		the \emph{semi-coarse structure inductively generated by the function $g:X \to Y$}, and $(Y,\mathcal{V}_g)$ 
		will likewise also be called the \emph{semi-coarse space inductively generated by the function $g:X \to Y$.}
\end{observacion}

The following proposition justifies this choice of terminology.

\begin{proposicion}
\label{prop:QC from function}
Let $(X,\mathcal{V})$, $Y$, $g:X\rightarrow Y$, and $\mathcal{V}_g$ be as in \autoref{def:Quotient structure}.
Then $(Y,\mathcal{V}_g)$ is a semi-coarse space. 
Moreover, $\mathcal{V}_g$ is the coarsest semi-coarse structure which makes the function $g$ bornologous.
\end{proposicion}

\begin{proof}
    Let $(X,\mathcal{V})$, $Y$, $g:X \to Y$, and $\mathcal{V}_g$ be as in the statement of the proposition.
    We verify that $\mathcal{V}_g$ satisfies the axioms of a semi-coarse structure.
    \begin{enumerate}[wide,label=(sc\arabic*)]
\item Since $g$ is surjective, $g(X)=Y$, and therefore $(g\times g)(\Delta_X)=\Delta_Y$, 
    which gives that $\Delta_Y\in \mathcal{V}_g$.
\item Let $B\in\mathcal{V}_g$ and suppose that $A\subset B$. 
    Then there is a set $W\in\mathcal{V}$ such that $(f\times f)(W)=B$, and for each $(y,y')\in A$, 
    there is an $(x,x')\in W$ such that $g(x)=y$ and $g(x')=y'$. Defining
$$A_g := \{ (x,x')\in W : (g\times g)(x,x')\in A \},$$
we have that $A_g\subset W$ and $(g\times g)(A_g)=A$. Since $A_g\in\mathcal{V}$, it follows
that $A\in\mathcal{V}_g$.
\item Let $A,B\in\mathcal{V}_g$. Then there are $A',B'\in\mathcal{V}$ such that 
    $(g\times g)(A')=A$ and $(g\times g)(B')=B$. 
    Now note that $A\cup B = (g\times g)(A')\cup (g\times g)(B')= (g\times g)(A'\cup B')$, so 
    $A\cup B\in \mathcal{V}_g$.
\item Let $A\in\mathcal{V}_g$, then there is $W\in\mathcal{V}$ such that $(g\times g)(A)=W$. We observe that
\begin{align*}
(g\times g)(A^{-1}) = & \{ (x,y) \mid \exists (x',y')\in A^{-1}, (g\times g)(x',y')=(x,y) \}\\
                    = & \{ (x,y) \mid \exists (y',x')\in A, (g\times g)(x',y')=(x,y)\} \\
                    = & \{ (y,x) \mid \exists (x',y') \in A, (g \times g)(x',y') = (x,y) \}\\
                    = & ((g\times g)(A))^{-1}.
\end{align*}
We conclude that $\mathcal{V}_g$ is a semi-coarse structure on $Y$.

We now show that $\mathcal{V}_g$ is the coarsest semi-coarse structure making $g$ bornologous.
By definition, $g$ is a bornologous function iff for each $A\in\mathcal{V}$ we have
that $(g\times g)(A)\in\mathcal{V}_g$. However, by definition,
$\{(g\times g)(A):A\in\mathcal{V}\}=\mathcal{V}_g$, so it follows that any semi-coarse structure $\mathcal{W}$
making $g$ bornologous must contain $\mathcal{V}_g$, and therefore $\mathcal{V}_g$ is the coarsest such structure.
\end{enumerate}
\end{proof}

The following property follows directly from the definition of 
the quotient structure.

\begin{proposicion}
    \label{thm:Bornologous quotient maps}
Let $(Y,\mathcal{V}_g)$ be a semi-coarse space inductively generated by the function $g:X \to Y$
and let $(Z,\mathcal{Z})$ be a semi-coarse space. A function $f:(Y,\mathcal{V}_g) \rightarrow (Z,\mathcal{Z})$
is bornologous iff $f\circ g:X\rightarrow Z$ is bornologous.
\end{proposicion}

\begin{proof}
    Suppose that $f:(Y,\mathcal{V}_g) \to (Z,\mathcal{Z})$ is a bornologous function and let $A \subset X$,
    $A \in\mathcal{V}$. By definition, $(g\times g)(A)\in\mathcal{V}_g$, so $(f\times f)\circ(g\times g)(A)=
    (f\circ g\times f\circ g)(A) \in \mathcal{Z}$.

    Now suppose that $f\circ g:(X,\mathcal{V}) \to (Z,\mathcal{Z})$ is a bornologous function, and let 
    $A\in\mathcal{V}_g$. Then there is a set $W\in\mathcal{V}$ such that $(g\times g)(W)=A$. It follows
    that $(f\times f)(A)=(f\circ g\times f\circ g)(W)\in\mathcal{Z}$, and we conclude
    that $f$ is bornologous.\qedhere
\end{proof}

\begin{definicion}
	Let $p:(X,\mathcal{V}_X)\to (Y,\mathcal{V}_Y)$ be a surjective bornologous map. We say that $p$ is a \emph{quotient map} if,
	in addition, $\mathcal{V}_Y = \mathcal{V}_p$, where $\mathcal{V}_p$ is the semi-coarse structure defined in \autoref{def:Quotient structure}.
\end{definicion}

\begin{observacion}
	\label{rem:Quotient maps}
	It follows immediately from the definition of quotient map that saying that $g: (X,\mathcal{V}_X)\rightarrow (Y,\mathcal{V}_Y)$ is a quotient map is equivalent to the statement that $B\in \mathcal{V}_Y$ if and only if there exists $A\in \mathcal{V}_X$ such that $B=(g\times g)(A)$.
\end{observacion}

The following theorem expresses the universal property of the 
quotient space $(Y,\mathcal{V}_g)$.

\begin{teorema}
	Let $(X,\mathcal{V}_X)$, $(Y,\mathcal{V}_g)$, $g:X\rightarrow Y$ be as in \autoref{def:Quotient structure}. Suppose that $(Z,\mathcal{V}_Z)$ is a semi-coarse space and that there exists a (not necessarily bornologous) map $h:X \to Z$ which is constant on each set $g^{-1}(\{y\})$. Then $h$ induces a map $f:Y \to Z$ such that $f \circ g = h$. The induced map $f$ is bornologous from $(Y,\mathcal{V}_g)\to (Z,\mathcal{V}_Z)$ iff $h$ is bornologous, and $f$ is a quotient map iff $h$ is a quotient map.  
\end{teorema}

\begin{proof}
We first observe that the part of the theorem saying that $f$ is bornologous if and only if $h$ is bornologous is a direct implication of \autoref{thm:Bornologous quotient maps}.

Now we prove the statement about quotient maps. Suppose that $f$ is a quotient map, then by \autoref{rem:Quotient maps}, for every $C\in \mathcal{V}_Z$ there exists $B\in \mathcal{V}_g$ such that $C=(f\times f)(B)$. At the same time, $g$ is a quotient map by hypothesis, so there exists $A\in \mathcal{V}_X$ such that $B=(g\times g)(A)$. Thus, there exists $A\in \mathcal{V}_X$ such that $C=(f\times f)\circ(g\times g)(A) = (h\times h)(A)$, obtaining that $h$ is a quotient map. Observe that this fragment of the proof also says that the composition of quotient maps is a quotient map.

On the other hand, let $h$ be a quotient map. Then for every $C\in\mathcal{V}_Z$ there exists $A\in \mathcal{V}_X$ such that $C=(h\times h)(A)$. Since $g$ is a bornologous map, then $(g\times g)(A)\in \mathcal{V}_g$. Since $C = (h\times h)(A) = (f\times f)\circ (g\times g)(A)$, we may conclude that for every $C\in \mathcal{V}_Z$ there exists $B=(g\times g)(A)\in \mathcal{V}_g$ such that $C=(f\times f)(B)$. Therefore, $f$ is a quotient map, as desired.
\end{proof}

An important example of an inductively generated semi-coarse space is the quotient space generated by 
an equivalence relation.

\begin{ejemplo}[Quotient Semi-Coarse Space]
\label{def:Quotient qc space}
Let $(X,\mathcal{V})$ be a semi-coarse space and let ${\sim}$ be an equivalence relation on $X$. Let 
$p:X \to X/{\sim}$ be the map $x \mapsto [x]$ sending each point $x \in X$ to its equivalence class 
$[x] \in X/{\sim}$. Then the semi-coarse space $(X/{\sim},\mathcal{V}_p)$ is called the 
\emph{quotient space of $X$ induced by the equivalence relation ${\sim}$}.
\end{ejemplo}

It will be useful in the following to have an alternate formulation of the quotient semi-coarse structure
induced by an equivalence relation, which we provide in the next proposition.

\begin{proposicion}
    Let $(X,\mathcal{V})$ be a semi-coarse space and let ${\sim}$ be an equivalence relation on $X$. 
    Furthermore, extend the relation ${\sim}$ to $X \times X$ by defining $(x,y) \sim (x',y')$
    iff $x\sim x'$ and $y\sim y'$. We denote by $[x,y]$ the equivalence class of $(x,y)$, and
    for a subset $A\subset X\times X$, we define
    \begin{equation*}
        [A]:=\{ [x,y]\mid (x,y)\in A\},
    \end{equation*}
    and we let $\mathcal{V}/{\sim}$ denote the collection
    \begin{equation*}
        \mathcal{V}/{\sim} := \{ [B] \mid B\in\mathcal{V} \}.
    \end{equation*}
    Finally, let $p:X \to X/{\sim}$ be the map $x \mapsto [x]$ sending a point $x \in X$
    to its equivalence class $[x] \in X/{\sim}$.
    
    Then $\mathcal{V}_p = \mathcal{V}/{\sim}$.
\end{proposicion}

\begin{proof}
    First, let $A\in\mathcal{V}/{\sim}$. Then there is a set $A'\in\mathcal{V}$ such that $A=[A']$. Furthermore, 
    $A=[A']=(p\times p)(A')$, and we conclude that $A\in \mathcal{V}_p$, so $\mathcal{V}/{\sim} \subset 
    \mathcal{V}_p$.

    Now suppose that $A\in\mathcal{V}_p$. Then there is a set $A'\in\mathcal{V}$ such that $(p\times p)(A')=A$. 
    However, by the definition of $p$, $(p\times p)(A')=[A']=A$, and therefore $A\in \mathcal{V}/\sim$, 
    so $\mathcal{V}_p \subset \mathcal{V}/{\sim}$ as well.\qedhere
\end{proof}

\section{Homotopy}
\label{sec:Homotopy}

In this section, we develop the basics of homotopy theory in semi-coarse spaces. One 
of the key difficulties in doing so is that there is no natural semi-coarse structure on 
the topological interval $[0,1]$. Naively, one could try to endow the set $[0,1]$ with any of the semi-coarse structures $\mathcal{V}_r$,
$\mathcal{V}^<_r$, $0<r<1$ defined in \ref{ex:Qc from metric}, and then
use the product $(X,\mathcal{V}) \times ([0,1],\mathcal{V}_r)$ as the cylinder for 
the space $(X,\mathcal{V})$. However, in the following proposition,
we see that the concatenation of two semi-coarse intervals
$([0,1],\mathcal{V}_{r'})$ and $([0,1],\mathcal{V}_{r''})$ does not result 
in an interval of the form $([0,1],\mathcal{V}_r)$.
A consequence of this is that proceeding in this manner will not allow us to define
a product between two paths by simply concatenating them, complicating
the development of homotopy theory for semi-coarse spaces.

\begin{proposicion}
\label{prop:ClosedIntervalinSC}
	Let $r,r',r''\in (0,1)$. Any bijective map 
	\[\phi:([0,1],\mathcal{V}_{r}) \to ([0,1],\mathcal{V}_{r'})\sqcup_{1\sim 1}([1,2],\mathcal{V}_{r''}),\]
	 where the latter space is the pushout of the two intervals, glued at one endpoint, is not bornologous at the point $1\in [0,2]$ where the two intervals are glued. Then $\phi$ cannot be a bornologous map, and the spaces $([0,1],\mathcal{V}_{r})$ and $([0,1],\mathcal{V}_{r'})\sqcup_{1\sim 1}([1,2],\mathcal{V}_{r''})$ are not homeomorphic.
\end{proposicion}

\begin{proof}
	Denote $([0,1],\mathcal{V}_{r'}) \sqcup_{1\sim 1}([1,2],\mathcal{V}_{r''})$ by $(Z,\mathcal{V}_Z)$.
	Suppose that $\phi:([0,1],\mathcal{V}_r) \to (Z,\mathcal{V}_Z)$ is
	bijective and  bornologous.
	By the Archimedean property, there exists an integer $n$ such that $(n-1)(r/2)< 1 \leq n(r/2)$; in particular, we can compute $n$ as the least integer that is greater or equal to $2/r$.
	
	Consider the set $\{x_0, x_1, \ldots, x_{n-1}, x_n\}$ with $x_i= i(r/2)$ if $i<n$ and $x_n=1$.
	Observe that $|x_i-x_{i+1}|<r$ for every $i$, then \[(\phi \times \phi)([x_i,x_{i+1}] \times [x_i,x_{i+1}]) \in \mathcal{V}_r\]
	There is no pair $z_i,z'_i\in [x_i,x_{i+1}]$ such that $\phi(z_i)<1$ and $\phi(z'_i)>1$ because $\{ (\phi(z_i), \phi(z'_i)) \}\notin \mathcal{V}_Z$ and $\phi$ would not be a bornologus map.
	Thus, we obtain that either $\phi([x_i,x_{i+1}])\subset [0,1]$ or $\phi([x_i,x_{i+1}])\subset [1,2]$.
	Without loss of generality, suppose that $\phi([0,x_1])\subset [1,2]$; this implies that every $\phi([x_i,x_{i+1}]) \subset [1,2]$ and then $\phi([0,1])\subset [1,2]$, from which we conclude that $\phi$ is not bijective, a contradiction.
	
	Therefore $\phi$ is either not bornologous or not bijective.
\end{proof}

 In order to resolve the problem highlighted by the \autoref{prop:ClosedIntervalinSC}, we adapt the construction
in \cite{Babson_etal_2006} to the semi-coarse category, using finite intervals in $\Z$ endowed
with their `nearest neighbor'
semi-coarse structure in place of the interval for homotopical constructions. In contrast
to \cite{Babson_etal_2006}, however, the product which we use throughout is the categorical product
for semi-coarse spaces, which gives the resulting theory a more simplicial, rather than
cubical, character.

\subsection{Homotopy and Homotopy Equivalence for Semi-Coarse Spaces}

Before we define homotopy, we define the semi-coarse structure on $\Z$ which
will be used throughout.

\begin{definicion}
\label{CanonicalIntQC}
We call the semi-coarse structure on $\Z$ generated by the graph $G = (\Z,E)$, where 
$E \coloneqq \bigcup_{z\in\mathbb{Z}} \{(z,z-1),(z,z),(z,z+1)\}$, i.e.
 
\begin{center}
    \vspace{5pt}
\begin{tikzpicture}
	\node (0) at (-1,0) {};
    \filldraw [black] (0,0) circle (2 pt);
    \node (1) [label=south:-2] at (0,0) {}; 
    \filldraw [black] (1,0) circle (2 pt);
    \node (2) [label=south:-1] at (1,0) {};
    \filldraw [black] (2,0) circle (2 pt);
    \node (3) [label=south:0] at (2,0) {};
    \filldraw [black] (3,0) circle (2 pt);
    \node (4) [label=south:1] at (3,0) {};
    \filldraw [black] (4,0) circle (2 pt);
    \node (5) [label=south:2] at (4,0) {};
    \node (6) at (5,0) {};
    \draw (1) -- (2) -- (3) -- (4) -- (5);
    \draw[dashed] (0) -- (1);
    \draw[dashed] (5) -- (6);
\end{tikzpicture}
\end{center}
the \emph{the canonical semi-coarse 
structure of $\mathbb{Z}$}, and we denote this structure by $\mathcal{Z}$.
\end{definicion}
\begin{observacion}
    Note that the canonical semi-coarse structure on $\Z$ is the roofed semi-coarse
    structure (\autoref{DefRoof}) on $\Z$ with roof $E$.
\end{observacion}

With the semi-coarse structure on $\Z$ in place, we now define homotopy between bornologous maps. 

\begin{definicion}
\label{HomotopiaPseCoar}
Let $(X,\mathcal{V})$ and $(Y,\mathcal{W})$ be semi-coarse spaces and let $f,g:(X,\mathcal{V})
\rightarrow (Y,\mathcal{W})$ be bornologous functions. We will say that \emph{f is homotopic
to g}, written $f\simeq_{sc} g$, iff there is a bornologous function $H:(X\times\mathbb{Z},
\mathcal{V}\times\mathcal{Z})\rightarrow (Y,\mathcal{W})$ and integers $N,M\in\mathbb{Z}$ with
$M<N$, where $H(x,n)=f(x)$ if $n\leq M$ and $H(x,n)=g(x)$ if $n\geq N$.

We will say that two semi-coarse spaces $(X,\mathcal{V})$ and $(Y,\mathcal{W})$
are \emph{homotopy equivalent} when there are bornologous functions $f:X\rightarrow Y$ and 
$g:Y\rightarrow X$ such that $g\circ f\simeq_{sc}id_X$ and $f\circ g\simeq_{sc}id_Y$.
\end{definicion}

\begin{proposicion}
Homotopy equivalence $\simeq_{sc}$ is an equivalence relation on bornologous functions.
\end{proposicion}

\begin{proof}
We check that $\simeq_{sc}$ satisfies the axioms of an equivalence relation.
\begin{itemize}[wide]
    \item[\emph{Reflexivity:}] Let $(X,\mathcal{V})$ and $(Y,\mathcal{W})$ be semi-coarse spaces, 
    let $f:(X,\mathcal{V}) \rightarrow (Y,\mathcal{W})$ be a bornologous function, and suppose that 
$H:X\times\Z\rightarrow Y$ is a map such that $H(x,z)=f(x)$ for each $z\in\mathbb{Z}$.
Since $H(x,n) = x$ for all $n \in \Z$, 
to see that $f \simeq_{sc} f$, it is enough to prove that $H$ is bornologous.
From the definition of $H$, we have
\begin{align*}
    (H\times H)(V\boxtimes \mathfrak{R}(\mathcal{Z})) = (f\times f)(V) \in\mathcal{W},
\end{align*}
so that $H$ is bornologous by \autoref{cor:Roofed product bornologous}. It follows that
$f\simeq_{sc}f$.
\item[\emph{Symmetry:}] Suppose that $f\simeq_{sc}g$. Then there exists a bornologous function 
$H:X\times\mathbb{Z}\rightarrow Y$ and integers $M<N$ such that $H(x,n) = f(x)$ for all $n \leq M$ and
$H(x,x) = g(x)$ for all $n \geq N$. Defining $H'(x,n)=H(x,-n)$, we have that $H'(x,n)=g(x)$ for 
each $x\in X$ if $n\leq -N$, and $H'(x,n)=f(x)$ for each $x\in X$ is $n'\geq -M$. 
Moreover, the function $h:(X \times \Z,\mathcal{V} \times \mathcal{Z})\to (X \times \Z,\mathcal{V}
\times \mathcal{Z})$ given by $h(x,n) \coloneqq (x,-n)$ is bornologous by 
\autoref{cor:Roofed product bornologous}, since
$(h \times h)(V \times \mathfrak{R}) = V \times \mathfrak{R} \in \mathcal{V} \times \mathcal{Z}$.
Therefore, $H' = H \circ h$ is a bornologous function, from which it follows that $g\simeq_{sc}f$.

\item[\emph{Transitivity:}] Suppose that $f\simeq_{sc}g$ and $g\simeq_{sc}h$. 
    Then there are bornologous functions $H_{fg}:(X\times\mathbb{Z},\mathcal{V} \times \mathcal{Z}
    \rightarrow Y$ and $H_{gh}:X\times\mathbb{Z}\rightarrow Y$ and pairs of integers 
    $M_f<N_g$ and $M_g<N_h$
    such that
    \begin{align*}
        H_{fg}(x,n)=& \begin{cases}
            f(x) &  n \leq M_f\\
            g(x) & n \geq N_g,
        \end{cases} \text{ and }\\
        H_{gh}(x,n) = &\begin{cases}
            g(x) & n \leq M_g\\
            h(x) & n \leq N_h.
        \end{cases}
    \end{align*}
    Without loss of generality, we take $M_f<-1=N_g$ and $M_g=1<N_h$.

Define $H_{fh}:X\times\mathbb{Z}\rightarrow Y$ by
\begin{equation*}
    H_{fh}(x,n) = \begin{cases}
        H_{fg}(x,n) & n \leq 0\\
        H_{gh}(x,n) & n> 0,
    \end{cases}
\end{equation*}
and observe that $H_{fh}(x,n) = f(x)$ for all $n\leq M_f$, $H_{fh}(x,n) = h(x)$ for all $n \geq N_h$, 
and $H_{fh}(x,n)= g(x) = H_{fg}(x,n) = H_{gh}(x,n)$ for $ -1 = N_g \leq n \leq M_g= 1$. 
$H_{fh}$ is bornologous by \autoref{prop:Bornologous function from subspaces}, where 
$X_1 = X \times ((-\infty,0] \times \Z)$ and $X_2 = X \times ([0,\infty) \cap \Z)$, 
and $X_1$ and $X_2$ are endowed with the respective subspace semi-coarse structures induced from $\mathcal{V}
\times \mathcal{Z}$. 

Let $\Z^2_- \coloneqq \{ (a,b)\in\mathbb{Z}^2 \mid a \leq 0 \}$,
and observe that, for any $V \in \mathcal{V}$,
if $(x,n,x',n')\in V \boxtimes (\mathfrak{R}(\mathcal{Z})\cap \Z^2_-)$ 
then $n'=n+k$, where $k\in\{-1,0,1\}$. In this case, we also have
\begin{align*}
H_{fh}(V \boxtimes (\mathfrak{R}(\mathcal{Z})\cap \Z_-^2) \}))=
H_{fg}(V \boxtimes (\mathfrak{R}(\mathcal{Z})\cap \Z_-^2) \},
\end{align*}
and therefore $H_{fh}(V \boxtimes (\mathfrak{R}(\mathcal{Z}) \cap \mathbb{Z}^2_-))\in\mathcal{W}$.

Similarly, let $\Z^2_+ \coloneqq \{(a,b) \in \Z^2 \mid a \geq 0\}$, and observe that, for any 
$V \in \mathcal{V})$,
if $(x,z,x',z')\in V \boxtimes (\mathfrak{R}(\mathcal{Z})\cap \mathbb{Z}^2_+)$, then $z'=z+k$, 
where $k\in\{-1,0,1\}$. In this case, we also have
\begin{align*}
    H_{fh}(V \boxtimes (\mathfrak{R}(\mathcal{Z})\cap \mathbb{Z}^2_+))=H_{gh}(V \boxtimes 
(\mathfrak(\mathcal{Z})\cap \mathbb{Z}^2_+)),
\end{align*}
and therefore $H_{fh}(V \boxtimes (\mathfrak{R}(\mathcal{Z})\cap \mathbb{Z}^2_+))\in\mathcal{W}$. 
Since
\begin{align*}
H_{fh}(V\boxtimes \mathfrak{R}(\mathbb{Z})) = H_{fh}(V \boxtimes (\mathfrak{R}(\mathcal{Z})\cap \mathbb{Z}^2_-))\cup
    H_{fh}(V \boxtimes (\mathfrak{R}(\mathcal{Z})\cap \mathbb{Z}^2_+))
\end{align*}
for each $V\in\mathcal{V}$, it follows from \autoref{cor:Roofed product bornologous} that $H_{fh}$ 
is a bornologous function. 
We conclude that $f\simeq_{sc} h$.\qedhere
\end{itemize}
\end{proof}
The following proposition gives a useful condition for determining whether a potential homotopy is bornologous.

\begin{proposicion}
    Let $(X,\mathcal{V})$ and $(Y,\mathcal{W})$ be semi-coarse spaces, let
    $f,g:(X,\mathcal{V}) \to (Y,\mathcal{W})$ be bornologous maps, and suppose that $H:X \times \Z \to Y$ 
    is a (possibly non-bornologous) map such that there exist integers
    $M_f < N_g$ with $H(x,n) = f(x)$ for all $n \leq M_f$
    and $H(x,n) = g(x)$ for all $n \geq N_g$. 

    If $(H\times H)(V\boxtimes (n,n+1))$, 
    $(H\times H)(V \boxtimes (n,n))$,
    and $(H\times H)(V\boxtimes (n,n-1))$ are controlled sets of $(Y,\mathcal{W})$
    for each $V\in \mathcal{V}$ and $n\in\mathbb{Z}$, then $H$ is bornologous.
\end{proposicion}

\begin{proof}
    Let $\Z^2_{fg} \coloneqq \mathfrak{R}(\mathcal{Z}) \cap 
        \{(a,b) \in \Z^2 \mid M_f \leq a \leq N_g\}$.
    Since $H(x,n) = f(x)$ for $n \leq M_f$ and $H(x,n) = g(x)$ for $n \geq N_g$, it follows that,
    for any $V \in \mathcal{V}$,
    \begin{align*}
        (H\times H)(V \boxtimes \mathfrak{R}(\mathcal{Z}))= & 
        (H\times H)(V \boxtimes \Z^2_{fg})\\
        = & (H\times H)\left(\bigcup_{(a,b) \in \Z^2_{fg}}(V \boxtimes (a,b))\right).
    \end{align*}
    Since the finite union of sets in $\mathcal{W}$ are in $\mathcal{W}$, $H$ is bornologous
    iff, for all $(a,b) \in \Z^2_{fg}$, $(H \times H)(V \boxtimes (a,b)) \in \mathcal{W}$, which
    is true iff $(H\times H)(V\boxtimes (n,n+1))$, 
    $(H\times H)(V \boxtimes (n,n))$,
    and $(H\times H)(V\boxtimes (n,n-1))$ are controlled sets of $(Y,\mathcal{W})$ for every
    $n \in \Z$.
\end{proof}

In the following, we define the homotopy of maps of pairs and triples for semi-coarse spaces and show that 
they are also an equivalence relation.

\begin{definicion}[Homotopy of Maps of Pairs and Triples]
\label{HomotopiaRelPseCoar}
Let $(X,\mathcal{V})$ and $(Y,\mathcal{W})$ be semi-coarse spaces, let $B\subset A\subset X$ and 
$D \subset C \subset Y$ be endowed with 
the subspace semi-coarse structures, and let $f,g:(X,A,B;\mathcal{V})\rightarrow (Y,C,D;\mathcal{W})$ be bornologous 
functions of a triple, i.e. such that $f|_A \subset C$ and $f|_B \subset D$.
We say that \emph{relatively homotopic to $g$} and write $f\simeq_{sc} g$ iff there is a homotopy
$H:(X\times\mathbb{Z},\mathcal{V}\times \mathcal{Z})\rightarrow (Y,\mathcal{W})$ such that $H|_{A\times \Z} 
\subset C$ and $H|_{B\times \Z} \subset D$.

We define a homotopy between maps of pairs $f,g:(X,A;\mathcal{V}) \to (Y,C;\mathcal{W})$ 
to be the homotopy between maps of a triple as above with $B =A$ and $C = D$.
\end{definicion}

\begin{lema}
Let $(X,\mathcal{V})$ and $(Y,\mathcal{W})$ be semi-coarse spaces, 
and $B\subset A\subset X$, $D\subset C \subset Y$ be endowed with the subspace semi-coarse structures. 
Then the relative homotopy of triples $\simeq_{sc}$ is an equivalence relation.
\end{lema}

\begin{proof}
The proof is nearly identical to the proof of \autoref{HomotopiaPseCoar}, with the addition that
we need to observe that the $H_{fh}|_A\subset C$ (and $H_{fh}|_D\subset B$). However, this follows the fact 
that $H_{fg}$ and $H_{gh}$ are homotopies of triples and from the definition of $H_{fh}$.
\end{proof}

We end this section with the definition of homotopy equivalence.

\begin{definicion}
    Two semi-coarse spaces $(X,\mathcal{V})$ and $(Y,\mathcal{W})$ are said to be 
    \emph{homotopy equivalent}, written $(X,\mathcal{V})\simeq (Y,\mathcal{W})$, iff
    there exist bornologous maps $f:(X,\mathcal{V})\to (Y,\mathcal{W})$ and 
    $g:(Y,\mathcal{W}) \to (X,\mathcal{V})$ such that $f\circ g \simeq_{sc} 1_Y$
    and $g\circ f \simeq_{sc} 1_X$.
\end{definicion}

\subsection{Homotopy Groups}

As mentioned in Section \ref{cap.introduccion}, to define the homotopy groups for semi-coarse spaces,
we will largely follow the ideas in the construction of the $A$-homotopy groups 
from \cite{Babson_etal_2006}, with the 
difference that the product we use is the categorical product rather than a semi-coarse version
of the graph product. We begin this section with the construction of the relative homotopy classes
of maps of cubes in $\Z^n$ to a semi-coarse space $(X,\mathcal{V})$. We set some notation in the
next definition.

\begin{definicion}
\label{def:nmCubos}
Let $n$ and $m$ be natural numbers, let $I_m$ be the set $\{0,1,\cdots,m\}$, and let $\mathcal{Z}_m$
be the subspace semi-coarse structure induced by the inclusion of $I_m$ in 
$(\Z,\mathcal{Z})$. We denote by $(I_m,\mathcal{Z}_m)$ 
the corresponding semi-coarse space, and we denote by $(I_m^n,\mathcal{Z}^n_m)$ the Cartesian product of 
$n$ copies of $(I_m,\mathcal{Z}_m)$. We write $a_i$ for the $i$-th coordinate 
of a point $a\in \Z^n$, where $i\in\{1,\cdots,n\}$, unless otherwise stated. Finally, we 
denote by \[
\partial I_m^n \coloneqq \{(x_1,\ldots,x_n)\in I_m^n\mid \exists i\in\{1,\ldots,n\} \text{ such that } x_i\in\{0,m\}  \},\]
the boundary of the discrete $n$-cube with sides of length $m$, 
and we define $J^{n-1}_m \subset \partial I_m^n$ to be the ``discrete (upside-down) open $n$-dimensional box with sides of length $m$", i.e. the set 
\begin{equation*}
    J^{n-1}_m = \{ (x_1,\dots,x_n) \in \partial I^n_m \mid \text{If } x_n = 0,\text{ then }\exists i\neq n \text{ such that } x_i\in\{0,m\}\}.
\end{equation*}
\end{definicion}

To define the relative homotopy groups of a pointed semi-coarse space pair $(X,A,*;\mathcal{V})$ (i.e.
$* \in A \subset X$),
we will consider maps of triples
\[f:(I_m^n,\partial^n_m,J^{n-1}_m;\mathcal{Z}_m^n) \to (X,A,*;\mathcal{V}),\] and
we define $\pi_n^{sc}(X,A,*;\mathcal{V})$ to be a certain direct limit which we construct below.
First, however, we require the following definition and lemmas.

\begin{definicion}
    Let $m < m'$ be two positive integers. We define the retraction $\phi^{m'}_m:I^n_{m'} \to I^n_m$ by
    \begin{equation*}
        (\phi^{m'}_m(x))_k \coloneqq \begin{cases}
            x_k & x_k \leq m,\\
            m & m < x_k \leq m',
        \end{cases}
    \end{equation*}
    where $(\phi^{m'}_m(x))_k$ is the $k$-th coordinate of $\phi^{m'}_m(x)$.
    We will typically abuse notation and write $\phi^{m'}_m$ for the retraction $I^n_{m'} \to I^n_m$ for 
    any dimension $n$, as the dimension will usually be clear from context.
\end{definicion}

\begin{lema}
    The function $\phi^{m'}_m:(I^n_{m'},\mathcal{Z}^n_{m'})\to (I^n_m,\mathcal{Z}^n_m)$ is bornologous.
\end{lema}
\begin{proof} The result follows from \autoref{prop:Bornologous function from subspaces}, 
    \autoref{prop:Born on roof is born}, and the definition of the semi-coarse structures $\mathcal{Z}^n_{m'}$
    and $\mathcal{Z}^n_m$.
\end{proof}

\begin{lema}
    \label{lem:Extension}Let $f:(I^n_m,\partial I^n_m,J^{n-1}_m;\mathcal{Z}^n) \to (X,A,*;\mathcal{V})$ 
    be a bornologous function of triples. Then for any $m'\in \Z$ with $m'>m$, $f$ has a bornologous 
    extension of triples
    \begin{equation*}
        f':(I^{n}_{m'}, \partial I^{n}_{m'},J^{n-1}_{m'};\mathcal{Z}^n) \to (X,A,*;\mathcal{V})
    \end{equation*}
    given by $f'(z) = f\circ \phi^{m'}_m(z)$.

    Furthermore, if $f \simeq_{sc} g$ 
    as maps of triples $(I^n_m,\partial I^n_m,J^{n-1}_m;\mathcal{Z}_m^n) \to (X,A,*;\mathcal{V})$, 
    then for the extensions $f'$ and $g'$ to $I^n_{m'}$ of $f$ and $g$, respectively, we have 
    $f' \simeq_{sc} g'$ as well.
\end{lema}

\begin{proof} The extension $f'$ of $f$ is bornologous since the composition of bornologous functions is
    bornologous. 

    If $f \simeq_{sc} g$ then there exists a homotopy of triples
    $H:(I^n_m \times \Z,\mathcal{Z}^n_m \times \mathcal{Z})\to (X,\mathcal{V})$ between $f$ and $g$. 
    We define $H':(I^n_{m'} \times \Z,\mathcal{Z}^n_m \times\mathcal{Z}) \to (X,\mathcal{V})$ by
    \begin{equation*}
        H'(x,n) \coloneqq H(\phi^{m'}_m(x),n).
    \end{equation*}
    Since $\phi^{m'}_m \times Id$ and $H$ are bornologous, it follows that $H'$ is bornologous, 
    and therefore $H'$ is a homotopy from $f'$ to $g'$ by definition.
\end{proof}

\autoref{lem:Extension} implies that extending a map from $I_m^n$ to $I_{m'}^n$ induces 
a map on the relative homotopy classes 
\[ i^{m'}_m:[(I^n_m,\partial I^n_m,J^{n-1}_m;\mathcal{Z}_m^n),
(X,A,*;\mathcal{V})] \to 
[(I^n_{m'},\partial I^n_{m'},J^{n-1}_m;\mathcal{Z}_m^n),(X,A,*;\mathcal{V})]\]
such that $i_m^{m''} = i_{m'}^{m''}
\circ i_m^{m'}$. The maps $i_m^{m'}$ therefore make the homotopy
classes into a directed system of sets $\left([(I^n_m,\partial I^n_m,\mathcal{Z}^n),
(X,*,\mathcal{V})],i^{m'}_m,\N\right)$. 
We define the classes $\pi_n^{sc}(X,A*;\mathcal{V})$ in the following manner. 

\begin{definicion}
    Given the set $\{*,1\}$ with the diagonal semi-coarse structure $\mathcal{D}_{*,1}$, i.e. $\mathcal{D}_{*,1}=\mathcal{P}\left(\{(*,*),(1,1)\}\right)$.
    We define the set $\pi_0^{sc}(X,A,*;\mathcal{V})$ by
    \begin{equation*}
        \pi_0^{sc}(X,A,*;\mathcal{V}) \coloneqq [(\{*,1\},*,*;\mathcal{D}_{0,1}),(X,A,*;\mathcal{V})].
    \end{equation*}
    For $n \geq 1$, we define $\pi_n^{sc}(X,A,*;\mathcal{V})$ by
    \begin{equation*}
        \pi_n^{sc}(X,A,*;\mathcal{V}) \coloneqq \varinjlim 
        \left([(I^n_m,\partial I^n_m,J^{n-1}_m;\mathcal{Z}_m^n),(X,A,*;\mathcal{V})],i^{m'}_m,\N\right),
    \end{equation*}
    For a homotopy class $[f] \in [(I_m^n,\partial I_m^n,J^{n-1}_m;\mathcal{Z}_m^n),(X,A,*;\mathcal{V})]$, we write
$\langle f \rangle$ for its image in $\pi_n^{sc}(X,*,\mathcal{V})$.
    
When $A = *$, we write $\pi_n^{sc}(X,*,*;\mathcal{V})$ as $\pi_n^{sc}(X,*;\mathcal{V})$ or $\pi_n^{sc}(X;\mathcal{V})$.
\end{definicion}

\begin{observacion}
    One may also construct the homotopy classes $\pi_n^{sc}(X,*,\mathcal{V})$ as a 
    direct limit of homotopy classes of maps of pairs 
    $f:(\Z^n,J_m^n,\mathcal{Z}^n) \to (X,*,\mathcal{V})$,
    where $J_m^n \coloneqq \Z - I_{m-1}^n$. One sees that these formulations are equivalent
    by extending the maps $f:(I_m^n,\partial I_m^n,\mathcal{Z}_m^n) \to (X,*,\mathcal{V})$
    to $\Z^n$, where the extension $f'$ sends all of $J_m^n$ to $* \in X$. The same
    argument as in \autoref{lem:Extension} shows that $f'$ is bornologous, and it follows from
    the definitions that
    \begin{equation*}
        f \simeq_{sc} g  \iff
        f'\simeq_{sc} g'.
    \end{equation*}
    Therefore the homotopy classes $\pi_n^{sc}(X,*;\mathcal{V})$ may be defined as
    \begin{equation*}
        \pi_n^{sc}(X,*;\mathcal{V}) \coloneqq \varinjlim \left( [(\Z^n,E^n_m,\mathcal{Z}^n),
        (X,*,\mathcal{V})],i_m^{m'},\N\right),
    \end{equation*}
    where, as before, the $i_{m}^{m'}$ are the maps on the respective homotopy classes induced
    by the interpreting a map $f:(\Z^n,E_m^n,\mathcal{Z}^n) \to (X,*,\mathcal{V})$ as a map
    $f:(\Z^n,E_{m'}^n,\mathcal{Z}^n) \to (X,*,\mathcal{V})$.
\end{observacion}

\color{black}
\begin{observacion}
    \label{rem:Equivalence of classes}
    Let $\langle f \rangle,\langle g \rangle\in \pi_n^{sc}(X,A,*;\mathcal{V})$. 
    Observe that $\langle f\rangle=\langle g \rangle$ 
    iff there exist $m\in\mathbb{N}$ and bornologous functions $\tilde{f},\tilde{g}:
    (I_m^n,\mathcal{Z}^n_m) \rightarrow (X,\mathcal{V})$ such that $[\tilde{f}]\in \langle f \rangle$, 
    $[\tilde{g}]\in\langle g\rangle$ and $\tilde{f}\simeq_{sc} \tilde{g}$.
\end{observacion}

We will now define an operation $\star$ on the homotopy classes $\pi_n^{sc}(X,*;\mathcal{V})$ and
$\pi(X,A,*;\mathcal{V})$
and show this operation makes these classes into a group $n\geq 1$ and $n\geq 2$, respectively.

\begin{definicion}[Operation $\star$]
\label{def:OpStar}
Let $X$ be a set, let $n\in\mathbb{N}$, $m$ and $m'$ be non-negative integers, and suppose that
$f:I^n_m \rightarrow (X,\mathcal{V})$ and $g:I_{m'}^n\rightarrow X$
are bornologous functions. Let $\alpha = (\alpha_1,\alpha_2,\dots,\alpha_n) \in I^n_{m+m'}$, and let $e_1 = (1,0,\dots,0) \in \Z^n$. We define the $\star$-product $f\star g: I_{m+m'}^n\rightarrow X$ such that
\begin{align*}
    f\star g(\alpha) = \begin{cases}
        f(\phi_m^{m+m'}(\alpha))   & \mbox{ if }  0\leq \alpha_1 \leq m\\
        g(\phi_{m'}^{m+m'}(\alpha-m e_1)) & \mbox{ if }  m<\alpha_1 \leq m+m'.
    \end{cases}
\end{align*}
\end{definicion}
The following proposition tells us that the result of applying the $\star$ operation to a pair of
bornologous functions which agree on the boundary is bornologous. 
\begin{proposicion}
\label{BornologousStar}
Let $n\in\mathbb{N}$, let $m$ and $m'$ be non-negative integers, and suppose that
$f:(I^n_m,\mathcal{Z}^n_m)\rightarrow (X,\mathcal{V})$ and $g:(I_{m'}^n,\mathcal{Z}^n_m)\rightarrow (X,\mathcal{V})$
are bornologous functions such that 
\begin{equation*}
    f((m,\alpha_2,\dots,\alpha_n)) = g((0,\alpha_2,\dots,\alpha_{n})) 
\end{equation*}
is satisfied for all $\alpha_2,\dots,\alpha_n$ where both $f((m,\alpha_2,\dots,\alpha_n))$ and 
$g((0,\alpha_2,\dots,\alpha_n))$ are defined. Then $f\star g$ is bornologous.
\end{proposicion}
\begin{proof}
    Let $n\in\mathbb{N}$, $m$ and $m'$ be non-negative integers, and suppose that $f:(I^n_m,\mathcal{Z}^n_m)
    \rightarrow (X,\mathcal{V})$ and $g:(I_{m'}^n,\mathcal{Z}^n_{m'})\rightarrow (X,\mathcal{V})$ satisfy the
    hypothesis of the proposition. Let $K_1$ and $K_2$ be the sets 
    \begin{align*}
        K_1 =& \{\alpha \in I^n_{m+m'} \mid 0 \leq \alpha_1 \leq m\\
        K_2 =& \{\alpha \in I^n_{m+m'} \mid m \leq \alpha_1 \leq m'+m.
    \end{align*}
    Then the restrictions $f\star g|_{K_1}$ and $f\star g|_{K_2}$ are bornologous, and $K_1$ and $K_2$ satisfy the 
    hypothesis of \autoref{prop:Bornologous function from subspaces}. Therefore $f \star g$ is 
    bornologous by \autoref{prop:Bornologous function from subspaces}.
\end{proof}

\begin{corolario}
    \label{cor:Star of triples born}
    Suppose that $f:(I^n_m,\partial I^n_m,J^{n-1}_m;\mathcal{Z}^n_m) \to (X,A,*;\mathcal{V})$ and 
    $g:(I^n_{m'},\partial I^n_{m'},J^{n-1}_{m'};\mathcal{Z}^n_{m'}) \to (X,A,*;\mathcal{V})$ are bornologous
    functions of triples. Then $f\star g:(I^n_{m+m'},\partial I^n_{m+m'}, J^{n-1}_{m+m'};\mathcal{Z}^n_{m+m'})
    \to (X,\mathcal{V})$ is bornologous. 
\end{corolario}

\begin{proof} Since $f$ and $g$ are maps of triples, it follows by definition that
    \[f((m,\alpha_2,\dots,\alpha_n)) = * = g(-m',\alpha_2,\dots,\alpha_n))\]
    whenever both sides are defined, so the hypotheses of \autoref{BornologousStar} are satisfied. The conclusion follows.
\end{proof}

\begin{corolario}
    Let $f:(I^n_m,\partial I^n_m;\mathcal{Z}^n_m) \to (X,*;\mathcal{V})$ and 
    $g:(I^n_{m'},\partial I^n_{m'};\mathcal{Z}^n_{m'}) \to (X,*;\mathcal{V})$ be bornologous
    functions of triples. Then $f\star g:(I^n_{m+m'},\partial I^n_{m+m'};\mathcal{Z}^n_{m+m'})
    \to (X,\mathcal{V})$ is bornologous. 
\end{corolario}
    
\begin{proof}
    The result follows from \autoref{cor:Star of triples born} with $A =*$.
\end{proof}

The next lemma enables us to to extend the definition of $\star$ from functions to homotopy classes.
\begin{lema}
\label{WellDefinedClass}
Let $n\geq 2$ and suppose that $\langle f\rangle,\langle g \rangle\in \pi_n^{sc}(X,A,*)$, then
the product
$\langle f\rangle\star \langle g\rangle:= \langle f\star g\rangle$ is well-defined. 
Similarly, if $n\geq1$ and $\langle f\rangle,
\langle g \rangle\in \pi_n^{sc}(X,*)$, then the product $\langle f\rangle\star \langle g\rangle:= \langle f\star g\rangle$
is well-defined.
\end{lema}

\begin{proof}
    Let $n\geq 2$, and suppose that $\langle f \rangle,\langle f'\rangle, \langle g \rangle, \langle g'
    \rangle \in \pi_n^{sc}(X,A,*;\mathcal{V})$ such that $\langle f \rangle = \langle f'\rangle$ and
    $\langle g \rangle = \langle g' \rangle$. Then, by \autoref{rem:Equivalence of classes}, there 
    exist non-negative integers $m$ and $m'$ and functions $\tilde{f},\tilde{f}':
    (I^n_{m},\partial I^n_{m},J^{n-1}_{m};\mathcal{Z}^n_m)
    \rightarrow (X,A,*;\mathcal{V})$ and $\tilde{g},\tilde{g}':
    (I^n_{m'},\partial I^n_{m'},J^{n-1}_{m'};\mathcal{Z}^n_{m'})\to (X,A,*;\mathcal{V})$
    such that $\tilde{f}\simeq_{sc} \tilde{f}'$ and $\tilde{g} \simeq_{sc} \tilde{g}'$.
    We wish to show that $\tilde{f}\star \tilde{g} \simeq_{sc} \tilde{f}'\star\tilde{g}'$.

    Let $H_1$ be a homotopy between $\tilde{f}$ and $\tilde{f}'$ and $H_2$ be a homotopy between
    $\tilde{g}$ and $\tilde{g}'$. Define a function 
    \[
        H:(I^n_{m+m'}\times \Z,\partial I^n_{m+m'}\times \Z,J^{n-1}_{m+m'}\times \Z;
        \mathcal{Z}^n_{m+m'}\times \mathcal{Z})
    \to (X,A,*;\mathcal{V})
    \] 
    by
    \begin{equation*}
        H(\alpha,n) =\begin{cases}
            H_1(\alpha+m'e_1,n) & -m'-m \leq \alpha_1 \leq -m'+m\\
            H_2(\alpha+me_1,n) & -m'+m < \alpha_1 \leq m'+m.
        \end{cases}
    \end{equation*}
    Then $H$ is a map of triples by definition, and $H$ is bornologous 
    by \autoref{prop:Bornologous function from subspaces} and $H$ is a homotopy
    between $\tilde{f} \star \tilde{g}$ and $\tilde{f}'\star \tilde{g}'$. By \autoref{rem:Equivalence of classes},
    it follows that $\langle f \rangle \star \langle g \rangle = \langle f \star g \rangle$ is well-defined.

    For $\pi_n^{sc}(X,*;\mathcal{V})$, the above proof also applies for any $n\geq 1$, with the modification 
    that all of the maps are bornologous maps of pairs instead of maps of triples.
\end{proof}

Given \autoref{WellDefinedClass}, we make the following definition.

\begin{definicion}
For $n\geq 2$ and $\langle f \rangle,\langle g \rangle \in \pi_n^{sc}(X,A,*;\mathcal{V})$, we define
$\langle f \rangle \star \langle g \rangle\coloneqq  \langle f\star g\rangle$.

Similarly, for $n\geq 1$ and $\langle f \rangle,\langle g \rangle \in \pi_n^{sc}(X,*;\mathcal{V})$, we define
$\langle f \rangle \star \langle g \rangle\coloneqq  \langle f\star g\rangle$.

\end{definicion}

We will now show that the homotopy classes $\pi_n^{sc}$ with the $\star$ operation form a group given assumptions on $n$ identical to those in the topological case. This theorem is stated below as \autoref{thm:HomotopyGroups}, and the proof is an adaptation of the proof for topological spaces. The semi-coarse setting introduces a number of technical subtleties, however, which must be dealt with in order to give a complete proof of this theorem. While not difficult, they are nonetheless intricate, and are handled in \autoref{lem:DispPlate}-\autoref{lem:DispWakeWStepBlock} below. The next definition establishes notation for these lemmas.

\begin{definicion}
	In the following definitions, we denote by $\mathbf{a}$ the point \[\mathbf{a}=(a_1^\leftarrow,a_1^\rightarrow,a_2^\leftarrow,a_2^\rightarrow,\cdots,a_n^\leftarrow,a_n^\rightarrow)\in I_m^{2n},\] where $a_i^\leftarrow\leq a_i^\rightarrow$ for all $i\in\{1,\ldots,n\}$. 
	\begin{enumerate}[leftmargin=*,label=(\roman*)]
		\item We call the set $\mathcal{B}_\mathbf{a}$
		\begin{align*}
			\mathcal{B}_\mathbf{a}:= \{ a\in I_m^n: a_i^\leftarrow\leq a_i\leq a_i^\rightarrow, i\in\{1,\ldots,n\} \}.
		\end{align*}
		the \emph{block delimited by $\mathbf{a}$}. We write $\mathcal{B}_{\mathbf{a}}^{\hat{n},k}$ for the set
		$\{ a\in\mathcal{B}_{\mathbf{a}}: a_{\hat{n}}=k\}$. When $a_{\hat{n}}^\leftarrow=a_{\hat{n}}^\rightarrow$ for some $\hat{n} \in \{1,\dots,n\}$, we denote $\mathcal{B}_{\mathbf{a}}$ by $\mathcal{P}_{\mathbf{a}}^{\hat{n}}$, which we call a \emph{plate}. 
		
		\item Let $\mathcal{B}_\mathbf{a}$ be the block delimited by $\mathbf{a}$. We denote by $\partial\mathcal{B}_{\mathbf{a}}$ the set
		\begin{align*}
			\{a\in\mathcal{B}_{\mathbf{a}}: \exists k \in \{1,\dots,n\} \mbox{ with } a_k\in\{a_k^\leftarrow,a_k^\rightarrow\}\}
		\end{align*}
		which we call the \emph{boundary of $\mathcal{B}_\mathbf{a}$}.
		\item Let $\mathcal{B}_\mathbf{a}$ be the block delimited by $\mathbf{a}$, let $\hat{n}\in\{1,\cdots,n\}$, and suppose that $\epsilon \in \{\leftarrow, \rightarrow\}$. We will call \emph{$J^{\hat{n},\epsilon}_{\mathbf{a}}\subset \partial \mathcal{B}_\mathbf{a}$ the $(\hat{n},\epsilon)$-open box of $\mathcal{B}_\mathbf{a}$}, which we define to be the set
		\begin{align*}
			J^{\hat{n},\epsilon}_{\mathbf{a}}:=\{a\in\partial \mathcal{B}_\mathbf{a} \mid \text{If } a_{\hat{n}} = a_{\hat{n}}^\epsilon, \text{ then } \exists k\neq \hat{n} \text{ such that } a_k \in \{a_k^\leftarrow,a_k^\rightarrow\}\}.
		\end{align*}
		In particular, $J^{n-1}_m$ from \autoref{def:nmCubos} is $J^{n,\leftarrow}_\mathbf{a}$, where $\mathbf{a}$ is given by $\mathbf{a} = \{0,m,0,m,\dots,0,m\}$.
		\item $e_k\in I_m^n$ is such that $(e_k)=\delta_{jk}$, the Kroenecker delta. For $a,b\in I_m^n$, $a+b$ is such that $(a+b)_k=a_k+b_k$. For $A\subset I_m^n$ and $a\in I_m^n$, $A+a:=\{b+a:b\in A\}$.
		\item 
		Let $\rho_\epsilon(\cdot):\{\rightarrow,\leftarrow\}\rightarrow \{-1,1\}$ be the set function defined by
		\begin{align*}
			\rho(\epsilon) = \begin{cases}	-1, & \epsilon = \;\leftarrow\\
				1, & \epsilon = \;\rightarrow
			\end{cases} 
		\end{align*}
		We will often refer to $\rho(\epsilon)$ as $\rho_\epsilon$.
		\item For every element $\beta\in I^n_m$, let $\text{Nbh}(\beta)\coloneqq \{(\gamma,\beta)\mid (\gamma,\beta) \in \mathcal{I}^n_m \}$.
	\end{enumerate}
\end{definicion}

When we move a plate one step, this plate will leave a copy of the plate in the original position, as well as shifting a second copy one step. The next lemma shows that this process is bornologous.

\begin{lema}
	\label{lem:DispPlate}
	Let $X$ be a semi-coarse space, let $f:I_m^n \rightarrow (X,\mathcal{V})$ be a bornologous function, and let $\mathcal{P}_\mathbf{a}$ be a plate on $n^*$ delimited by $\mathbf{a}$. Suppose that $\epsilon \in \{\leftarrow,\rightarrow\}$. Suppose that, for every $\alpha\in\mathcal{P}_\mathbf{a}$, $f(\alpha)$ is adjacent to every point in the set $\{ f(\beta) \mid \beta \text{ is adjacent to } \alpha + \rho_\epsilon e_{n^*}\}$, and define 
	\begin{align*}
		f^\epsilon_{\mathcal{P}_\mathbf{a}}(\beta)\coloneqq\begin{cases}
			f(\beta -\rho_\epsilon e_{n^*}) & \text{ if } \beta\in\mathcal{P}_\mathbf{a}+\rho_\epsilon e_{n^*}\\
			f(\beta)                    & \text{ if } \beta\notin\mathcal{P}_\mathbf{a}+\rho_\epsilon e_{n^*}.
		\end{cases}
	\end{align*} 
	Then $f^\epsilon_{\mathcal{P}_{\mathbf{a}}}:I_m^n\rightarrow X$ bornologous and $f \simeq_{sc} f^\epsilon_{\mathcal{P}_\mathbf{a}}$.
	
\end{lema}

\begin{proof}
	We show that $f^\epsilon_{\mathcal{P}_\mathbf{a}}$ is bornologous. We check that $f^{\epsilon}_{P_{\mathbf{a}}}$ sends adjacent points in $I^n_m$ to adjacent points in $X$.
	
	First, let $\beta\notin \mathcal{P}_{\mathbf{a}}\cup (\mathcal{P}_{\mathbf{a}}+\rho_\epsilon e_{n^*})\cup (\mathcal{P}_{\mathbf{a}} +2\rho_\epsilon e_{n^*})$. Then $f^\epsilon_{\mathcal{P}_{\mathbf{a}}}|_{\text{Nbh}(\beta)}=f|_{\text{Nbh}(\beta)}$.
	
	Now suppose that $\beta\in \mathcal{P}_{\mathbf{a}}$.
	\begin{align*}
		f^\epsilon_{\mathcal{P}_{\mathbf{a}}}(\text{Nbh}(\beta)) = &\, f^\epsilon_{\mathcal{P}_{\mathbf{a}}}(\{(\gamma,\beta)\in \text{Nbh}(\beta)\mid \gamma\in \mathcal{P}_{\mathbf{a}}+\rho_\epsilon e_{n^*}\})\\
		&\cup f^\epsilon_{\mathcal{P}_{\mathbf{a}}}(\{(\gamma,\beta)\in \text{Nbh}(\beta)\mid \gamma\notin \mathcal{P}_{\mathbf{a}}+\rho_\epsilon e_{n^*}\})  \\
		= &\, f(\{(\gamma-\rho_\epsilon e_{n^*}),\beta) \mid (\gamma,\beta)\in \text{Nbh}(\beta), \gamma\in \mathcal{P}_{\mathbf{a}}+\rho_\epsilon e_{n^*}\}) \\
		&\cup f(\{(\gamma,\beta)\in \text{Nbh}(\beta)\mid \gamma\notin \mathcal{P}_{\mathbf{a}}+\rho_\epsilon e_{n^*}\}) \\
		\subset  &\, f(\text{Nbh}(\beta)) \in \mathcal{V}.
	\end{align*}
	
	Next, let $\beta\in \mathcal{P}_{\mathbf{a}} + \rho_\epsilon e_{n^*}$. Then
	\begin{align*}
		f^\epsilon_{\mathcal{P}_{\mathbf{a}}}(\text{Nbh}(\beta)) = &\, f^\epsilon_{\mathcal{P}_{\mathbf{a}}}(\{(\gamma,\beta)\in \text{Nbh}(\beta)\mid \gamma\in \mathcal{P}_{\mathbf{a}}\cup \mathcal{P}_{\mathbf{a}}+2\rho_\epsilon e_{n^*}\})\\
		& \cup f^\epsilon_{\mathcal{P}_{\mathbf{a}}}(\{(\gamma,\beta)\in \text{Nbh}(\beta)\mid \gamma\in \mathcal{P}_{\mathbf{a}}+\rho_\epsilon e_{n^*}\})\\
		= &\, \{(f(\gamma),f(\beta-\rho_\epsilon e_{n^*}))\mid (\gamma,\beta)\in \text{Nbh}(\beta),\gamma\in \mathcal{P}_{\mathbf{a}} \cup \mathcal{P}_{\mathbf{a}}+2\rho_\epsilon e_{n^*}\}\\
		&\, \cup \{(f(\gamma-\rho_\epsilon e_{n^*}),f(\beta-\rho_\epsilon e_{n^*}))\mid (\gamma,\beta)\in \text{Nbh}(\beta),\gamma\in \mathcal{P}_{\mathbf{a}}\}\\
		\subset & \,f(\text{Nbh}(\beta - \rho_\epsilon e_{n^*})),
	\end{align*}
	where the last set is in $\mathcal{V}$ by hypothesis. Then
	
	Finally, let $\beta\in \mathcal{P}_{\mathbf{a}} + 2\rho_\epsilon e_{n^*}$.
	\begin{align*}
		f^\epsilon_{\mathcal{P}_{\mathbf{a}}}(\text{Nbh}(\beta)) = & \, f^\epsilon_{\mathcal{P}_{\mathbf{a}}}(\{(\gamma,\beta)\in \text{Nbh}(\beta)\mid \gamma\in \mathcal{P}_{\mathbf{a}}+\rho_\epsilon e_{n^*}\})\\
		& \cup f^\epsilon_{\mathcal{P}_{\mathbf{a}}}(\{(\gamma,\beta)\in \text{Nbh}(\beta)\mid \gamma\notin \mathcal{P}_{\mathbf{a}}+\rho_\epsilon e_{n^*}\})  \\
		= &\, \{(f(\gamma-\rho_\epsilon e_{n^*}),f(\beta))\mid (\gamma,\beta)\in \text{Nbh}(\beta), \gamma\in \mathcal{P}_{\mathbf{a}}+\rho_\epsilon e_{n^*}\} \\
		& \, \cup f(\{(\gamma,\beta)\in \text{Nbh}(\beta)\mid \gamma\notin \mathcal{P}_{\mathbf{a}}+\rho_\epsilon e_{n^*}\})\\
		\subset & \, f(\text{Nbh}(\beta)),
	\end{align*}
	where the last set is in $\mathcal{V}$ by hypothesis. Since $I_m^n$ is a roofed semi-coarse space and has a finite number of elements, the above verifications imply that $f^{\epsilon}_{P_\mathbf{a}}$ is bornologous by \autoref{prop:Bornologous function from subspaces}.
	
	It remains to show that $f\simeq_{sc}f^\epsilon_{\mathcal{P}_\mathbf{a}}$. Let $H:I_m^n\times\mathbb{Z}\rightarrow X$ such that $H(\cdot,k)=f(\cdot)$ for every $k\leq 0$ and $H(\cdot,k)=f^\epsilon_{\mathcal{P}_\mathbf{a}}(\cdot)$ for every $k\geq 1$. We wish to show that $H$ is bornologous. By \autoref{prop:Bornologous function from subspaces}, it is enough to show that the restrictions of $H$ to $I_m^n\times \{z\leq 0\}$, $I_m^n\times \{0,1\}$, and $I_m^n\times \{z\geq 1\}$ are bornologous. By definition, $H(\mathcal{I}_m^n\boxtimes \{z\leq 0\}^2)=f(\mathcal{I}_m^n)$ and $H(\mathcal{I}_m^n\boxtimes \{z\geq 1\}^2)=f^\epsilon_{\mathcal{P}_\mathbf{a}}(\mathcal{I}_m^n)$, so $H$ is bornologous on these regions. It remains to check the behavior of $H$ on $I_m^n\times \{0,1\})$. Since $I_m^n\times \{0,1\}$ has a finite number of elements, we proceed as above by checking that $H$ sends adjacent elements to adjacent elements.
	
	Note that $H(\mathcal{I}_m^n \boxtimes (0,0))=f(\mathcal{I}_m^n)$ and $H(\mathcal{I}_m^n \boxtimes (1,1))=f^\epsilon_{\mathcal{P}_\mathbf{a}}(\mathcal{I}_m^n)$, so $H$ is also bornologous on these regions. Therefore, it only remains to check $H$ on $\mathcal{I}_m^n \boxtimes (0,1)$. Let $\beta\in I_m^n$. We consider the following two cases separately: either $\beta \in \mathcal{P}_\mathbf{a}+\rho_\epsilon e_{n^*}$  or $\beta\notin\mathcal{P}_\mathbf{a}+\rho_\epsilon e_{n^*}$.
	\begin{enumerate}
		\item Case 1: $\beta \in \mathcal{P}_\mathbf{a} + \rho_\epsilon e_{n*}$. Then $H(\text{Nbh}(\beta)\boxtimes (0,1)) = \{(f(\gamma),f_{\mathcal{P}_\mathbf{a}}^\epsilon(\beta))\mid(\gamma,\beta)\in \text{Nbh}(\beta)\}=\{(f(\gamma),f_(\beta-\rho_\epsilon e_{n^*}))\mid(\gamma,\beta)\in \text{Nbh}(\beta)\}\in\mathcal{V}$, which was shown in the first part of the proof.
		
		\item Case 2: $\beta\notin\mathcal{P}_\mathbf{a}+\rho_\epsilon e_{n^*}$. Then $H(\text{Nbh}(\beta)\boxtimes (0,1)) = \{(f(\gamma),f_{\mathcal{P}_\mathbf{a}}^\epsilon(\beta))\mid(\gamma,\beta)\in \text{Nbh}(\beta)\}= \{(f(\gamma),f(\beta))\mid(\gamma,\beta)\in \text{Nbh}(\beta)\} \in \mathcal{V}$.
	\end{enumerate} 
	From this, we conclude that the image of $H$ of the roof of $\mathcal{I}_m^n \boxtimes \mathcal{Z}$ is controlled, so $H$ is bornologous.
\end{proof}

\begin{lema}[Displacement with wake of a block.]
	\label{DispWakeBlock}
	Let $(X,\mathcal{V})$ be a semi-coarse space, $f:I_m^n\rightarrow X$ a bornologous function, let $\mathcal{B}$ be a block delimited by $\mathbf{a}$, and let $n^*\in\{1,\cdots,n\}$.
	Define
	\begin{align*}
		f_\mathcal{B}^\epsilon(\beta)=\left\lbrace\begin{array}{rr}
			f(\beta-\rho_\epsilon e_{n^*}) & \mbox{ if } \beta\in\mathcal{B}+\rho_\epsilon e_{n^*}\\
			f(\beta)         & \mbox{ if } \beta\notin\mathcal{B}+\rho_\epsilon e_{n^*}
		\end{array}\right.
	\end{align*}
	Then $f^*:I_m^n\rightarrow X$ is bornologous and $f^*\simeq_{sc} f$ if for each $\alpha\in J^{n^*,\epsilon}_{\mathbf{a}}$, $f(\alpha)$ is adjacent to the image under $f$ for each adjacent element to $\alpha+ \rho_\epsilon e_{n^*}$ which are not in $\mathcal{B}$.
	
\end{lema}

\begin{proof}
	Let $(X,\mathcal{V})$, $f:I_m^n\rightarrow X$, $\mathcal{B}$, and
	$n^*\in\{1,\ldots,n\}$ be as in the statement of the lemma. For $j\in\{a_{n^*}^\leftarrow,a_{n^*}^\leftarrow+1,\ldots,a_{n^*}^\rightarrow\}$, define $f_j:I_m^n\rightarrow X$ by
	\begin{align*}
		f_j(\beta)=\left\lbrace \begin{array}{rl}
			f_{j+1}(\beta-\rho_\epsilon e_{n^*}) & \mbox{ if } \beta\in\mathcal{B}_{n^*}^{j}+\rho_\epsilon e_{n^*}\\
			f_{j+1}(\beta)         & \mbox{ if } \beta\notin\mathcal{B}_{n^*}^{j}+\rho_\epsilon e_{n^*}
		\end{array}\right.
	\end{align*}
	where $f_{a_{n^*}^\rightarrow+1}:=f$. 
	
	We proceed by induction.
	
	\textbf{Base case:} By hypothesis, for each $\alpha\in\mathcal{B}_{n^*}^{a_{n^*}^\epsilon}$, $f(\alpha)$ is adjacent to the image under $f$ of adjacent elements of $\alpha+\rho_\epsilon e_{n^*}$. By \autoref{lem:DispPlate}, $f_{a_{n^*}^\epsilon}$ is bornologous and homotopy equivalent to $f$.
	
	\textbf{Inductive step:} Suppose that $f_{J+1}$ is bornologous and homotopy equivalent to $f$. We show that this implies that the same is true for $f_J$. Note that for $\alpha\in\mathcal{B}_{n^*}^{J+2}$ we have that $f_{J+1}(\alpha)=f_{J+1}(\alpha-\rho_\epsilon e_{n^*})=f_J(\alpha-\rho_\epsilon e_{n^*})$, and therefore $f_J$ restricted to $\mathcal{B}_{n^*}^J$ satisfies \autoref{lem:DispPlate}. It follows that $f_J$ is bornologous and $f_J\simeq_{sc} f_{J+1} \simeq_{sc} f$.
\end{proof}

We now use the above to move blocks a finite number of steps.

\begin{lema}[Displacement with wake $k$-steps of a block.]
	\label{lem:DispWakeWStepBlock}
	Let $(X,\mathcal{V})$ be a semi-coarse space, $f:I_m^n\rightarrow X$ a bornologous function, $\mathcal{B}$ be a block delimited by $\mathcal{A}$, $k\in\mathbb{N}$ and $n^*\in\{1,\cdots,n\}$.
	
	Define 
	
	\begin{align*}
		f^\epsilon_\mathcal{B}(\beta)=\left\lbrace\begin{array}{rl}
			f(\beta-k \rho_\epsilon e_{n^*}) & \mbox{ if } \beta\in\mathcal{B}+k \rho_\epsilon e_{n^*}\\
			f(\beta-v \rho_\epsilon e_{n^*}) & \mbox{ if } \beta\in\mathcal{B}_{n^*}^{a_{n^*}^\leftarrow}+v \rho_\epsilon e_{n^*},\ v\in\{1,\ldots,k-1\}\\
			f(\beta) & \mbox{ otherwise}
		\end{array}\right.
	\end{align*}
	Then $f^\epsilon_\mathcal{B}:I_m^n\rightarrow X$ is bornologous and $f^\epsilon_\mathcal{B}\simeq_{sc} f$ if for each $\alpha\in J_{\mathbf{a}}^{n^*,\epsilon}$ we have that $f(\alpha)$ is adjacent to the image under $f$ of the set 
	\begin{equation*}\begin{array}{ll}
			\bigcup_{v\in \{1,\dots,k-1\}} \text{Nbh}(\alpha+v \rho_\epsilon e_{n^*}) & \text{ if }\alpha\in \mathcal{B}_\mathbf{a}^{n^*,a_{n^*}^\epsilon}\\
			\bigcup_{v\in \{1,\dots,k-1\}}\text{Nbh}(\alpha+v \rho_\epsilon e_{n^*}) - \bigcup_{v\in \{1,\dots,k-1\}}(\mathcal{B}+v\rho_\epsilon e_{n^*}) & \text{ if }\alpha\notin \mathcal{B}_\mathbf{a}^{n^*,a_{n^*}^\epsilon}.
		\end{array}
	\end{equation*}
	
\end{lema}

\begin{proof}
	Let $(X,\mathcal{V})$, $f:I_m^n\rightarrow X$, $\mathcal{B}$, $k\in\mathcal{N}$, and $n^*\in\{1,\ldots, n\}$ be as in the hypothesis of the lemma. Assume that
	for each $\alpha\in J^{n^*,\epsilon}_{\mathbf{a}}$ we have that $f(\alpha)$ is adjacent to the image under $f$ of the adjacent elements of $\alpha+ve_{n^*}$, $v\in\{1,\cdots,w\}$.
	
	For $j\in\{0,\ldots,k-1\}$, we define
	\begin{align*}
		\mathcal{B}^{\langle j \rangle}=& \mathcal{B}+j \rho_\epsilon e_{n^*}\\
		f_{j+1} := & \left\lbrace \begin{array}{rl}
			f_{j}(\beta-\rho_\epsilon e_{n^*}) & \mbox{if } \beta\in\mathcal{B}^{\langle j\rangle}+\rho_\epsilon e_{n^*}\\
			f_{j}(\beta) & \mbox{if } \beta\notin\mathcal{B}^{\langle j\rangle}+\rho_\epsilon e_{n^*}
		\end{array}\right.
	\end{align*}
	Let $f_0:=f$ and observe that $\mathcal{B}^{\langle 0 \rangle}=\mathcal{B}$. We proceed by induction.
	
	\textbf{Base case:} Observe that, for $v=1$, the functions $f$, $f_1$ and the block $\mathcal{B}$ satisfy the hypotheses of \autoref{DispWakeBlock}, from which it follows that $f_1$ is bornologous and $f_1\simeq_{sc} f$.
	
	\textbf{Inductive step:} Assume that $f_{J-1}$ is bornologous and homotopy equivalent to $f$. Observe that, if $v=J$, then $f_{J-1}$, $f_J$ and $\mathcal{B}^{\langle J-1\rangle}$ satisfy the hypotheses of \autoref{DispWakeBlock}, from which it follows that $f_J$ is bornologous and $f_J\simeq_{sc} f_{J-1}$.
	
	Since $f^\epsilon_\mathcal{B}=f_w$, we conclude that $f^\epsilon_\mathcal{B}$ is bornologus and $f^\epsilon_\mathcal{B}\simeq_{sc} f$.\qedhere
\end{proof}

The above lemmas now allow us to prove the following theorem. 
The proof is analogous to that in the topological case, replacing the $n$-disk $D^n$ with the blocks $I^n_m$.

\begin{teorema}[Semi-coarse Homotopy Groups]
\label{thm:HomotopyGroups}
Let $(X,\mathcal{V})$ be a semi-coarse space, $A\subset X$, and $n\in\mathbb{N}$.
\begin{itemize}
\item If $n\geq 1$, then $(\pi_n^{sc}(X,*),\star)$ is a group.
\item If $n\geq 2$, then $(\pi_n^{sc}(X,*),\star)$ is an abelian group.
\item If $n\geq 2$, then $(\pi_n^{sc}(X,A,*),\star)$ is a group.
\item If $n\geq 3$, then $(\pi_n^{sc}(X,A,*),\star)$ is an abelian group.
\end{itemize}
\end{teorema}

\begin{proof}
Let $(X,\mathcal{V})$ be a semi-coarse space,
let $n\geq 2$, and suppose that there are $m,m',m''$ such that $f:(I_{m}^n,\partial I_{m}^n,J_{m}^{n-1})\rightarrow (X,A,*)$, $g:(I_{m'}^n,\partial I_{m'}^n,J_{m'}^{n-1})\rightarrow (X,A,*)$ and $h:(I_{m''}^n,\partial I_{m''}^n,J_{m''}^{n-1})\rightarrow (X,A,*)$. Recall from \autoref{def:OpStar},
that $(f\star g)\star h=f\star(g\star h)$. Thus, by \autoref{WellDefinedClass},
\begin{align*}
\langle f\star g \rangle\star \langle h \rangle=\langle (f\star g)\star h \rangle=\langle f\star(g\star h) \rangle=\langle f\rangle \star \langle g\star h \rangle,
\end{align*}
so the product $\star$ is associative.

 Define $e:I_0^n\rightarrow X$ to be the constant map $e(x)=*$, and consider the function $f:(I_m^n,\partial I_m^n,J_m^{n-1})\rightarrow (X,A,*)$, a bornologous map of triples. By \autoref{def:OpStar} of $\star$, $f\star e = f = e \star f$. Thus, $\langle f \star e \rangle = \langle f \rangle=\langle e\star f \rangle$.

 We now show that every bornologous map $f:(I_m^n,\partial I_m^n,J_m^{n-1})\rightarrow (X,A,*)$ has a homotopy inverse. Define $g:(I_m^n,\partial I_m^n,J_m^{n-1})\rightarrow (X,A,*)$ to be the map $g(\alpha)=f(\beta)$, where $\alpha,\beta \in I_m^n$ and the coordinates of $\alpha$ and $\beta$ satisfy 
 \begin{align*}
 \alpha_1&=m-\beta_1,\\ 
 \alpha_i&=\beta_i \text{ for }i\neq 1.
 \end{align*} 

We show that $\langle f\star g\rangle=\langle e \rangle$ by induction, and we remark that we can proceed in the same way to show that $\langle g \star f\rangle = \langle e \rangle$. Define $f_m\coloneqq f\star g$, and define $f_j:(I_{2m}^n,\partial I_{2m}^n,J_{2m-1}^n)\rightarrow (X,A,*)$ for $j\in\{1,\ldots,m-1\}$ recursively by
\begin{align*}
f_{j-1}(\alpha):= \begin{cases}
f_{j}(\alpha), & \mbox{if }0\leq \alpha_1\leq j-1\\
f_{j}(-\alpha+2(j-1)e_1), & \mbox{if } j\leq\alpha_1\leq 2(j-1)\\
*, & \mbox{elsewhere}.
\end{cases}
\end{align*}

 Let $\mathcal{B}^{\langle j\rangle}$ the block delimited by 
\begin{align*}
\mathbf{a}^{\langle j \rangle}:=(j+1,2j,0,2m,\dots,0,2m) \in I_{2m}^n
\end{align*}
\textbf{Base case:} $f_{m-1}$ is a bornologous function and homotopy equivalent to $f'\star(f')^{-1}$ by \autoref{lem:DispWakeWStepBlock}, taking $\mathcal{B}=\mathcal{B}^{\langle m \rangle}$, $k=2$ and $f=f'\star(f')^{-1}$. To be able to use that lemma, observe that for all $\alpha$ such that $\alpha_1=m+1$, $f(\alpha)$ is adjacent to the image of elements adjacent to $\alpha-e_1$ (which is $*$) and to the image of elements adjacent to $\alpha-2e_1$.

\textbf{Induction step:} Suppose that $f_{j}$ is a bornologous function and that it is homotopy equivalent to $f_{j+1}$. By \autoref{lem:DispWakeWStepBlock}, $f_{j-1}$ is bornologus and homotopy equivalent to $f_j$, taking $\mathcal{B}=\mathcal{B}^{\langle j\rangle}$, $k=2$, and $f=f_j$ in the lemma. Also observe that, for all $\alpha$ such that $\alpha_1=j+1$, $f_j(\alpha)$ is adjacent to the image of elements adjacent to $\alpha-e_1$ and $\alpha-2e_1$.

The above implies that $f \star g \simeq_{sc} {*}$, and therefore $(\pi_{n}^{sc}(X,A,*),\star)$ is a group. 

Now suppose, in addition, that $n\geq 3$, and let $\langle f \rangle, \langle g \rangle \in \pi_n^{sc}(X,A,*)$ such that $m$ and $m'$ are integers with $f:(I_{m}^n,\partial I_{m}^n,J_{m}^{n-1})\rightarrow (X,A,*)$ and $g:(I_{m'}^n,\partial I_{m'}^n,J_{m'}^{n-1})\rightarrow (X,A,*)$.

For $i\in\{1,2,3,4,5\}$, let $\mathcal{B}^{\langle i \rangle}$ be the blocks delimited by $\mathbf{a}^{\langle i \rangle}$, where
\begin{align*}
\mathbf{a}^{\langle 1 \rangle} := & (m+1,m+m'-1,0,m'-1,0,m+m',\dots,0,m+m')\\
\mathbf{a}^{\langle 2 \rangle} := & (m+1,m+m',m+1,m'-1,0,m+m',\dots,0,m+m')\\
\mathbf{a}^{\langle 3 \rangle} := & (1,m'-1,m+1,m+m',0,m+m',\dots,0,m+m')\\
\mathbf{a}^{\langle 4 \rangle} := & (0,m-1,1,m-1,0,m+m',\dots,0,m+m')\\
\mathbf{a}^{\langle 5 \rangle} := & (m'+1,m+m'-1,1,m-1,0,m+m',\dots,0,m+m').
\end{align*}
Define
\begin{align*}
f_1(\beta):=\left\lbrace\begin{array}{rl}
f\star g(\beta-me_2) & \mbox{if }\beta\in\mathcal{B}^{\langle 1 \rangle}+me_2\\
f\star g(\beta-ve_2) & \mbox{if }\beta\in(\mathcal{B}^{\langle 1 \rangle})_2^{0}+ve_2,\ v\in\{1,\ldots,m-1\}\\
f\star g(\beta) & \mbox{anywhere else}
\end{array}\right.
\end{align*}
Observe that $f\star g((\mathcal{B}^{\langle 1\rangle}\cup\mathcal{B}^{\langle 4\rangle})^c)=\{*\}$, then by \autoref{lem:DispWakeWStepBlock} we have that $f_1$ is bornologous and $f_1\simeq_{sc} f\star g$. Let's define
\begin{align*}
f_2(\beta):=\left\lbrace\begin{array}{rl}
f_1(\beta+me_1) & \mbox{if }\beta\in\mathcal{B}^{\langle 2 \rangle}-me_1\\
f_1(\beta+ve_1) & \mbox{if }\beta\in(\mathcal{B}^{\langle 2 \rangle})^{m+m'}_1-ve_1,\ v\in\{1,\ldots,m-1\}\\
f_1(\beta) & \mbox{anywhere else}
\end{array}\right.
\end{align*}
Observe that $f_1((\mathcal{B}^{\langle 2 \rangle}\cup\mathcal{B}^{\langle 4 \rangle})^c)=\{*\}$, then by \autoref{lem:DispWakeWStepBlock} we have that $f_2$ is bornologous and $f_2\simeq_{sc}f_1$. Let's define
\begin{align*}
f_3(\beta):=\left\lbrace\begin{array}{rl}
f_2(\beta-m'e_1) & \mbox{if }\beta\in\mathcal{B}^{\langle 4 \rangle}+m'e_1\\
f_2(\beta-ve_1) & \mbox{if }\beta\in(\mathcal{B}^{\langle 4 \rangle})_1^0+ve_1,\ v\in\{1,\ldots,m'-1\}\\
f_2(\beta) & \mbox{anywhere else}
\end{array}\right.
\end{align*}
Observe that $f_2((\mathcal{B}^{\langle 3 \rangle}\cup\mathcal{B}^{\langle 4 \rangle})^c)=\{*\}$, then by \autoref{lem:DispWakeWStepBlock} we have that $f_3$ is bornologous and $f_3\simeq_{sc}f_2$. Let's define
\begin{align*}
f_4(\beta):=\left\lbrace\begin{array}{rl}
f_3(\beta+me_2) & \mbox{if }\beta\in\mathcal{B}^{\langle 3 \rangle}-me_2\\
f_3(\beta+ve_2) & \mbox{if }\beta\in(\mathcal{B}^{\langle 3 \rangle})_2^{m+m'}-ve_2, v\in\{1,\ldots,m-1\}\\
f_3(\beta) & \mbox{anywhere else}
\end{array}\right.
\end{align*}
Observe that $f_3((\mathcal{B}^{\langle 3 \rangle}\cup\mathcal{B}^{\langle 5 \rangle})^c)=\{*\}$, then by \autoref{lem:DispWakeWStepBlock} we have that $f_4$ is bornologous and $f_4\simeq_{sc}f_3$. Finally, note that $f_4=g\star f$, and we conclude that $\pi_n^{sc}(X,A,*)$ is commutative.

If $A=\{*\}$, then we have shown that $\pi_n^{sc}(X,*)$ is a group for $n\geq 2$ and abelian for $n\geq 3$. Consider the case where $n=1$, and let $\langle f \rangle,\langle g \rangle, \langle h \rangle\in \pi_1^{sc}(X,*)$.  Then there exist non-negative integers $m,m',m''$ such that $f:(I_{m},\partial I_{m})\rightarrow (X,*)$, $g:(I_{m'},\partial I_{m'})\rightarrow (X,*)$ and $h:(I_{m''},\partial I_{m''})\rightarrow (X,*)$. By definition, $(f\star g)\star h=f\star(g\star h)$, and therefore, by \autoref{WellDefinedClass},
\begin{align*}
[\langle f'\star g' \rangle]\star [\langle h' \rangle]=[\langle (f\star g)\star h \rangle]=[\langle f\star(g\star h) \rangle]=[\langle f'\rangle]\star[\langle g'\star h' \rangle],
\end{align*}
so $\star$ is associative.

Define $e:I_0\rightarrow X$ to be the map which sends $I_0$ to $*$. Therefore, $\langle e \rangle\in \pi_1^{sc}(X,*)$. Let $\langle f \rangle \in \pi_1^{sc}(X,*)$, so there exists a non-negative integer $m$ with $f:(I_m^n,\partial I_m^n)\rightarrow(X,*)$. By definition of $\star$, $f\star e = f = e \star f$. Thus, $\langle f \star e \rangle = \langle f \rangle=\langle e\star f \rangle$ and $\langle e \rangle$ is the identity in $\pi_1^{sc}(X,*)$.

Let $\langle f \rangle \in \pi_1^{sc}(X,*)$. Then there exists a non-negative integer $m$ with $f:(I_m,\partial I_m)\rightarrow (X,*)$. Define $(f')^{-1}$ such that $(f')^{-1}(\alpha)=f'(\beta)$ with $\alpha=m-\beta$.

We show that $\langle f'\star (f')^{-1}\rangle=\langle e \rangle$ by mathematical induction and remark that we can proceed in the same way for $\langle (f')^{-1}\star f'\rangle=\langle e \rangle$. Define $f_j:I_{2m}\rightarrow X$ for $j\in\{1,\ldots,m\}$ by
\begin{align*}
f_{j-1}(\alpha):=\left\lbrace \begin{array}{rl}
f'_{j}(\alpha) & \mbox{if }0\leq \alpha\leq j-1\\
f'_{j}(-\alpha+2(j-1)e_1) & \mbox{if } j\leq\alpha\leq 2(j-1)\\
* & \mbox{anywhere else}
\end{array}\right.
\end{align*}
with $f_m=f'\star(f')^{-1}$, and observe that $f_{j-1}:(I_{2m},\partial I_{2m})\rightarrow (X,*)$ is bornologous. Let $\mathcal{B}^{\langle j\rangle}$ be the block delimited by 
\begin{align*}
\mathbf{a}^{\langle j \rangle}:=(j+1,2j).
\end{align*}

\textbf{Base case:} $f_{m-1}$ is a bornologous function and homotopy equivalent to $f'\star(f')^{-1}$ by \autoref{lem:DispWakeWStepBlock} with $\mathcal{B}=\mathcal{B}^{\langle m \rangle}$, $w=2$ and $f=f'\star(f')^{-1}$. Observe, in addition, that for $\alpha=m+1$, $f(\alpha)$ is adjacent to the image of elements adjacent to $\alpha-e_1$ (which is $*$) and to the image of elements adjacent to $\alpha-2e_1$.

\textbf{Induction step:} Suppose that $f_{j}$ is a bornologous function and homotopy equivalent to $f_{j+1}$. By \autoref{lem:DispWakeWStepBlock}, $f_{j-1}$ is bornologous and homotopy equivalent to $f_j$, taking $\mathcal{B}=\mathcal{B}^{\langle j\rangle}$, $w=2$, and $f=f_j$ in the lemma. Observe that, for $\alpha=j+1$, $f_{j}(\alpha)$ is adjacent to the image of elements adjacent to $\alpha-e_1$ and $\alpha-2e_1$. 

We conclude that $(\pi_1^{sc}(X,*),\star)$ is a group.

Finally, if $n=2$, then we may show that $\pi_2^{sc}(X,*)$ is abelian by repeating the above proof of this fact for $\pi_n^{sc}(X,A,*)$ for $n\geq 3$. Observe that the proof is not valid for $\pi_2^{sc}(X,A,*)$ because $f\star g(\alpha)\in A$ when $\alpha_2=0$. However, if $A=\{*\}$, we have that $f\star g(\alpha)=*$ when $\alpha_2=0$. 
\end{proof}

In the case of coarse spaces, the next theorem shows that the semi-coarse homotopy groups are trivial.

\begin{teorema}
\label{theo:HomotopyClassesCoarse}
Let $(X,\mathcal{V})$ be a coarse space with $A\subset X$. Then for any $n\in\mathbb{N}$, we have $\pi_n^{sc}(X,A,*)\cong\pi_n^{sc}(X,*)\cong \{1\}$.
\end{teorema}

\begin{proof}
Let $(X,\mathcal{V})$ be a coarse space, $n\in\mathbb{N}$, and suppose that $h:I_m^n\rightarrow X$ is the constant map $h(\alpha)=*$ for all $\alpha \in I^n_m$.

Let $\langle f \rangle\in\pi_n^{sc}(X,A,*)$. By definition, $f(J^{n-1})=\{*\}$, in particular $f(\alpha)=*$ when the first coordinate $\alpha_1=0$. Note that, for all $\alpha\in I_m^n$ satisfies that $f(\alpha-ke_1)$ is adjacent to $f(\alpha-(k+1)e_1)$ for any $k\in\{0,\ldots,\alpha_1-1\}$, that is
\begin{align*}
\{(f(\alpha),f(\alpha-e_1)),(f(\alpha-e_1),f(\alpha-2e_1)),\ldots,(f(\alpha-(\alpha_1-1)e_1,f(\alpha-\alpha_1e_1))\}\in\mathcal{V}
\end{align*}
because $f$ is bornologous. Thereby, since $(X,\mathcal{V})$ is a coarse space and $f(\alpha-\alpha_1e_1)=*$, we conclude that $\{(f(\alpha),*)\}\in\mathcal{V}$ for every $\alpha\in I^n_m$, and therefore
\begin{align*}
f \simeq_{sc} *.
\end{align*}
Therefore, $\pi_n^{sc}(X,A,*)\cong \{1\}$.
The proof for $\pi_n^{sc}(X,*)$ is analogous.
\end{proof}

While the above theorem says that semi-coarse homotopy groups are trivial for coarse spaces, we expect there to be other semi-coarse invariants, perhaps defined relative to infinity, which are non-trivial for both coarse and non-coarse semi-coarse spaces. However, we leave this question open for future work.
\subsection{Connectedness}

\begin{definicion}[Semi-Coarse $n$-Connected Space]
\label{QCnConn}
Let $(X,\mathcal{V})$ be a semi-coarse space. We will say that $(X,\mathcal{V})$ is an 
\emph{$n$-connected} (semi-coarse) space iff $\pi_{n}^{sc}(X,*)\cong \{0\}$. 
In particular, a $0$-connected space will be called a \emph{connected} (semi-coarse) space.
\end{definicion}

In the next proposition we will observe that, when the semi-coarse structure is coarse, 
then semi-coarse connectedness is equivalent to coarse connectedness. We first recall the notion
of coarse connectedness. See also \cite{Roe_2003} for more details on coarse connectedness.

\begin{definicion}[Coarse connected space]
\label{CConn}
Let $(X,\mathcal{V})$ be a coarse space. We say that $(X,\mathcal{V})$ is \emph{coarsely connected} iff
every point in $X\times X$ belongs to some controlled set $V \in \mathcal{V}$.
\end{definicion}

\begin{proposicion}
    Let $(X,\mathcal{V})$ be a coarse space. Then $(X,\mathcal{V})$ is coarsely connected iff
    $\pi_0^{sc}(X,*;\mathcal{V}) \cong \{0\}$.
\end{proposicion}

\begin{proof}
Let $(X,\mathcal{V})$ be a coarse space.

$(\Leftarrow)$ If $(X,\mathcal{V})$ is coarsely connected, 
then for each $(x,y)\in X\times X$ there exists $E_{x,y}\in\mathcal{V}$ such that 
$(x,y)\in E_{x,y}$, that is, $\{(x,y)\}\in\mathcal{V}$ for any $x,y\in X$.
Let $f:(\{*,1\},*;\mathcal{D}_{\{*,1\}})\rightarrow (X,*;\mathcal{V})$ and $g:(\{*,1\},*;\mathcal{D}_{\{*,1\}})
\rightarrow (X,*;\mathcal{V})$ be bornologous maps. Note that, since $\mathcal{D}_{\{*,1\}}$ is the 
diagonal semi-coarse structure, $f(1)$ and $g(1)$ may be arbitrary points of $X$.
Define $H:\{*,1\}\times\mathbb{Z}\rightarrow X$ such that $H(\alpha,n)=f(\alpha)$ 
when $n\leq 0$ and $H(\alpha,n)=g(\alpha)$ when $n\geq 1$. Since $X$ is coarsely connected,
we have that $\{(f(1),g(1))\}\in\mathcal{V}$, and therefore
$H$ is bornologous. It follows that $f\simeq_{sc} g$. However, $f$ and $g$ were arbitrary, and we conclude that 
$\pi_0^{sc}(X,*)\cong\{0\}$.

$(\Rightarrow)$ If $(X,\mathcal{V})$ is a coarse space with $\pi_0^{sc}(X,*;\mathcal{V})=\{0\}$,
then for each $f:(\{*,1\},*)\rightarrow (X,*)$ and 
$g:(\{*,1\},*)\rightarrow (X,*)$ there exist $N,M\in\mathbb{Z}$ with $N<M$ and a bornologous 
function $H:\{*,1\}\times\mathbb{Z}\rightarrow X$ such that $H(e,z)=f(z)$ when $z\leq N$ and 
$H(e,z)=g(z)$ when $z\geq M$. Given that $H$ is bornologous, then $H(1,z)$ is adjacent to 
$H(1,z+1)$ for each $z\in\mathbb{Z}$. Therefore
\begin{align*}
\{(f(1),H(1,N+1))\},\{(H(1,N+1),H(1,N+2))\},\ldots \{(H(1,M-1),g(1))\}\in\mathcal{V}
\end{align*}
Then, $\{(f(1),g(1))\}\in\mathcal{V}$. How $f$ and $g$ are arbitrary, $f(1)$ and $g(1)$ are any element of $x$. Thus, for each $(x,y)\in X\times X$, $\{(x,y)\}\in\mathcal{V}$.
\end{proof}

The highlight in this section is being able to note that the base point does not matter when we have a connected space. We are going to prove a lemma which makes that easy.

\begin{lema}
Let $(X,\mathcal{V})$ be a semi-coarse space. If $x,y\in X$ are adjacent elements, then $\pi_n^{sc}(X,x)\cong \pi_n^{sc}(X,y)$ for $n\geq 1$.
\end{lema}

\begin{proof}
Let $(X,\mathcal{V})$ be a semi-coarse space, $n\geq 1$ and $x,y\in X$ be adjacent elements. Let $[\langle f \rangle]:(I^n,\partial I^n)\rightarrow (X,x)$, then there exists a non-negative integer $m$ such that $f':(I_m^n,\partial I_m^n)\rightarrow (X,x)$ and $f'\in \langle f\rangle$. By definition, $(\phi_m^{m+2})^n f'\in \langle f \rangle$.

Let's define $\mathcal{B}$ the block delimited by $\mathcal{A}$ such that $a_i^\leftarrow=1$ and $a_i^\rightarrow=m-1$ and
\begin{align*}
f_1(\beta):= \left\lbrace\begin{array}{rl}
(\phi_m^{m+2})^n f'\left( \beta - \sum_{i=1}^n e_i \right) & \mbox{if }\beta\in\mathcal{B}+\sum_{i=1}^n e_i\\
x & \mbox{anywhere else}
\end{array}\right.
\end{align*}
We get that $(\phi_m^{m+2})^n f'$ is bornologous and $(\phi_m^{m+2})^n f'\simeq_{sc} f_1$. Let's define
\begin{align*}
f_2(\beta):=\left\lbrace\begin{array}{rl}
y & \mbox{if }\beta\in\partial I_{m+2}^n\\
f_1(\beta) & \mbox{anywhere else}
\end{array}\right.
\end{align*}
which is bornologous because $x$ is adjacent to $y$.

Observe that we get a function $\Psi:\pi_n^{sc}(X,x)\rightarrow \pi_n^{sc}(X,y)$ which is well-defined by $(\phi_m^{m+2})^n f'\simeq_{sc} f_1$. We need to prove that this is a group homomorphism.

Let $[\langle f \rangle],[\langle g \rangle]\in \pi_n^{sc}(X,x)$, then there exist non-negative integers $m,m'$ such that $f':(I_m^n,\partial I_m^n)\rightarrow (X,x)$, $g':(I_{m'}^n,\partial I_{m'}^n)\rightarrow (X,x)$, $f'\in \langle f \rangle$ and $g'\in\langle g \rangle$. We just need to prove that $\langle\Psi(f'\star g')\rangle\simeq_{sc} \langle\Psi(f')\star\Psi(g')\rangle$. Define $\mathcal{B}^{\langle 1\rangle}$ the block delimited by $\mathcal{A}^{\langle 1\rangle}$ such that
\begin{align*}
a_1^\leftarrow=m+m'+2,\ a_1^\rightarrow=m+m'+3\\
a_i^\leftarrow=1,\ a_i^\rightarrow=m'+1 &\ i\neq 1
\end{align*}
and $\mathcal{B}^{\langle 2\rangle}$ the block delimited by $\mathcal{A}^{\langle 2\rangle}$ such that
\begin{align*}
a_1^\leftarrow=m+1,\ a_1^\rightarrow=m+m'+1\\
a_i^\leftarrow=1,\ a_i^\rightarrow=m'+1 &\ i\neq 1
\end{align*}
Let's define $h_1:I_{m+m'+4}\rightarrow X$ such that
\begin{align*}
h_1(\beta):= \left\lbrace\begin{array}{rl}
x & \mbox{if }\beta\in \mathcal{B}^{\langle 1\rangle}\\
(\phi_{m+m'+2}^{m+m'+4})^n\Psi(f'\star g')(\beta) & \mbox{anywhere else}
\end{array}\right.
\end{align*}
Note that $h_1\simeq_{sc}(\phi_{m+m'+2}^{m+m'+4})^n\Psi(f'\star g')$. Let's define
\begin{align*}
h_2(\beta):=\left\lbrace\begin{array}{rl}
h_1(\beta-2e_1) & \mbox{if }\beta\in \mathcal{B}^{\langle 2\rangle}+2e_1\\
h_1(\beta-e_1) & \mbox{if }\beta\in(\mathcal{B}^{\langle 2\rangle})_{1}^{m+1}+e_1\\
h_1(\beta) & \mbox{anywhere else}
\end{array}\right.
\end{align*}
which is a bornologous function and $h_1\simeq_{sc}h_2$ 
. Let's define $\mathcal{B}^{\langle 3 \rangle}$ and $\mathcal{B}^{\langle 4 \rangle}$ delimited by $\mathcal{A}^{\langle 3 \rangle}$ and $\mathcal{A}^{\langle 4 \rangle}$ such that
\begin{align*}
\mathcal{A}^{\langle 3 \rangle}:=& \{ \alpha\in I_{m+m'+4}^n: 1\leq\alpha_i\leq m+1 \}\\
\mathcal{A}^{\langle 4 \rangle}:=& \{ \alpha\in I_{m+m'+4}^n: m+3\leq\alpha_1\leq m+m'+3,\ 1\leq\alpha_i\leq m'+1, i\neq 1 \}
\end{align*}
Let's define
\begin{align*}
h_3(\beta):=\left\lbrace\begin{array}{rl}
h_2(\beta) & \mbox{if }\beta\in\mathcal{B}^{\langle 3 \rangle}\cup\mathcal{B}^{\langle 4 \rangle}\\
y &  \mbox{anywhere else}
\end{array}\right.
\end{align*}
which is bornologous and $h_3\simeq_{sc}h_4$ because $x$ is adjacent to $y$. Finally, observe that $h_3=\Psi(f'\star g')$, getting what we want. Thus, $\Psi$ is a group homomorphism.

Under the same argument, we are able to define $\Theta:\pi_n^{sc}(X,y)\rightarrow \pi_n^{sc}(X,x)$. So, we just need to prove that $\Psi^{-1}=\Theta$ and we will get our isomorphism.

Let $[\langle f \rangle]\in\pi_n^{sc}(X,x)$, then there exists a non-negative integer $m$ such that $f':(I_m,\partial I_m)\rightarrow (X,x)$ and $f'\in\langle f\rangle$. Since $x$ and $y$ are adjacent, we are able to replace every $y$ for $x$ in $\Theta\Psi(f')$ calling that function $h$, then $h$ is a bornologous function, $h\simeq_{sc} \Theta\Psi(f')$ and $h\in \langle f'\rangle=\langle f \rangle$. Thus, $\Theta\Psi[\langle f \rangle]=[\langle f\rangle]$, that is, $\Theta\Psi=1_{\pi_n^{sc}(X,x)}$. Under the same argument, $\Psi\Theta=1_{\pi_n^{sc}(X,y)}$, getting that $\Psi^{-1}=\Theta$.
\end{proof}

\begin{teorema}
Let $(X,\mathcal{V})$ be a connected semi-coarse space. Then, for every $x,y\in X$ we get $\pi_n^{sc}(X,x)\cong \pi_n^{sc}(X,y)$.
\end{teorema}

\begin{proof}
Let $(X,\mathcal{V})$ be a connected semi-coarse space and $x,y\in X$. Then, there exist $N,M\in\mathbb{Z}$ such that $N<M$ and a bornologous function $H:\{*,1\}\times\mathbb{Z}\rightarrow X$ such that $H(*,z)=*$ when $z\in\mathbb{Z}$, $H(1,z)=x$ when $z\leq N$ and $H(1,z)=y$ when $z\geq M$. Since $H$ is bornologous, then $\{(x,H(1,N+1)\},\{(H(1,N+1),H(1,N+2))\},\ldots \{(H(1,M-1),y)\}\in\mathcal{V}$. Thus,
\begin{align*}
\pi_n^{sc}(X,x)\cong \pi_n^{sc}(X,H(1,N+1))\cong \ldots \cong \pi_n^{sc}(X,H(1,M-1))\cong \pi_n^{sc}(X,y).
\end{align*}
Getting that $\pi_n^{sc}(X,x)\cong \pi_n^{sc}(X,y)$.
\end{proof}

\subsection{The Semi-Coarse Fundamental Group of Cyclic Graphs}

In this section, we will compute the semi-coarse fundamental group of cyclic graphs with different structures, and, in particular, this will
provide a class of examples where semi-coarse homotopy is non-trivial in dimensions greater than zero. (For the $0$-homotopy class, we get that $\#(\pi_0^{sc}(X,*))=\# X$,
the number of semi-coarse connected components of $X$.)

\begin{definition}
	Let $C_n:=\{0,1,2,\ldots,n-1\}$ and $m$ a positive integer. Let's define $\mathcal{C}_n^m$ as the graph with vertices $C_n$ and edges
	$(k,k-i)\text{ mod }n,(k,k+i)\text{ mod }n$ for every $k\in C_n$ and $i\in\{1,2,\ldots,m\}$. When $m=1$, we call $\mathcal{C}_n^1$ an $n$-cycle and write $\mathcal{C}_n$.
\end{definition}

\begin{ejemplo}
		\label{4Cycle}
		For $n=4, m=1$ the resulting graph is
		
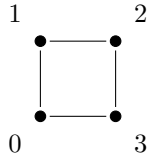
\begin{center}
			\begin{tikzpicture}
				\node (0) [label=south west:0] at (0,0) {};
				\filldraw [black] (0,0) circle (2 pt); 
				\node (1) [label=north west:1] at (0,1) {};
				\filldraw [black] (0,1) circle (2 pt); 
				\node (2) [label=north east:2] at (1,1) {};
				\filldraw [black] (1,1) circle (2 pt);
				\node (3) [label=south east:3] at (1,0) {};
				\filldraw [black] (1,0) circle (2 pt);
				\draw (3) -- (0) -- (1) -- (2);
				\draw (2) -- (3);
			\end{tikzpicture}
		\end{center}
\end{ejemplo}
For the remainder of this section, we denote by $C_{n,m}$ the semi-coarse space $(C_n,\mathcal{C}_n^m)$, where $C_n$ is
given the semi-coarse structure induced by the graph $\mathcal{C}_n^m$. 
We denote by $c$ to the map from $I_{n+1}$ onto $C_n$ such that $c(i)=i\ mod(n)$, and by $c^{-1}$ the map such that $c^{-1}(i)=c(n-i)$.

We will now define a class of functions which we can see as being endowed with an orientation. They will be particularly useful in 
the following.

\begin{definition}
	We say that $f:I_k\rightarrow C_{n,m}$ is \textbf{unidirectional} if for every $i\in\{1,\ldots,k-2\}$ we have that $f(i)$ is not a neighbor of $f(i+2)$ (i.e. $\{(f(i),f(i+2))\}\not\in \mathcal{C}_n^m$).
\end{definition}

\begin{lema}
	For every bornologous function $f:I_k\rightarrow C_n$ there exists a bornologous function $f':I_{k'}\rightarrow C_n$ such that $\langle f \rangle \simeq \langle f'\rangle $ and $f'$ is unidirectional.
\end{lema}

\begin{proof}
	Assume $f$ is not a unidirectional function, then there exists $i\in\{1,\ldots,k-2\}$ such that $f(i)$ is neighbor of $f(i+2)$. We are able to define $f_1:I_k\rightarrow C_n$ such that $f_1(j)=f(j)$ with $j\neq i+1$, $f_1(i+1)=f(i+2)$. Clearly $f_1\simeq_{sc} f$.
	
	We now define $\hat{f}_1:I_{k-1}\rightarrow C_n$ such that $\hat{f}_1(j)=f_1(j)$ if $j\leq i$, $\hat{f}_1(j)=f_1(j+1)$ if $j>i$. It
	is not difficult to show that $\langle f \rangle \simeq_{sc} \langle \hat{f}_1\rangle$.
	
	We repeat this process on $\hat{f}_1$, and so on, until one is left with a unidirectional function. Observe this is a finite process 
	since each iteration reduces the length of th domain of $\hat{f}_1$ by $1$.
\end{proof}

\begin{observacion}
	If $\lceil \frac{n}{m}\rceil \leq 2$, then no function $f:I_k\to C_n$ is unidirectional for $k>0$, because $\mathcal{C}_n^m$ is a complete graph. 
\end{observacion}

The following result shows that we are able to reduce the map $c$ to a standard form.

\begin{proposicion}
	 Let $\hat{c}:I_{\lceil \frac{n}{m}\rceil}\rightarrow C_n$ be the map $\hat{c}(k)=(k-1)m$ if $1\leq i \leq \lceil \frac{n}{m}\rceil$, and $\hat{c}(\lceil \frac{n}{m}\rceil)=0$. Then $\hat{c}\simeq_{sc} c$.
\end{proposicion}

\begin{proof}
	Define $\dot{c}:I_{n}\rightarrow C_n$ as follows
	\begin{align*}
		\dot{c}(i):=\left\lbrace\begin{array}{ll}
			0 &\ 0\leq i \leq m-1\ or\ i=n\\
			m &\ m \leq i \leq 2m-1\\
			2m &\ 2m\leq i \leq 3m-1\\
			\vdots &\ \vdots\\
			(\lceil\frac{n}{m}\rceil-1)m &\ (\lceil\frac{n}{m}\rceil-1)m\leq i \leq n-1
		\end{array} \right.
	\end{align*}
	The $\mathcal{C}_n^m$ structure gives immediately that $c\simeq_{sc}\dot{c}$. It now follows 
	that 
	$\langle c\rangle \simeq_{sc} \langle \hat{c}\rangle $
\end{proof}

Our final proposition shows us that every unidirectional function has a standard form, up to homotopy.

\begin{proposicion}
	Every unidirectional function $f:I_k \to C_n$ is homotopy equivalent to a function of one of the following forms: $1$, $c\star c\star\ldots\star c$ or $c^{-1}\star\ldots \star c^{-1}$.
\end{proposicion}

\begin{proof}
	We will assume that $\frac{n}{m}>2$, otherwise we have complete graph which is homotopy equivalent a point, and the conclusion
	of the proposition is satisfied. We also assume that $1\neq f$.
	
	Let $f:I_k\rightarrow C_n$ be unidirectional. Then $f(0)=0$, and either $f(1)\in\{1,\ldots,m\}$ or 
	$f(1)\in\{n-1,\ldots,n-m\}$. We will work with the first case. The other case is analogous.
	
	Let $\hat{f}_1:I_{k+f(1)}\to C_n$ be defined by
	\begin{align*}
		\hat{f}_1(i):=\left\lbrace\begin{array}{ll}
			i &\ 1\leq i\leq f(1)\\
			f(i-f(1)) &\ f(1)< i \leq f(1)+k
		\end{array} \right.
	\end{align*}
	We observe that $\langle f \rangle \simeq_{sc} \langle \hat{f}_1 \rangle$ 
	
	We call the transformation of $f$ into $\hat{f}_1$ a \emph{flattening of $f$ from $i=0$ to $i=1$}. We now flatten
	$\hat{f}_1$ from $f(1)$ to $f(1)+1$, and we repeat this process until we arrive at a function $\hat{f}$ such that $\hat{f}(i)-\hat{f}(i+1)=1$.
	Note that this takes a finite number of steps, and that $\langle f \rangle \simeq_{sc} \langle \hat{f} \rangle$, which proves
	the result.
\end{proof}

Combining the above two propositions, we get that the first homotopy group of $(C_n,\mathcal{C}_n^m)$ is isomorphic to a the
fundamental group of a different cyclic graph with a simpler structure, i.e.
 
\begin{corolario}
	$\pi_1^{sc}(C_n,\mathcal{C}_n^m)\cong \pi_1^{sc}(C_{\lceil \frac{n}{m}\rceil},\mathcal{C}_{\lceil \frac{n}{m}\rceil})$.
\end{corolario}

\begin{observacion}
Let $(C_n,\mathcal{C}_n)$ be the $n$-cycle and $n\in\{1,2,3\}$, then it follows that $\mathcal{C}_n$ is a connected coarse structure, precisely $\mathcal{C}_n=\mathcal{P}(C_n\times C_n)$. So, their homotopy groups are trivial.
\end{observacion}

In the rest of this subsection, we compute the fundamental group of $(C_n,\mathcal{C}_n)$.

\begin{teorema}
\label{NoTrivHom}
Let $(C_n,\mathcal{C}_n)$ be the semi-coarse space induced by the $n$-cycle graph. Then, $\pi^{sc}_1(C_n,\mathcal{C}_n)\cong \mathbb{Z}$.
\end{teorema}

Before the proof, we are going to need proving the following three lemmas, for which we define the function $p:(\mathbb{Z},\mathcal{Z})\rightarrow (C_n,\mathcal{C}_n)$ by $p(k):= k \mod n$ This function will be called \emph{projection of the integers onto $n$-cycle} or simply the \emph{projection}. We note that $p$ is a bornologous function.

\begin{lema}
\label{Lema1NoTrivHom}
Let $f:I_{k}\rightarrow C_n$, $n\geq 4$, be a bornologous function with $f(0)=0$. Then, for each $m\in p^{-1}(0)$ there exists a unique bornologous function $\hat{f}:I_{k}\rightarrow\mathbb{Z}$ such that $\hat{f}(0)=m$ and $f=p\circ\hat{f}$.
\end{lema}

\begin{proof}
Let $f:I_{k}\rightarrow C_n$ be a bornologous function with $f(0)=0$, and let $m\in p^{-1}(x_0)$.

We begin by proving that if there are two bornologous functions $\hat{f},\hat{g}:I_k\rightarrow (\mathbb{Z},\mathcal{Z})$ such that $\hat{f}(0)=\hat{g}(0)=m$ and 
$f=p\circ \hat{f}=p\circ \hat{g}$, then $\hat{f}=\hat{g}$. Given that $p\circ \hat{f}=p\circ \hat{g}$, then $\hat{f}(i)=\hat{g}(i)+n \cdot i$ for each 
$i\in\{0,\cdots,k-1\}$. Also, by hypothesis, $\hat{f}(0)=\hat{g}(0)$. Since both functions are bornologous, we have that $\hat{f}(1)=\hat{f}(0)+i_0$ 
and $\hat{g}(1)=\hat{g}(0)+j_0$ with $i_0,j_0\in\{-1,0,1\}$, and $i_0 \mod n = j_0 \mod n$, which implies that $i_0=j_0$. We conclude that $\hat{f}(1)=\hat{g}(1)$. Inductively, 
if $\hat{f}(i-1)=\hat{g}(i-1)$, then $\hat{f}(k)=\hat{g}(k)$ 
by the same argument made for $i=1$. Therefore, $\hat{f}(i)=\hat{g}(i)$ for each $i\in\{0,\cdots,m-1\}$, and $\hat{f}=\hat{g}$ as desired.

We now construct a function $\hat{f}:I_k \to \Z$ which satisfies $\hat{f}(0) = m$ and $f = p\circ \hat{f}$. First, define $\hat{f}(0)=m$, and $i\in\{1,2,\cdots,m-1\}$ we define inductively $\hat{f}(i)=\hat{f}(i-1)+i_{i-1}$ where
\begin{align*}
i_{i-1} = \left\lbrace \begin{array}{rl}
 1 & \mbox{ si } f(k) \equiv (f(k-1)+1)\ mod\ n \\
 0 & \mbox{ si } f(k) \equiv (f(k-1))\ mod\ n \\
-1 & \mbox{ si } f(k) \equiv (f(k-1)-1)\ mod\ n
\end{array}\right.
\end{align*}
So it follows by construction that $\hat{f}$ is a bornologous function, that $\hat{f}(0)=m$ and that $p\circ \hat{f}=f$.
\end{proof}

\begin{lema}
\label{Lema2NoTrivHom}
Let $f,g:I_k\rightarrow C_n$ be bornologous functions such that $f(0)=g(0)=0$ and $f(k)=g(k)=x_k$. Assume we have a homotopy $H:I_k\times \mathbb{Z}\rightarrow C_n$, between $f$ and $g$ such that $H(0,z)=x_0$ and $H(k,z)=x_k$ for each $z\in\mathbb{Z}$. Then, for each $\hat{x}_0\in p^{-1}(0)$, there exists a unique lift $\hat{H}: I_k\times \mathbb{Z}\rightarrow \mathbb{Z}$ of $H$ such that $p\circ\hat{H}=H$, $\hat{H}(0,z)=\hat{x}_0$, and there exists an $N\in \Z$ with $\hat{H}(k,N)=\hat{H}(k,z)$ for each $z\in\mathbb{Z}$.
\end{lema}

\begin{proof}
Let $f,g:I_k\rightarrow C_n$ be bornologous functions such that $f(0)=g(0)=x_0$ and $f(k)=g(k)=x_k$. Assume we have $H:I_k\times \mathbb{Z}\rightarrow C_k$, a homotopy between
$f$ and $g$. By definition, there exist integers $N,M \in \Z$, $N<M$, such
that $H(m,i)=f(m)$ for $i\leq N$ and $H(m,i)=g(m)$ for $i\geq M$. Therefore $H(0,z)=x_0$ and $H(k,z)=x_k$ for each $z\in\mathbb{Z}$. Consider a point $\hat{x}_0\in p^{-1}(x_0)$.

For each $z\in\mathbb{Z}$, we define $f_z:I_k \to C_n$ by $f(x)\coloneqq H(x,z)$ with $x\in I_k$. By \autoref{Lema1NoTrivHom}, then there are unique bornologous functions 
$\hat{f}_z:I_m\rightarrow \mathbb{Z}$ for each $z\in\mathbb{Z}$ such that $p\circ \hat{f}_z=f_z$ and $\hat{f}_z(0)=\hat{x}_0$, so we can define 
$\hat{H}(x,z)\coloneqq \hat{f}_z(x)$ for each $z\in\mathbb{Z}$ and $x\in I_k$. Under this construction, it follows that $p\circ \hat{H}=H$, $\hat{H}(0,z)=\hat{x}_0$ 
for each $z\in\mathbb{Z}$, so it remains to show that $H$ is bornologous.

Let's observe first that by construction $\hat{H}(x,z)=\hat{f}(x)$ if $z\leq N$, $\hat{H}(x,z)=\hat{g}(x)$ if $z\geq M$ and $\{(\hat{H}(x,z),\hat{H}(x+1,z))\}\in\mathcal{Z}$ for each $z\in\mathbb{Z}$ and $x\in I_k$. So, to prove $\hat{H}$ is bornologous is enough to show that
\begin{align*}
\{(\hat{H}(x,z),\hat{H}(x,z+1))\}\in \mathcal{Z}, &
\{(\hat{H}(x+1,z),\hat{H}(x+1,z+1))\}\in \mathcal{Z},\\
\{(\hat{H}(x,z),\hat{H}(x+1,z+1))\}\in \mathcal{Z}, &
\{(\hat{H}(x,z+1),\hat{H}(x+1,z))\}\in \mathcal{Z}.
\end{align*}
for each $x\in I_k$ and $z\in\{N,N+1,\cdots,M-2,M-1\}$.

Let $x \in I_k$ and $z\in\{N,N+1,\cdots,M-2,M-1\}$. For each of the pairs $(x,z),(x',z')$ considered, above, since $H$
is bornologous, we have $p \circ \hat{H}(x,z) = p\circ \hat{H}(x',z') + i_{x,z,x',z'}$, where $i_{x,z,x',z'}\in \{-1,0,1\}$. Then $H(x,z) = H(x',z') + 
\hat{i}_{x,z,x',z'}$, where $\hat{i}_{x,z,x',z'} \in \{-1,0,1\} \mod n$. For $x=0$, note that $\hat{f}_z(0)=\hat{f}_{z+1}(0)=\hat{x}_0$, and therefore $\hat{i}_{0,z,0,z'}=0$ for any
$z\in \Z$ and $z'\in \{z,z+1\}$. Since $\hat{f}_z(x) = \hat{f}_z(x+1) + \hat{i}_{x,z,x+1,z}$ by construction, we also have that $\hat{i}_{x,z,x+1,z}=i_{x,z,x'z}\in \{-1,0,1\}$ for all $x\in I_k$. It follows that, for $x,x' \in \{0,1\}$, any $z$, and $z'\in\{z,z+1\}$, we have $\hat{i}_{x,z,x',z'} \in \{-1,0,1\}$.
We use the same argument and induction to show that $\hat{i}_{x,z,x',z'}=i_{x,z,x'z'}\in \{-1,0,1\}$ for any $x,x'\in I_k$, $z\in \Z,$ and $z'\in \{z,z+1\}$. It follows that $\hat{H}$
is bornologous.
\end{proof}

\begin{observacion}
The homotopy is unique because $\hat{f}_z$ is unique for each $z\in\mathbb{Z}$. If we change any element, the relations $p\circ \hat{H}=H$, $\hat{f}_z(0)=\hat{x}_0$ will no longer be true.
\end{observacion}

\begin{lema}
\label{Lema3NoTrivHom}
Let $f:I_p\rightarrow \mathbb{Z}$ and $g:I_{q}\rightarrow \mathbb{Z}$ be bornologous functions such that $f(0)=g(0)$ and $f(p)=g(q)$. Then $\langle f \rangle \simeq_{sc} \langle g \rangle$.
\end{lema}

\begin{proof}
Let $f:I_p\rightarrow \mathbb{Z}$ and $g:I_{q}\rightarrow \mathbb{Z}$ be bornologous functions such that $f(0)=g(0)$ and $f(p)=g(q)$.

Observe that there exist unidirectional functions $f':I_{p'}\rightarrow \mathbb{Z}$, $f' \in \langle f \rangle$, and $g':I_{q'}\rightarrow \mathbb{Z}$, $g'\in \langle g \rangle$. If $f(0)=f(p)$, then $m'=n'=0$, and the result follows.

Now assume that $f(0)\neq f(p)$. Then $f'$ must be the function $f':I_{|f(p)-f(0)|}\rightarrow \mathbb{Z}$ such that $f'(i)=f(0)+sign(f(m)-f(0))i$ for every $i\in \{0,\ldots,|f(m)-f(0)|\}$, with $sign(\cdot)$ the sign function. Analogously for $g'$. However, by hypothesis, $f(p)-f(0)=g(q)-g(0)$, thus $f'=g'$. Therefore $\langle f \rangle =\langle g \rangle$.
\end{proof}

\begin{observacion}
Let $f,g:I_k\rightarrow X$ be bornologous functions such that $f(k)=g(0)$. By \autoref{BornologousStar}, $f\star g$ is bornologous.
\end{observacion}

\begin{proof}[\autoref{NoTrivHom}] Let $(C_n,\mathcal{C}_n)$ be the semi-coarse space induced by the $n$-cycle graph. We define $\phi:\mathbb{Z}\rightarrow \pi_1^{sc}(C_n,\mathcal{C}_n)$ by $\phi(z)\coloneqq [\langle c^{\star z}\rangle ]$ Note that, by construction, $c^{\star z}$ is a bornologous function for each $z\in\mathbb{Z}$. It is clear from the definition and 
the definition of the $\star$ operation that $\phi$ is a group homomorphism.

To prove that $\phi$ is surjective, consider $[\langle f \rangle]\in\pi_1^{sc}(C_n,\mathcal{C}_n)$, and, without loss of generality, assume that $f\in \langle f \rangle$ is
a bornologous function $f:I_{k}\rightarrow C_n$ with $f(0)=f(k)=0$, and note that $0\in p^{-1}(\{0\})=n\mathbb{Z}$. Thus, by \autoref{Lema1NoTrivHom} there is a unique
bornologous function $\hat{f}$ such that $\hat{f}(0)=0$ and $f=p\circ \hat{f}$. Since $(p\circ \hat{f})(k)=f(k)=0$, note that $\hat{f}(k)=nq$ for some $q\in\mathbb{Z}$. Thus, $\hat{f}\simeq_{sc}\widehat{c^{\star q}}$ by \autoref{Lema3NoTrivHom},
and therefore
\begin{align*}
f=p\circ \hat{f} \cong_{sc} p\circ \widehat{c^{\star n}} = c^{\star n}
\end{align*}
We conclude that $[\langle f\rangle]=[\langle \widehat{c^{\star q}}\rangle]=\phi(q)$, and therefore $\phi$ is surjective.

To prove that $\phi$ is injective, suppose that $\phi(p)=\phi(q)$ for some $p,q \in \Z$, so we have that $\langle c^{\star q}\rangle\simeq_{sc}\langle c^{\star p} \rangle$. Let $H:I_k \to
C_n$ be a homotopy from $c^{\star q}$ to $c^{\star p}$, where $k = |pqn|$, and, abusing notation, we understand $c^{\star p}$ and $c^{\star q}$ here to be the extension of each map to $I_{k}$ by defining $c^{\star p}(i) = 0$ if $i>|pn|$, and similarly for $c^{\star q}$. Note that $H(0,i)=0\in C_n$ for each $i\in\mathbb{Z}$ by definition. We also note that $0\in p^{-1}(\{0\})=n\mathbb{Z}$.

By \autoref{Lema2NoTrivHom}, there is a unique homotopy $\hat{H}:I_{pqn}\to C_n$ such that $H=p\circ \hat{H}$ and $\hat{H}(0,n)=0$ for each $n\in\mathbb{Z}$. By uniqueness of the
homotopy, we have that $\widehat{c^{\star p}}(i)=\hat{H}(i,N)$ and $\widehat{c^{\star q}}(i)=\hat{H}(i,M)$. Furthermore, since $\widehat{c^{\star q}}(|nq|)=\widehat{c^{\star p}}(k)=\hat{H}(k,N)=\hat{H}(k,M)=\widehat{c^{\star q }}(k)=\widehat{c^{\star q}}(|nq|)$, we obtain $p=q$, and we conclude that the function is injective.
\end{proof}

\subsection{Long Exact Sequence in Homotopy}

With the same goal, we will do the following definitions which will help us to prove the fundamental result: the long exact sequence in relative homotopy.

\begin{definicion}
\label{retraccion}
Let $(X,\mathcal{V})$ be a semi-coarse space and $(A,\mathcal{V}_A)$ be a semi-coarse subspace. Then,
\begin{itemize}
\item[•] The bornologous function $r:X\rightarrow A$ is a retraction if satisfies that $r\circ i=id_A$, where $i:A\rightarrow X$ such that $i(a)=a$ for each $a\in A$.
\item[•] A bornologous function $F:(X\times \mathbb{Z},\mathcal{V}\times\mathcal{Z})\rightarrow (X,\mathcal{V})$ is a deformation retract of $X$ onto $A$ if there are $N<0<M$ such that, for each $x\in X$ and $a\in A$, $F(x,z)=x$ if $z\leq N$, $F(x,z)\in A$ and $F(a,z)=a$ if $z\geq M$.
\item[•] If $F$ also satisfies that $F(a,z)=a$ for each $z\in\mathbb{Z}$, then $F$ is called strong deformation retraction.
\end{itemize}
\end{definicion}

\begin{lema}
\label{CompReng}
Let $f:I_m^n\rightarrow X$ be a bornologous function, $k\in I_m^n$ and $w\in\{1,\ldots,n\}$. Then $f\simeq g\simeq h$ such that
\begin{align*}
g(y):= \left\lbrace\begin{array}{rl}
g(y) & \mbox{if } y_w\leq k_w\\
g(y-(y_w-k_w)e_w) & \mbox{if }y_w>k_w
\end{array}\right.
\end{align*}
and
\begin{align*}
h(y):= \left\lbrace\begin{array}{rl}
h(y) & \mbox{if } y_w\geq k_w\\
h(y-(y_w-k_w)e_w) & \mbox{if }y_w<k_w
\end{array}\right.
\end{align*}
\end{lema} 

\begin{proof}
Let $f:I_m^n\rightarrow X$ be a bornologous function, $k\in I_m^n$, $w\in\{1,\ldots,n\}$,
\begin{align*}
g(y):= \left\lbrace\begin{array}{rl}
g(y) & \mbox{if } y_w\leq k_w\\
g(y-(y_w-k_w)e_w) & \mbox{if }y_w>k_w
\end{array}\right.
\end{align*}
and
\begin{align*}
h(y):= \left\lbrace\begin{array}{rl}
h(y) & \mbox{if } y_w\geq k_w\\
h(y-(y_w-k_w)e_w) & \mbox{if }y_w<k_w
\end{array}\right.
\end{align*}
We will just prove for $g$, the proof for $h$ is analogous. We are going to proceed by math induction:

\textbf{Base case:} Let $k$ such that $k_w=m-1$, $f\simeq g$ such that 
\begin{align*}
g(y)=\left\lbrace\begin{array}{rl}
g(y) & \mbox{if } y_w\leq m-1\\
g(y-e_w) & \mbox{if } y_w=m
\end{array}\right. .
\end{align*}

\textbf{Induction step:} Assume that $f\simeq g_i$ for every $i\in\{1,\ldots,k_w-1\}$ where
\begin{align*}
g_i(y):= \left\lbrace\begin{array}{rl}
g_i(y) & \mbox{if } y_w\leq i\\
g_i(y-(y_w-i)e_w) & \mbox{if }y_w>i
\end{array}\right.
\end{align*}
once again $f\simeq g$ such that 
\begin{align*}
g(y):= \left\lbrace\begin{array}{rl}
g(y) & \mbox{if } y_w\leq k_w\\
g(y-(y_w-k_w)e_w) & \mbox{if }y_w>k_w
\end{array}\right.
\end{align*}
\end{proof}

Observe that the last lemma makes us being able to replace either the plates ``above''  or ``below'' $(I_m^n)_w^{k_w}$. Thereby, the following lemma is a corollary. With the intention of making the proof clear, allow us call $g$ as $(f)_w^\uparrow$ and $h$ as $(f)_w^\downarrow$.

\begin{lema}
\label{StrongDef}
Let $n$ be a natural number and $m$ be a non-negative integer number, then there is a strong deformation retraction $I^n_m$ onto $\{*\}$.
\end{lema}

\begin{proof}
Let $n$ be natural number and $m$ be a non-negative integer number, let's consider $*\in I^n_m$, then $*:=k=(k_1,\cdots, k_n)$ where $k_i\in\{0,1,\cdots,m\}$ fixed for each $i\in\{1,\cdots,n\}$. Denote $1_{I_m^n}$ the identity function, then $1_{I_m^n}\simeq(1_{I_m^n})_1^\rightarrow \simeq ((1_{I_m^n})_1^\rightarrow)_1^\leftarrow)=:g_1$ by \autoref{CompReng}. Therefore, $g_i:= ((g_{i-1})_i^\rightarrow)_i^\leftarrow)\simeq g_{i-1}$ for each $i\{2,\ldots,n\}$ by \autoref{CompReng}. Thus, $g_n\simeq f$. The strong deformation that we are looking for is precisely that homotopy.
\end{proof}

The following lemma is a tool to prove the Compression Criterion, which one is our goal to finally get the long exact sequence in our homotopy.

\begin{lema}
\label{PreCritComp}
Let $f,g: (I_m^n,\partial I_m^n, J_m^{n-1})\rightarrow (X,A,*)$ such that $f\simeq_{sc} g\ rel(\partial I^n)$. Let $H$ the homotopy with $N<M$ and $N'>0$, then $ \langle f \rangle \simeq_{sc} \langle h_{N'} \rangle\ rel(\partial I^n)$ with $h_{N'}:(I_{m+N'}^n,\partial I_{m+N'}^n, J_{m+N'}^{n-1}\rightarrow (X,A,*)$ such that
\begin{align*}
h_{N'}(x):=\left\lbrace\begin{array}{rl}
H(x-N'e_n,N+N') & x\in I_m^n+N'e_n \\
H(x-ke_n,N+k) & x\in (I_m^n)_n^0+ke_n \mbox{ and } k\in \{0,\ldots,N'-1\}\\
* & \mbox{otherwise}
\end{array}\right.
\end{align*}
\end{lema}
where $(I_m^n)_n^0\subset I_m^n\subset I_{m+N'}^n$.

\begin{proof}
Let $f,g: (I_m^n,\partial I_m^n, J_m^{n-1})\rightarrow (X,A,*)$ such that $f\simeq_{sc} g\ rel(\partial I^n)$. Let $H$ the homotopy with $N<M$. We are going to do the proof by mathematical induction:

\textbf{Base case:} 
We have that $\langle f \rangle \simeq_{sc} \langle g_1 \rangle\ rel(\partial I^n)$ where $g_1:(I_{m+1}^n,\partial I_{m+1}^n,J_{m+1}^{n-1})\rightarrow (X,A,*)$ 
\begin{align*}
g_1(x)=\left\lbrace\begin{array}{rl}
f(x-e_n) & x\in I_m^n+e_n\\
f(x)     & x\in (I_m)_n^0\\
*        & \mbox{otherwise}
\end{array}\right.
\end{align*}
Since $H$ is a bornologous function, we are able to replace $f(x-e_n)$ by $H(x-e_n,N+1)$, calling that function $h_1$, and we get $g\simeq_{sc} h_1\ rel(\partial I_{m+1}^n)$. Note that $f(x)=H(x,N')$ when $x\in (I_m)_n^0$. So, $\langle f \rangle \simeq_{sc} \langle h_1\rangle\ rel(\partial I^n)$.

\textbf{Induction step:} Suppose that for $w>0$ we have $\langle f \rangle \simeq_{sc} \langle h_w\rangle\ rel(\partial I^n)$ with $h_{N'}:(I_{m+w}^n,\partial I_{m+w}^n, J_{m+w}^{n-1}\rightarrow (X,A,*)$ such that
\begin{align*}
h_{w}(x):=\left\lbrace\begin{array}{rl}
H(x-we_n,N+w) & x\in I_m^n+we_n \\
H(x-ke_n,N+k) & x\in (I_m^n)_n^0+ke_n \mbox{ and } k\in \{0,\ldots,w-1\}\\
* & \mbox{otherwise}
\end{array}\right.
\end{align*}
We have that $\langle h_w\rangle \simeq_{sc}\langle g_{w+1}\rangle\ rel(\partial I^n)$ where $g_{w+1}:(I_{m+w+1}^{n},\partial I_{m+w+1}^{n}, J_{m+w+1}^{n-1})\rightarrow (X,A,*)$ such that
\begin{align*}
g_{w+1}(x):=\left\lbrace\begin{array}{rl}
H(x-(w+1)e_n,N+w) & x\in I_m^n+(w+1)e_n \\
H(x-we_n,N+w) & x\in (I_m^n)_n^0+we_n\\
H(x-ke_n,N+k) & x\in (I_m^n)_n^0+ke_n \mbox{ and } k\in \{0,\ldots,w-1\}\\
* & \mbox{otherwise}
\end{array}\right.
\end{align*}
Since $H$ is a bornologous function, we are able to replace $H(x-(w+1)e_n,N+W)$ by $H(x-(w+1),N+w+1)$, calling that function $h_{w+1}$, and we get $g_{w+1}\simeq_{sc} h_{w+1}\ rel(\partial I_{m+w+1}^n)$.
\end{proof}

\begin{lema}[Compression Criterion]
\label{CritComp}
Let $(X,\mathcal{V})$ be a semi-coarse space and $*\in A\subset X$. Then a function $f:(I^n,\partial I^n,J^{n-1})\rightarrow (X,A,*)$ represents the zero in $\pi_n^{sc}((X,\mathcal{V}),A,*)$ if, and only if, is relative homotopic in $\partial I^n$ to a function with its image contained in $A$.
\end{lema}

\begin{proof}
Let $(X,\mathcal{V})$ be a semi-coarse space, $*\in A\subset X$ and $f:(I^n,\partial I^n,J^{n-1})\rightarrow (X,A,*)$.
\begin{itemize}
\item[($\Rightarrow$)] If $[f]=0$, then there exist $m\in\mathbb{N}$ with $f:(I_m^n,\partial I_m^n,J^{n-1})\rightarrow (X,A,*)$ such that $f\simeq_{sc} *$. Let $H$ be the homotopy with $N<M$, then by \autoref{PreCritComp}, we get $h_{M-N}$ as we write in such lemma which has its image contained in $A$ and $\langle f \rangle \simeq_{sc} \langle h_{M-N}\rangle\ rel(\partial I^n)$.

\item[($\Leftarrow$)] If $f$ is relative homotopic in $\partial I^n$ to a function with its image contained in $A$, let's say $g$. As we do in \autoref{StrongDef}, with $w=(0,\ldots,0,m)$, then $g\simeq_{sc} (g)_n^\downarrow=*$. Getting what we wanted.
\end{itemize}
\end{proof}

\begin{definicion}
Let $(X,\mathcal{V})$ and $(Y,\mathcal{W})$ be semi-coarse spaces and $f:((X,\mathcal{V}),A,*)\rightarrow ((Y.\mathcal{W}),B,*)$ be a Bornologous function. Then, we define $f_*:\pi^{sc}_*((X,\mathcal{V}),A,*)\rightarrow \pi^{sc}_*((Y,\mathcal{W}),B,*)$ as $f_*([\langle h \rangle])=[\langle f\circ h \rangle]$ for each $[\langle h \rangle]\in \pi^{sc}_*((X,\mathcal{V}),A,*)$.
\end{definicion}

\begin{observacion}
Let $(X,\mathcal{V})$ and $(Y,\mathcal{W})$ be semi-coarse spaces and $f:((X,\mathcal{V}),A,*)\rightarrow ((Y.\mathcal{W}),B,*)$ be a (based) Bornologous function. Then $f_*$ is well-defined and is an homomorphism. It is well-defined because, if we take $\langle g\rangle\in [\langle h\rangle]\in\pi^{sc}_*((X,\mathcal{V}),A,*)$, then there is a homotopy $H$ between them, so $f\circ H$ is a homotopy between $\langle f\circ h\rangle$ and $\langle f\circ g\rangle$. If $[\langle h \rangle],[\langle g \rangle]\in \pi^{sc}_*((X,\mathcal{V}),A,*)$. It is a homomorphism because $f\circ (g\star h)=f\circ g(k_1,\cdots k_n)$ if $0\leq k_i \leq m$ for each $i\in\{1,\cdots,n\}$, $f\circ (g\star h)=f\circ h(k_1,\cdots k_n)$ if $m<k_1\leq 2m$ $0\leq k_i \leq m$ for each $i\in\{2,\cdots,n\}$, and $f\circ (g\star h)=f(*)$ anywhere else, which is the same as $(f\circ g) \star (f\circ h)$.
\end{observacion}

\begin{teorema}[Long Exact Sequence in Homotopy]
\label{SecLarExact}
Let $(X,\mathcal{V})$ be a semi-coarse space and $*\in A\subset X$. As well $i:(A,*)\rightarrow (X,*)$, $j:(X,*)\rightarrow (X,A)$ inclusions. Then there is a homomorphism $\partial_n:\pi^{sc}_n((X,\mathcal{V}),A)\rightarrow \pi^{sc}_n(A,*)$ such that the long sequence

\[
\xymatrix@1{
\cdots \ar[r]^{j_*} & \pi^{sc}_{n+1}(X,A) \ar[r]^{\partial_{n+1}} & \pi^{sc}_{n}(A,*) \ar[r]^{i_*} & \pi^{sc}_n(X,*) \ar[r]^{j_*} & \cdots\\
\cdots \ar[r]^{\partial_2} & \pi^{sc}_{1}(A,*) \ar[r]^{i_*} & \pi^{sc}_1(X,*) \ar[r]^{j_*} & \pi^{sc}_1(X,A)
}
\]

is exact.
\end{teorema}

\begin{proof}
Let $(X,\mathcal{V})$ be a semi-coarse space and $*\in A\subset X$. As well $i:(A,*)\rightarrow (X,*)$, $j:(X,*)\rightarrow (X,A)$ inclusions. We are going to define $\partial_n:\pi^{sc}_n((X,\mathcal{V}),A)\rightarrow \pi^{sc}_n(A,*)$ as $f\mapsto f|_{I^{n-1}\times \{0\}}$, which is well-defined and is a homomorphism.

Now we will prove that the long sequence is exact:
\begin{itemize}
\item[•] \textbf{($Im(i_*)\subset Ker(j_*)$)} If $n\geq 1$ and $f:(I^n,\partial I^n, J^{n-1})\rightarrow (A,*,*)\in\pi_n^{sc}(A,*)$, then $j\circ i\circ f\in\pi_n^{sc}(X,A)$, so by compression criterion we have that $[\langle f \rangle]=[\langle * \rangle]$.
\item[•] \textbf{($Im(i_*)\supset Ker(j_*)$)} If $n\geq 1$ and $[\langle f \rangle]\in \pi_n^{sc}((X,\mathcal{V}),*)$ such that $j_*[\langle f \rangle]=[\langle *\rangle]$. Then by compression criterion, $\langle f \rangle \simeq_{sc} \langle g \rangle\ rel(\partial I^n)$, where the image of $g$ is contained in $A$, thus $[\langle g \rangle]=[\langle f \rangle]\in\pi_n^{sc}((X,\mathcal{V}),*)$ is in the image of $i_*$.
\item[•] \textbf{($Im(j_*)\subset Ker(\partial)$)} If $n\geq 2$ and $[\langle f \rangle]\in \pi_n^{sc}((X,\mathcal{V}),*)$, then $\partial j\langle f \rangle=\langle f|_{I^{n-1}\times \{0\}} \rangle= \langle * \rangle$, so $\partial j_*[\langle f \rangle]=[\langle * \rangle]$.
\item[•] \textbf{($Im(j_*)\supset Ker(\partial)$)} If $n\geq 2$ and $[\langle f \rangle]\in\pi_n^{sc}((X,\mathcal{V}),A)$ such that $\partial[\langle f \rangle]=[\langle*\rangle]$. Then $f|_{I^{n-1}\times\{0\}}$ is homotopic to $\langle * \rangle$ by the compression theorem through a homotopy $H:I^{n-1}_m\times \mathbb{Z}\rightarrow A\ rel (\partial I_m^{n-1})$, $N<0<M$ such that $H(x,z)=*$ if $z\leq N$ and $H(x,z)=f|_{I^{n-1}\times\{0\}}(x,z)$ if $z\geq M$. Let's define $g:(I_{m+M-N}^n,\partial I_{m+M-N}^n, J_{m+M-N}^{n-1})\rightarrow ((X,\mathcal{V}),A,*)$ such that
\begin{align*}
g(k_1,\cdots,k_n)= H(k_1,\cdots,k_{n-1},z=k_n)
\end{align*}
with $0\leq k_i\leq m$ when $i\in\{1,\cdots,n-1\}$ and $0\leq k_n\leq M-N$,
\begin{align*}
g(k_1,\cdots,k_n)= f(k_1,\cdots,k_{n-1},k_n-M+N)
\end{align*}
with $0\leq k_i\leq m$ when $i\in\{1,\cdots,n-1\}$ and $M-N+1\leq k_n\leq M-N+m$, and $g(k_1,\cdots,k_n)$ anywhere else. So, we note that $g(\partial I_{m+M-N}^n)=*$, that is, $[\langle g \rangle]\in Im(j_*)$, and $[\langle f \rangle]=[\langle g\rangle]$.
\item[•] \textbf{($Im(\partial)\subset Ker(i_*)$)} If $n\geq 2$ and $f:(I_m^n,\partial I_m^n,J_m^{n-1})\rightarrow((X,\mathcal{V}),A,*)$. Then, defining $H:I_m^{n-1}\times Z\rightarrow (X,\mathcal{V})$ as $H(k_1,\cdots,k_n)=f(k_1,\cdots,k_n)$ with $0\leq k_n \leq m$, $H(k_1,\cdots,z)=f(k_1,\cdots,k_{n-1},0)$ if $z\leq 0$, and $H(k_1,\cdots,z)=f(k_1,\cdots,k_{n-1},m)$ if $z\geq m$, we get that $f|_{I_m^{n-1}\times \{0\}}$ is relative homotopic to $*$ through $H$. So, $i_*\partial[\langle f \rangle]=[\langle * \rangle]$ by \autoref{CritComp}.
\item[•] \textbf{($Im(\partial)\supset Ker(i_*)$)} If $n\geq 2$ and $f:(I_m^n,\partial I_m^n,J_m^{n-1})\rightarrow(A,*,*)$ such that $i_*([\langle f \rangle])=[\langle*\rangle]$, then we have a homotpy between $f$ and $*$, which gives us a function $F:(I_{M-N}^{n+1},\partial I_{M-N}^{n+1}, J_{M-N}^{n})\rightarrow((X,\mathcal{V}),A,*)$ such that $\partial[\langle F\rangle]=[\langle f \rangle]$.
\end{itemize}
\end{proof}

\section{Homology}
\label{sec:Homology}

In the last section, we will construct a Vietoris-Rips homology for semi-coarse spaces and show that it is homotopy invariant. We finish with the observation that the homology defined here for a finite semi-coarse space is isomorphic to the homology of the clique complex of a finite graph.

\subsection{Simplicial Homology}

\begin{definicion}
\label{Simplejo}
Let $(X,\mathcal{V})$ be a semi-coarse and $E$ be a controlled set, and let $\mathcal{R}$ be a relation on $X$ such that $x\mathcal{R}y$ iff $(x,y)\in E$. We now define the following
\begin{itemize}
\item $\Sigma_{E}^{(0)}:=\{\{x\}:x\in X\}$,
\item $\Sigma_{E}^{(n)}:=\{\{x_0,\cdots,x_n\} \subset X\mid \forall i,j\in\{0,\cdots,n\} x_i\mathcal{R}x_j \text{ and } ((i\neq j)\implies x_i\neq x_j) \}$
\item We call the collection $\cup_{i=0}^n \Sigma_{E}^{(i)}$ the \emph{$n$-skeleton} of $\Sigma$.
\item $\Sigma_{E} := \cup_n \Sigma_{E}^{(n)}$ will be called \emph{the simplicial complex associated to $E$}.
\end{itemize}
\end{definicion}

\begin{observacion}
We note that $\Sigma_{E}$ satisfies the definition of simplicial complex for any element $E\in\mathcal{V}$ in the semi-coarse structure $X$, \autoref{Simplejo}, since $\Sigma_{E}^{(0)}$ contains all sets with one vertex and every subset of a simplex is a simplex having all their elements related.
\end{observacion}

\begin{definicion}
We define $C_q(X,E)$ to be the free abelian groups generated by ordered simplicial chains, denoted $[v_0,v_1,\dots,v_q]$, where $[v_0,v_1,\dots,v_q]=0$ if the vertices are not all pairwise different. We define the differential $\partial_q$ by \begin{equation*}
	\partial_q[v_0,v_1,\dots,v_q]:=\sum_{i=0}^q(-1)^i[v_0,\dots,\hat{v}_i,v_q]
\end{equation*}
and we denote the chain complex by $C_*(X,E):=\{C_{q}(X,E),\partial_q\}$ and the resulting homology groups by  $H_*(X,E)$.
\end{definicion}

The above definitions show how to construct homology groups from a single element of a semi-coarse structure $E\in\mathcal{V}$. The next
lemma will allow us to construct a directed system from these homology groups.
\begin{lema}
\label{LemaIncHom}
Let $(X,\mathcal{V})$ be a semi-coarse space and $E,E'$ controlled sets such that $E\subset E'$. Then there exists a homomorphism $i_*:H(X,E)\rightarrow H(X,E')$.
\end{lema}

\begin{proof}
Let $(X,\mathcal{V})$ be a semi-coarse space and let $E,E'$ be controlled sets such that $E\subset E'$. If $i:E\rightarrow E'$ is the inclusion from $E$ to $E'$, then for a generator $\sigma=[v_0,\dots,v_q]$ of $C_q(X,E)$, we define $(i_\#)_q:C_{q}(X,E)\rightarrow C_{q}(X,E')$ by
\[(i_\#)_q[x_0,\cdots,x_q]:=[i(x_0),\cdots,i(x_q)]=[x_0,\cdots,x_q]\in C_q(X,E').\] 
We then extend this by linearity. Note that the inclusion $E\subset E'$ ensures that $C_q(X,E)\subset C_q(X,E')$, since $x\mathcal{R}_E y \implies x\mathcal{R}_{E'}y$ for any $x,y\in X$. The inclusion $i_\#$ also satisfies $\partial_{E'}i_\#=i_\#\partial_{E}$, so $i_\#$ is a chain map. The induced map $i_*:H_*(X,E)\rightarrow H_*(X,E')$ such that $i_*[\sigma]=[i_\#(\sigma)]=[\sigma]\in H_*(X,E')$ is the desired map in homology.
\end{proof}

\begin{definicion}
\label{LimDirHomom}
Let $(X,\mathcal{V})$ be a semi-coarse space, and consider $\mathcal{V}$ with the partial order given by inclusion of sets. We define \[H_*(X,\mathcal{V}):=\lim\limits_\rightarrow\{H(X,E),\pi_E^{E'},\mathcal{V}\},\]
where $\pi_E^{E'}=i_*:H(X,E)\rightarrow H(X,E')$ from the last lemma. 
We call $H_*(X,\mathcal{V})$ the \emph{homology of the semi-coarse space $(X,\mathcal{V})$}. We will sometimes refer to $H_*(X,\mathcal{V})$ as the
\emph{Vietoris-Rips homology of $(X,\mathcal{V})$.}
\end{definicion}

We can note that the set of symmetric sets $E\in \mathcal{V}$ is cofinal in $\mathcal{V}$, i.e. for every $E\subset \mathcal{V}$, we have
$E \subset E\cup E^{-1}\cup \Delta_X \in \mathcal{V}$.  Denote by $\mathcal{V}_\mathcal{S}$ the collection $\{E\in\mathcal{V} \mid E=E^{-1}\}$. 
By the cofinality We have that
\[
\lim_\rightarrow \{ H(X,E),\pi_E^{E'},\mathcal{V}\}\cong \lim_\rightarrow \{H(X,E),\pi_E^{E'},\mathcal{V}_\mathcal{S}\},
\]
so it is enough to consider symmetric elements of $\mathcal{V}$ when constructing the homology of $(X,\mathcal{V})$.

We will now show the semi-coarse homology is a covariant functor. We will start by showing that homology is functorial for a fixed controlled set
$E\in\mathcal{V}$.

\begin{lema}
\label{FuntorialidadControlado}
Let $(X,\mathcal{V})$ and $(X,\mathcal{W})$ be semi-coarse spaces, $E$ controlled by $\mathcal{V}$ and $f:(X,\mathcal{V})\rightarrow (Y,\mathcal{W})$ be a bornologous function. Then,
\begin{itemize}
\item[(i)] $f_\#:C(X,E)\rightarrow C(Y,(f\times f)(E))$ where \[(f_\#)_n[x_0,\cdots,x_n]=[f(x_0),\cdots,f(x_n)]\] is a chain map.
\item[(ii)] $f_*:H(X,E)\rightarrow H(Y,(f\times f)(E))$ is a group homomorphism, where $(f_*)_n[\sigma^n]=[(f_\#)_n\sigma^n]$.
\end{itemize}
\end{lema}

\begin{proof}
Let $(X,\mathcal{V})$ and $(Y,\mathcal{W})$ be semi-coarse spaces, $E$ controlled by $\mathcal{V}$ and 
$f:(X,\mathcal{V})\rightarrow (Y,\mathcal{W})$ be a bornologous function. Since $f$ is bornologous, $\{(f(x),f(y))\}\in(f\times f)(E)\in\mathbb{W}$ 
if $\{(x,y)\}\in E$, so $f$ induces the simplicial map $f_\#$ defined by
\begin{align*}
(f_\#)_n[x_0,\cdots,x_n]= & [f(x_0),\cdots,f(x_n)].
\end{align*}
Since $f_\#$ is a chain map, it induces the map $f:_*:H_*(X,E)\to H_*(Y,(f\times f)(E))$ on homology by
\begin{align*}(f_*)_n[\sigma^n]= & [(f_\#)_n\sigma^n].\qedhere
\end{align*}
\end{proof}

\begin{observacion}
The previous result is also true for $f_\#:C(X,E)\rightarrow C(Y,A)$ and $f_*:H(X,E)\rightarrow H(Y,A)$ such that $(f\times f)(E) \subset A$ is controlled by $\mathcal{W}$. This follows from the previous lemma and \autoref{LemaIncHom}.
\end{observacion}

\begin{teorema}
\label{TeoFuntqc}
Semi-coarse homology is a covariant functor $H_*:\catname{SCoarse}\to\catname{Ab}$.
\end{teorema}

\begin{proof}
Let $(X,\mathcal{V})$ and $(Y,\mathcal{W})$ be semi-coarse spaces and let $f:(X,\mathcal{V})\rightarrow(Y,\mathcal{W})$ be a bornologous function. Let $\{H(X,E),\pi_E^{E'},\mathcal{V}\}$ and $\{H(Y,E),\pi_E^{E'},\mathcal{W}\}$ be the directed systems of the homology groups, where the directed sets $\mathcal{V}$ and $\mathcal{W}$ are partially ordered by inclusion.

Note that, if $A\subset B\subset X\times X$, then $(f\times f)(A)\subset (f\times f)(B)$ so $(f\times f)$ is preserves the preorders on $\mathcal{V}$ and $\mathcal{W}$. Therefore, for each $E,E'\in\mathcal{V}$ such that $E\subset E'$, the diagram
\[\xymatrix{
H(X,E) \ar[r]^{f_*} \ar[d]_{\pi_{E}^{E'}} & H(Y,(f\times f)(E)) \ar[d]_{\pi_{(f\times f)(E)}^{(f\times f)(E')}} \\
H(X,E') \ar[r]^{f_*} & H(Y,(f\times f)(E'))
}\]
commutes. It follows by definition that there exists a homomorphism $F:H_*(X,\mathcal{V})\to H_*(X,\mathcal{W})$ between the directed limits, that is,
\begin{align*}
F:\lim_\rightarrow\{H(X,E),\pi_E^{E'},\mathcal{V}\}\rightarrow\lim_\rightarrow\{H(Y,E),\pi_E^{E'},(f\times f)(\mathcal{V})\}.
\end{align*}
\end{proof}

We now show that homotopic maps induce the same map on semi-coarse homology

\begin{teorema}
\label{EquivHom}
Let $(X,\mathcal{V})$ and $(Y,\mathcal{W})$ be semi-coarse spaces. If $f,g:(X,\mathcal{V})\rightarrow(Y,\mathcal{W})$ are homotopic functions, then the induced homomorphisms $f_*,g_*:H_*(X,\mathcal{V})\rightarrow H_*(Y,\mathcal{W})$ are equal.
\end{teorema}

\begin{proof}
Let $(X,\mathcal{V})$ and $(Y,\mathcal{W})$ be semi-coarse spaces. If $f,g:(X,\mathcal{V})\rightarrow(Y,\mathcal{W})$ are homotopic functions, then there are a bornologous function
\begin{align*}
H:(X\times\mathbb{Z},\mathcal{V}\times \mathcal{Z})\rightarrow (Y,\mathcal{W})
\end{align*}
and $N<0<M$ integer numbers such that $H(x,z)=f(x)$ if $z\leq N,x\in X$ and $H(x,z)=g(x)$ if $z\geq M,x\in X$. We will denote by $h_{z}(x)\coloneqq H(x,z)$ for $x\in X$ and $z\in\{N+1,N+2,\cdots,M-2,M-1\}$.

We will observe what happens with the induced homomorphisms by $f$ and $h_{N+1}$. Then, for each $E$ controlled by $\mathcal{V}$, we define $(\Psi_E)_q:C_q(X,E)\rightarrow C(Y,H(E,\mathbb{Z}))$ such that
\begin{align*}
(\Psi_E)_q[x_0,\cdots,x_q]\coloneqq \sum_{i=0}^q(-1)^i[f(x_0),\cdots,f(x_i),h_{N+1}(x_i),\cdots, h_{N+1}(x_q)].
\end{align*}
A straightforward calculation shows that $\Psi_E$ is a chain homotopy between the induced chain maps $f_\#$ and $(h_{N+1})_\#$, so $f_*=(h_{N+1})_*:H(C(X,E))\rightarrow H(C(Y,H(E,\mathbb{Z})))$. Thereby, by similar arguments to the previous theorem, we get that
\begin{align*}
f_*=(h_{N+1})_*:H(C(X,\mathcal{V}))\rightarrow H(C(Y,\mathcal{W})).
\end{align*}
We will repeat the same for $h_i$ and $h_{i+1}$ with $i\in{N+1,\cdots M-1}$, finally arriving at
\begin{align*}
f_*=g_*:H(C(X,\mathcal{V}))\rightarrow H(C(Y,\mathcal{W})),
\end{align*}
as desired.
\end{proof}

\subsection{Graphs and Semi-Coarse Spaces}

Now that we have introduced the basic concepts of semi-coarse homology, we recall every graph is a roofed semi-coarse space (where roofed semi-coarse space is defined in \autoref{DefRoof}), and, given the similarities between the constructions of the semi-coarse homology and the Vietoris-Rips homology of a graph, is natural to ask whether they are isomorphic. We will answer this in the affirmative in this section, in addition to showing that the semi-coarse homology only depends on the roof
of the semi-coarse structure. We begin with the following lemma.
\begin{lema}\label{lem:Finite subset of roof in structure}
	Let $(X,\mathcal{V})$ be a (possibly non-roofed) semi-coarse space with roof $\mathfrak{A}$. Then every finite subset $A$ of $\mathfrak{A}$ is an element of $\mathcal{V}$.
\end{lema}
\begin{proof}
	The lemma is immediate if $(X,\mathcal{V})$ is roofed. Recall that, by definition, $\mathfrak{A} = \cup_{V \in \mathcal{V}} V$. If $A=\{a_0,\dots,a_k\}$ is a finite subset of $\mathfrak{A}$, then each element $a_i$ of $A$ is contained in some set $V_i \in \mathcal{V}$. By the axioms of a semi-coarse structure, this implies that $\{a_i\}\in\mathcal{V}$ for each $a_i\in A$, and therefore that $A = \cup_{i=0}^k \{a_i\} \in \mathcal{V}$. 
\end{proof}

We now use this to show that the semi-coarse homology only depends on the roof of the semi-coarse structure.

\begin{teorema}
\label{theo:HomologyRoof}
Let $(X,\mathcal{V})$ be a semi-coarse space. Then
\begin{align*}
H(X,\mathcal{V})\cong H(X,\mathfrak{R}(X,\mathcal{V})).
\end{align*}
\end{teorema} 

\begin{proof}
Let $(X,\mathcal{V})$ be a roofed semi-coarse space with roof $\mathfrak{A}$. Then $\{H(X,\mathfrak{A})\}$ is cofinal in the directed system $\{H(X,E),\pi_E^{E'},\mathcal{V}\}$, and the result follows.

Now assume that $\mathcal{V}$ is non-roofed with roof $\mathfrak{A}$, and let $\mathcal{W}$ be the the roofed semi-coarse structure with roof $\mathfrak{A}$. Since $\mathcal{V}\subset\mathcal{W}$, then there exists a homomorphism
\begin{align*}
\pi_{\mathcal{V}}^{\mathcal{W}}:\lim\limits_{\rightarrow}\{H(X,E),\pi_E^{E'},\mathcal{V}\}\rightarrow\lim\limits_{\rightarrow}\{H(X,E),\pi_E^{E'},\mathcal{W}\}
\end{align*}
such that $\pi_{\mathcal{V}}^{\mathcal{W}}\langle\sigma_E \rangle_\mathcal{V}=\langle\sigma_E\rangle_\mathcal{W}$.

\textbf{Surjectivity:} Let $n\in\mathbb{N}$ and $x\in\lim\limits_{\rightarrow}\{H_n(X,E),\pi_E^{E'},\mathcal{W}\}$, then there exists a set
 $E\in\mathcal{W}$ and $y\in H_n(X,E)$ such that $x=\langle y_E \rangle_{\mathcal{W}}$. If $E\in\mathcal{V}$, then $\pi_\mathcal{V}^\mathcal{W}[y_E]_\mathcal{V}=x$.

If $E\notin\mathcal{V}$, we need to do something else. As $y\in H(X,E)$, then there exists a cycle $z\in Ker(\partial_n^E)$ such that $y=\langle z\rangle\in H_n(X,E))$ and
\begin{align*}
z=\sum_{i=1}^q \alpha_i\sigma_i,\ \sigma_i=[z_0^i,\ldots,z_n^i], \alpha_i\in \Z,
\end{align*}
where, the $\sigma_i$ are the elements of the chain complex corresponding to the ordered simplices $[z^i_0,\dots,z^i_n]$. By 
\autoref{lem:Finite subset of roof in structure}, each subset of vertices $A_i\coloneqq \{z^i_0,\dots,z^i_n\} \in \mathcal{V}$ forming 
each ordered simplex $\sigma_i$ is in $\mathcal{V}$, and therefore $A:=\cup_i A_i \in
\mathcal{V}$ as well.

It now follows that $z\in Ker(\partial_n^{A})$, which the properties of the directed system imply that 
$\pi_A^{A\cup E}(\langle z\rangle_A)=\langle z\rangle_{A\cup E}=\pi_E^{A\cup E}(\langle z\rangle_E)$, and it follows that $\langle z \rangle_A$ and
$\langle z \rangle_{E}$ represent the same element of $H_*(X,\mathcal{W})$, i.e. $[\langle z\rangle_{A\cup E}]_\mathcal{W}=[\langle z\rangle_{A}]_\mathcal{W}=x$, i.e. $\langle z \rangle$. it now follows that $\pi_\mathcal{V}^\mathcal{W}[\langle z\rangle_A]_\mathcal{V}=x$. Thus, $\pi_\mathcal{V}^\mathcal{W}$ is surjective.

\textbf{Injectivity:} In this part of the proof we will abuse notation and denote by $\pi_\mathcal{U}^\mathcal{W}$ the homomorphism
\begin{align*}
\pi_\mathcal{V}^\mathfrak{A}: \lim\limits_{\rightarrow}\{H(X,E),\pi_E^{E'},\mathcal{V}\}\rightarrow H(X,\mathfrak{A})
\end{align*}
implicitly composing the original map with the isomorphism
\begin{align*}
\lim\limits_{\rightarrow}\{H(X,E),\pi_E^{E'},\mathcal{W}\}\cong H(X,\mathfrak{A})
\end{align*}
Let $n\in\mathbb{N}$ and $x\in\lim\limits_{\rightarrow}\{H_n(X,E),\pi_E^{E'},\mathcal{V}\}$ such that $\pi_\mathcal{V}^\mathfrak{A}(x)=0$. Then
there exists a set $E\in\mathcal{V}$ and $z\in Ker(\partial_n^E)$ such that $x=[\langle z\rangle_E]_\mathcal{V}$ and $\pi_E^\mathfrak{A}\langle z\rangle=0$. Therefore, $z\in Im(\partial_{n+1}^\mathfrak{A})$, that is, there exists an element $z'\in C_{n+1}^\mathfrak{A}$ such that $\partial_{n+1}^\mathfrak{A}(z')=z$.

Since $\pi_\mathcal{V}^\mathfrak{A}$ is surjective, there is a set $U\in\mathcal{V}$ such that $z'\in C_{n+1}^U$ and $\partial z'=z\in C_n^U$. 
This implies that $\langle z\rangle_U=0_U$ and $\pi_E^U\langle z\rangle_E=0_U$. Thus, $[\langle z\rangle_E]_\mathcal{V}=0\in\lim\limits_{\rightarrow}\{H(X,E),\pi_E^{E'},\mathcal{V}\}$, and $\pi_{\mathcal{V}}^{\mathcal{W}}$ is injective.

It now follows that $\lim\limits_{\rightarrow}\{H(X,E),\pi_E^{E'},\mathcal{V}\}\cong H(X,\mathfrak{A})$.
\end{proof}

We now compare the semi-coarse homology of the vertices of a graph with the semi-coarse structure induced by the graph and the homology of the clique complex of a graph. We begin by recalling the definition of the clique complex.

\begin{definicion}
\label{GraficasCompletas}
A graph is called complete if each pair of vertices is adjacent. A \emph{$k$-clique} in $G$ is a complete subgraph of $G$ with $k$ vertices, and it is a \emph{maximal $k$-clique} if it is not proper subgraph of another clique. 
\end{definicion}

\begin{definicion}
Given a graph $G = (V,E)$ the \emph{clique complex $\Sigma_{G}$} of $G$ is the simplicial complex such that $\Sigma^{(0)}_{G}=V$ and a finite set $\sigma:=\{v_0,\dots,v_k\}\subset V$ is
a $k$-simplex in $\Sigma_{G}$ iff the induced subgraph of $G$ on the vertices in $\sigma$ is a $(k+1)$-clique). The \emph{Vietoris-Rips homology of 
a graph $G$}, denoted $H_*^{VR}(G)$, is the simplicial homology of the simplicial complex $\Sigma_G$.
\end{definicion}

\begin{teorema}
\label{theo:HomologyGraph}
Let $G=(V,E)$ be a graph and let $(V,\mathcal{V}_G)$ be the semi-coarse space induced by $G$. Then $H_*(V,\mathcal{V}_G)\cong H^{VR}_*(G)$.
\end{teorema}

\begin{proof}
Let $G=(V,E)$ be a graph and $(V,\mathcal{V}_G)$ be semi-coarse space associated to $G$. Note that $\cup_{E\in\mathcal{V}}E = \mathfrak{R}(\mathcal{V})$, and
therefore
\begin{align*}
\lim_\rightarrow \{H(V,E),\pi_E^{E'},\mathcal{V}\}= H\left(V,\bigcup_{E\in\mathcal{V}}E\right),
\end{align*}
by thereby we only have to examine the chain complex for the roof of $\mathcal{V}$.

Let $C(G)$ denote the ordered chain complex generated by ordered simplices of the clique complex of $G$, and define $\iota:C(X,\cup_{E\in\mathcal{V}}E)\rightarrow C(G)$ such that $\iota_q ([x_0,\cdots, x_q])=[x_0,\cdots x_q]$ if all the elements are different and $\iota_q ([x_0,\cdots, x_q])=0$ it at least a pair of elements are equal. Note that $\iota$ is well-defined when all of elements are different because $x_i\mathcal{R}x_j$ for each $i,j\in\{0,1,\cdots,q\}$, thereby $\{x_i,x_j\}\in V$ when $i\neq j$, concluding $\{x_0,\cdots, x_q\}$ is a clique with $q+1$ elements in $G$.

On the other hand, let's define $\kappa: C(G)\rightarrow C(X,\cup_{E\in\mathcal{V}}E)$ such that \[\kappa_q([x_0,\cdots,x_q])=[x_0,\cdots,x_q].\] By the argument in the previous paragraph, it is clear that $\kappa$ is well-defined. Moreover, $\iota\circ\kappa = 1_{C(X,\cup_{E\in\mathcal{V}}E)}$ and $\kappa\circ\iota = 1_{C(G)}$, so $C(X,\cup_{E\in\mathcal{V}}E) \cong C(G)$, thus $H(V(G),\mathcal{V})\cong H(G)$.
\end{proof}

\section{Discussion and Future Work}

In this paper, we have begun the development of algebraic topology for semi-coarse spaces, a generalization of the category of coarse spaces which allows, in particular one to coarsen a metric space up to a predetermined scale. The axioms for a semi-coarse structure on a space are obtained by deleting the product axiom from the axioms for coarse spaces, and semi-coarse spaces receive their name in analogy with the loss of the triangle inequality when one passes from metric spaces to semi-metric spaces. Our main results are the development of a homotopy theory on semi-coarse spaces using the integers as an interval, and the construction of Vietoris-Rips homology for semi-coarse spaces. We further define roofed semi-coarse spaces, which clarify the precise structure which determines the Vietoris-Rips homology of a semi-coarse space.

Since we are trying to adapt the algebraic topology of topological spaces to this category, it is tempting to try to construct homotopy invariants using the unit interval as a cylinder object. However, our first observation is that the closed interval $[0,1]$ is not the best candidate for the interval in a semi-coarse homotopy theory:

\begin{customprop}{\ref{prop:ClosedIntervalinSC}}
	Let $r,r'\in (0,1)$. Any bijective map $\phi:([0,1],\mathcal{V}_{r'}) \to ([0,1],\mathcal{V}_{r'})\sqcup_{1\sim 1}([1,2],\mathcal{V}_r)$ where the latter space is the pushout of the two intervals, glued at one endpoint, is not bornologous at the point $1$ where the two intervals are glued.
\end{customprop}

After that, through a long construction, we obtain our homotopy classes and a binary operation for them, and show that they behave similarly to  homotopy classes from topological spaces with their usual product.

\begin{customthm}{\ref{thm:HomotopyGroups}}
Let $(X,\mathcal{V})$ be a semi-coarse space, $A\subset X$, and $n\in\mathbb{N}$.
\begin{itemize}
\item If $n\geq 1$, then $(\pi_n^{sc}(X,*),\star)$ is a group.
\item If $n\geq 2$, then $(\pi_n^{sc}(X,*),\star)$ is an abelian group.
\item If $n\geq 2$, then $(\pi_n^{sc}(X,A,*),\star)$ is a group.
\item If $n\geq 3$, then $(\pi_n^{sc}(X,A,*),\star)$ is an abelian group.
\end{itemize}
\end{customthm}

At the same time, in \autoref{theo:HomotopyClassesCoarse}, we observe that this homotopy becomes trivial on coarse spaces.

\begin{customthm}{\ref{theo:HomotopyClassesCoarse}}
Let $(X,\mathcal{V})$ be a coarse space with $A\subset X$. Then for any $n\in\mathbb{N}$, we have $\pi_n^{sc}(X,A,*)\cong\pi_n^{sc}(X,*)\cong \{1\}$.
\end{customthm}

We conclude the homotopy section with three important results: the existence of semi-coarse spaces whose homotopy is not trivial (the cyclic graphs, \ref{NoTrivHom}), the compression criterion (\ref{CritComp}) and the long exact sequence in homotopy (\ref{SecLarExact}).

\begin{customthm}{\ref{NoTrivHom}}
Let $(C_n,\mathcal{C}_n)$ be the semi-coarse space induced by the $n$-cycle graph. Then, $\pi^{sc}_1(C_n,\mathcal{C}_n)\cong \mathbb{Z}$.
\end{customthm}

\begin{customlem}{\ref{CritComp}}
Let $(X,\mathcal{V})$ be a semi-coarse space and $*\in A\subset X$. Then a function $f:(I^n,\partial I^n,J^{n-1})\rightarrow (X,A,*)$ represents the zero in $\pi_n^{sc}((X,\mathcal{V}),A,*)$ if, and only if, is relative homotopic in $\partial I^n$ to a function with its image contained in $A$.
\end{customlem}

\begin{customthm}{\ref{SecLarExact}}
Let $(X,\mathcal{V})$ be a semi-coarse space and $*\in A\subset X$. As well $i:(A,*)\rightarrow (X,*)$, $j:(X,*)\rightarrow (X,A)$ inclusions. Then there is a homomorphism $\partial_n:\pi^{sc}_n((X,\mathcal{V}),A)\rightarrow \pi^{sc}_n(A,*)$ such that the long sequence

\[
\xymatrix@1{
\cdots \ar[r]^{j_*} & \pi^{sc}_{n+1}(X,A) \ar[r]^{\partial_{n+1}} & \pi^{sc}_{n}(A,*) \ar[r]^{i_*} & \pi^{sc}_n(X,*) \ar[r]^{j_*} & \cdots\\
\cdots \ar[r]^{\partial_2} & \pi^{sc}_{1}(A,*) \ar[r]^{i_*} & \pi^{sc}_1(X,*) \ar[r]^{j_*} & \pi^{sc}_1(X,A)
}
\]

is exact.
\end{customthm}

After our discussion of semi-coarse homotopy groups, we develop a homology with the same flavor as the Vietoris-Rips homology of metric spaces at a given scale. We observe, in particular, that this homology is completely characterized by the roofed semi-coarse spaces (\autoref{theo:HomologyRoof}), which are essentially graphs, and we conclude with a theorem showing that, for any graph, their Vietoris-Rips homology is isomorphic to the semi-coarse homology of the semi-coarse spaces induced by the graph (\autoref{theo:HomologyGraph}). This latter theorem serves to show that the semi-coarse homology developed here is a strict generalization of the homology of the clique complex of a graph, and motivates the possibility of its use in topological data analysis. We further show that the semi-coarse homology is invariant with respect to semi-coarse homotopy equivalence which we defined earlier (\autoref{EquivHom}).

\begin{customthm}{\ref{theo:HomologyRoof}}
Let $(X,\mathcal{V})$ be a semi-coarse space. Then
\begin{align*}
H(X,\mathcal{V})\cong H(X,\mathfrak{R}(X,\mathcal{V})).
\end{align*}
\end{customthm}

\begin{customthm}{\ref{theo:HomologyGraph}}
Let $G=(V,E)$ be a graph and let $(V,\mathcal{V}_G)$ be the semi-coarse space induced by $G$. Then $H_*(V,\mathcal{V}_G)\cong H^{VR}_*(G)$.
\end{customthm}

\begin{customthm}{\ref{EquivHom}}
Let $(X,\mathcal{V})$ and $(Y,\mathcal{W})$ be semi-coarse spaces. If $f,g:(X,\mathcal{V})\rightarrow(Y,\mathcal{W})$ are homotopic functions, then the induced homomorphisms $f_*,g_*:H_*(X,\mathcal{V})\rightarrow H_*(Y,\mathcal{W})$ are equal.
\end{customthm}

There are many directions in which these results may be· developed in future work, particularly as we expect much of algebraic topology to be transferable to this setting. Some specific possibilities are Hurewicz-type theorems, relations between the homologies in semi-uniform spaces and semi-coarse spaces, excision for semi-coarse homology, the development of other semi-coarse homology theories, and the development of semi-coarse invariants which are non-trivial on coarse spaces. In the case of the latter, the second author has constructed a version of the fundamental groupoid for semi-coarse spaces in \cite{Trevino_2024_FundGroupoid_arXiv} and further showed that it may be used to distinguish coarse structures in addition to semi-coarse structure.

\appendix

\section{Additional Results on Point-Set Topology of Semi-Coarse Spaces and Their Relations with Semi-Uniform and Coarse Spaces}

\subsection{Semi-uniform Spaces, Graphs and Roofed Semi-coarse Spaces}

We introduce the notion of semi-uniform spaces 

Finally, we show how to generate a semi-coarse structure from a semi-uniform
space, generalizing Examples \ref{ex:Qc from graphs}, \ref{ex:Qc from metric},
and \ref{ex:Qc from semipseudometric} above.
 
Semi-uniform spaces (described in
detail in \cite{Cech_1966}, Chapter 23) are a generalization of uniform spaces
which are no longer necessarily topological. After the following preliminary
definition, we recall the definition of semi-uniform spaces and give several
examples. We then show how to construct a semi-coarse space from a
semi-uniform space.

\begin{definicion}
    \label{Filtro}A \emph{filter} $\mathcal{U}$ on a set $X$ is a non-empty collection of
    subsets of $X$ such that
\begin{enumerate}[label=(f\arabic*)]
    \item \label{item:Emptyset in filter} $\emptyset \notin \mathcal{U}$,
\item If $A,A'\in\mathcal{U}$, then $A\cap A'\in\mathcal{U}$,
\item If $A\in\mathcal{U}$ and $A\subset A'$, then $A'\in\mathcal{U}$.
\end{enumerate}
A subcollection $\mathcal{U}_0$ of $\mathcal{U}$ is a \emph{filter base} of $\mathcal{U}$ iff
each element of $\mathcal{U}$ contains some element of $\mathcal{U}_0$. 
\end{definicion}

\begin{observacion} A filter is sometimes defined in the literature without 
    \ref{item:Emptyset in filter} above,
    in which case a filter which also satisfies \ref{item:Emptyset in filter} is called a 
    \emph{proper filter}. We will not make this distinction in the present article, and we
    assume that a filter always satisfies \ref{item:Emptyset in filter}.
\end{observacion}

\begin{definicion}[Semi-Uniform Space; \cite{Cech_1966}, 23 A.3.]
\label{EspSU}
Let $X$ be a set and $\mathcal{U}$ be a filter on $X\times X$. We call $\mathcal{U}$ a 
\emph{semi-uniform structure on $X$} and the pair $(X,\mathcal{U})$ 
a \emph{semi-uniform space} iff
\begin{enumerate}[label=(su\arabic*)]
\item Each element $U \in \mathcal{U}$ contains the diagonal $\Delta$ of $X$,
\item \label{item:Symmetry in semi-uniform structure} If $A\in \mathcal{U}$, then $A^{-1}$ contains an element of $\mathcal{U}$,
\end{enumerate}
Since $\mathcal{U}$ is a filter, as noted in \cite{Cech_1966}, Axiom 
\ref{item:Symmetry in semi-uniform structure}
may be replaced by
\begin{enumerate}[label=(su\arabic*')]
    \setcounter{enumi}{1}
\item If $A\in\mathcal{U}$, then $A^{-1}\in\mathcal{U}$.
\end{enumerate}
A semi-uniform space $(X,\mathcal{U})$ is called a \emph{uniform space} iff, in addition to 
the above, $\mathcal{U}$ satisfies
\begin{enumerate}[resume,label=(su\arabic*)]
    \item For every $U \in \mathcal{U}$, there exists a $V \in \mathcal{U}$ such that
        $V \circ V \subset U$.
\end{enumerate}
When context allows it, we will just represent the semi-uniform space by its set $X$.

Let $(X,\mathcal{U})$ and $(Y,\mathcal{T})$ be semi-uniform spaces.
We will say that $f:(X,\mathcal{U})\rightarrow (Y,\mathcal{T})$ is \emph{uniformly continuous} iff
for each $B\in\mathcal{T}$ there exists a $A\in \mathcal{U}$ such that $(f\times f)(A)\subset (B)$.
\end{definicion}

\begin{definition} We denote by $\catname{SUnif}$ the category of semi-uniform spaces and 
    uniformly continuous maps.
\end{definition}

Following \cite{Cech_1966}, Example 23.A.7, we have the following construction of a semi-uniform space from a semi-pseudometric space.

\begin{ejemplo}
    \begin{enumerate}[wide,label=(\arabic*)]
    \label{ex:Semi-uniform from semi-pseudometric}
\item \label{item:SU from sp} Let $(X,d)$ be a semi-pseudometric space, and let $\mathcal{U}_d$ be the semi-uniform structure
on $X\times X$ generated by the collection of sets $\{B_q\}_{q>0}$, where $B_q=\{(x,y): d(x,y)<q\}$.
\item Given a metric space $(X,d)$ and a positive real number $r>0$, the construction in
    \autoref{item:SU from sp} above applied to
    the semi-pseudometrics $d_r, d^<_r$, and $d^\leq_r$ 
    \autoref{def:Metric to semipseudometric}
    gives examples of semi-uniform spaces which are not uniform.
    \end{enumerate}
\end{ejemplo}

\begin{definition}
    Given a semi-pseudometric space $(X,d)$, We call $\mathcal{U}_d$ the
    \emph{semi-uniform structure induced by the semi-pseudometric $d$}.
    If, in addition, $r\geq 0$ is a non-negative real number, then we write 
    $\mathcal{U}_r$ for the semi-uniform structure induced by the 
    semi-pseudometric $d_r$ from \autoref{def:Metric to semipseudometric}. 
\end{definition}

We now show how to build a semi-uniform structure from a semi-coarse structure and vice-versa. 
We first state the following simple lemma.

\begin{lema}
\label{invsub}
Let $X$ be a set and suppose that $A,B\in\mathcal{P}(X\times X)$. Then
\begin{enumerate}[label=(\roman*)]
\item $A=(A^{-1})^{-1}$, and 
\item If $A\subset B$, then $A^{-1}\subset B^{-1}$.
\end{enumerate}
\end{lema}

\begin{proof}
\begin{enumerate}[wide,labelwidth=!,labelindent=0pt,label=(\roman*)]
\item By definition, we have that
\begin{align*}
(A^{-1})^{-1} = & \{(x,y):(y,x)\in A^{-1} \}\\ 
              = & \{(x,y):(x,y)\in A\}\\
              = & A.
\end{align*}
\item Let $(x,y)\in A^{-1}$, then $(y,x)\in A$, so that $(y,x)\in B$. Therefore,
    $(x,y)\in B^{-1}$.
\end{enumerate}
\end{proof}

\begin{proposicion}[Semi-Coarse 
Space from a Semi-Uniform Space]
\label{prop:SU to QC}
Let $(X,\mathcal{U})$ be a semi-uniform space, and define the collection $\mathcal{U}^\downarrow 
\subset \mathcal{P}(X \times X)$ by
$$\mathcal{U}^\downarrow:=\left\{ B\subset X\times X: B\subset \bigcap_{A\in \mathcal{U}} A \right\}.$$
Then $(X,\mathcal{U}^\downarrow)$ is a semi-coarse space.
\end{proposicion}

\begin{proof}
Let $(X,\mathcal{U})$ and $\mathcal{U}^\downarrow$ be as in the statement of the Proposition.
We check that the axioms for semi-coarse structures are satisfied by $\mathcal{U}^\downarrow$.
\begin{enumerate}[wide,labelwidth=!,label=(sc\arabic*)]
\item By definition of a semi-uniform space, the diagonal $\Delta_X \subset U$ for
    each $U\in\mathcal{U}$. Therefore, $\Delta_X\in\mathcal{U}^\downarrow$.
\item Let $B\in \mathcal{U}^\downarrow$ and $B'\subset B$. Then $B'\subset B\subset A$ for each $A\in\mathcal{U}$. So $B'\in \mathcal{U}^\downarrow$.
\item Let $B,B'\in\mathcal{U}^\downarrow$. By definition, we have that $B,B'\subset A$ for
    every $A\in\mathcal{U}$, and therefore $B\cup B'\subset A$ for every $A\in\mathcal{U}$ as well. 
    Therefore, $B\cup B' \in\mathcal{U}^\downarrow$.
\item Let $B\in\mathcal{U}^\downarrow$, then $B\subset A$ for each $A\in\mathcal{U}$.
    By \autoref{invsub}, we have that $B^{-1}\subset A^{-1}$ for each $A\in\mathcal{U}$,
    and since $A\in\mathcal{U}$ implies $A^{-1}\in\mathcal{U}$, \autoref{invsub} further implies that 
    $B^{-1}\subset A$ for each $A\in\mathcal{U}$. Therefore, $B^{-1}\in\mathcal{U}^\downarrow$.
\end{enumerate}
It now follows that $(X,\mathcal{U}^\downarrow)$ is a semi-coarse space.
\end{proof}

\begin{definicion}\label{SemUnifPseCoar}
    The semi-coarse structure $\mathcal{U}^\downarrow$ in \autoref{prop:SU to QC} is called
\emph{the semi-coarse structure induced by the semi-uniform structure $\mathcal{U}$}.
\end{definicion}

We now construct a semi-uniform structure $\mathcal{V}^\uparrow$ from a semi-coarse 
space $(X,\mathcal{V})$.

\begin{proposicion}
    \label{prop:QC to SU}
    Let $(X,\mathcal{V})$ be a semi-coarse space, and define the collection 
    $\mathcal{V}^\uparrow \subset \mathcal{P}(X \times X)$ by 
    \[ \mathcal{V}^\uparrow=\left\{ U\subset X\times X : \left(\bigcup_{V \in \mathcal{V}} V\right) \subset U
        \right\}.\]
Then $(X,\mathcal{V}^\uparrow)$ is a semi-uniform space.
\end{proposicion}

\begin{proof}
Let $(X,\mathcal{V})$ and $\mathcal{V}^\uparrow$ be as in the statement of the
proposition. We first show that $\mathcal{V}^\uparrow$ is a filter.
\begin{enumerate}[wide,label=(f\arabic*)]
    \item Since $\Delta_X \in \mathcal{V}$, 
        $\Delta_X \subset \left(\cup_{V \in \mathcal{V}} V\right) \subset A \cap A'$, and, 
        in particular, $\emptyset \neq A \cap A'$.

    \item Suppose that $A,A'\in \mathcal{V}^\uparrow$
        Then $(\cup_{V \in \mathcal{V}} V) \subset A\cap A'$ by construction, so 
        $A \cap A' \in \mathcal{V}^\uparrow$.
    \item Suppose that $A\in\mathcal{V}^\uparrow$ and $A'\subset X\times X$ such that $A\subset A'$. 
        Then $(\cup_{V \in \mathcal{V}} V) \subset A\subset A'$, and therefore
        $A'\in\mathcal{V}^\uparrow$ as well.
\end{enumerate}
We now prove that $\mathcal{V}^\uparrow$ satisfies the axioms of a semi-uniform structure.
\begin{enumerate}[wide,label=(su\arabic*)]
    \item Since $\Delta_X\in\mathcal{V}$, then $\Delta_X\subset (\cup_{V \in \mathcal{V}} V) \subset 
        U$ for each $U\in\mathcal{V}^\uparrow$.
    \item[(su2')] Since $V \in \mathcal{V} \iff V^{-1} \in \mathcal{V}$, it follows that 
        $\cup_{V \in \mathcal{V}} V = \cup_{V \in \mathcal{V}} V^{-1}$. Now suppose that 
        $A\in\mathcal{V}^\uparrow$. Then $\cup_{V \in \mathcal{V}} V =
        \cup_{V \in \mathcal{V}} V^{-1} \subset A^{-1}$, so $A^{-1}\in\mathcal{V}^\uparrow$.
\end{enumerate}
It follows that $(X,\mathcal{V}^\uparrow)$ is a semi-uniform space.
\end{proof}

\begin{definition}
    We call the semi-uniform structure $\mathcal{V}^\uparrow$ in \autoref{prop:QC to SU} the
    \emph{Semi-uniform structure induced by the semi-coarse space $(X,\mathcal{V})$}.
\end{definition}

Propositions \ref{prop:SU to QC} and \ref{prop:QC to SU} motivate the following definitions.

\begin{definicion}{(Roof)}
\label{DefRoof}
Let $(X,\mathcal{V})$ be a semi-coarse space.
We call $\mathfrak{R}(\mathcal{V}) \coloneqq \bigcup_{V\in\mathcal{V}} V$ the \emph{roof} of 
    $\mathcal{V}$, and we say that the semi-coarse space $(X,\mathcal{V})$ is \emph{roofed} iff 
    $\mathfrak{R}(\mathcal{V})\in\mathcal{V}$. Otherwise, we say that $(X,\mathcal{V})$ is 
    \emph{non-roofed}.
\end{definicion}

\begin{ejemplo}\label{ex:Power set roof}
   \begin{enumerate}[wide]
   	\item The semi-coarse spaces in \autoref{lem:Qc from power set} are roofed with roof $W$.
	\item The semi-coarse spaces induced by graphs as defined in
	\autoref{def:SC from graphs} are roofed with roof $E$, the set of edges of the graph.
\end{enumerate}
\end{ejemplo}

\begin{proposicion}
	\label{RoofedGraphs}
	Let $\catname{SCoarse}$, $\catname{RSCoarse}$ and $\catname{Graphs}$ the categories of semi-coarse spaces, roofed semi-coarse spaces and graphs, respectively. Let $\Phi: \catname{SCoarse}\rightarrow \catname{Graphs}$ be a map of categories defined by 
	\begin{align*}
		\Phi(X,\mathcal{V})&=(V,E),\\
		\Phi(f) &= f,
	\end{align*} 
where $V=X$ and $E=\{\{u,v\} \mid (u\neq v) \text{ and }(u,v)\in \mathcal{V}\}$.

Given a graph $G=(V,E)$, define 
\begin{align*}
	\mathfrak{R}_E\coloneqq \Delta_{V} \cup \left( \bigcup_{\{u,v\}\in E} \{(u,v),(v,u)\} \right),
\end{align*}
and let $\Psi: \catname{Graphs}\rightarrow \catname{SCoarse}$ be the map of categories defined by 
\begin{align*}
	\Psi(V,E) &= (V,\mathcal{P}(\mathfrak{R}_E))\\
	\Psi(f) &= f
\end{align*}
where $\mathcal{P}(\mathfrak{R}_E)$ is the power set of $\mathfrak{R}_E$. 
	
	Then $\Phi$ and $\Psi$ are functors. Moreover, $\Phi|_{\catname{RSCoarse}}\circ \Psi=1_{\catname{Graphs}}$ and $\Psi\circ \Phi|_{\catname{RSCoarse}}=1_{\catname{RSCoarse}}$, so $\Psi(\catname{Graphs})=\catname{RSCoarse}$ and $\Phi(\catname{SCoarse})=\catname{Graphs}$.
\end{proposicion}

\begin{proof}
	
	We first need to show that $\Phi$ maps bornologous functions to graph morphisms, and that $\Psi$ maps graph morphisms to bornologous functions.
	
	Let $(X,\mathcal{V})$ and $(Y,\mathcal{W})$ be semi-coarse spaces, and let $f:X\rightarrow Y$ be a bornologous function. We recall that $f:G\to G'$ is a graph morphism iff, for all $(v,v')\in E_G$, either $f(v)=f(v')$ or $(f(v),f(v'))\in E_{G'}.$ Let $x,y\in X$ be such that $f(x)\neq f(y)$. Then $\{(f(x),f(y))\}\in\mathcal{W}$, so $\{f(x),f(y)\}\in E_Y$ by construction. Thus $\Phi(f)$ is a graph morphism. Since $\Phi$ preserves the identity and respects composition, it follows that $\Phi$ is a functor from $\catname{SCoarse}$ to $\catname{Graphs}$.
	
    Now let $G$ and $G'$ be graphs and $f:G\rightarrow G'$ a graph morphism. Suppose that $(v,w) \in G_E$. Then $(f(v),f(w))\in \mathfrak{R}_{E_{G'}}$ by definition. Furthermore, by the definition of $\mathfrak{R}_{E_G}$ and $\mathfrak{R}_{E_{G'}}$, it follows that $f(\mathfrak{R}_{E_G}) \subset \mathfrak{R}_{E_{G'}}$, and therefore $f$ is bornologous. Since $\Psi$ also preserves the identity and respects composition, it follows that $\Psi$ is a functor from $\catname{Graphs}$ to $\catname{SCoarse}$.
	
	Note that, by construction, $\Psi(V,E)$ is the roofed semi-coarse structure on $V$ with roof $\mathfrak{R}_E$. Finally, the last sentence of the theorem follows by construction.
\end{proof}

\begin{definicion}
\label{Cemented}
Let $(X,\mathcal{U})$ be a semi-uniform space. We call $\mathfrak{F}(\mathcal{U}) \coloneqq 
\bigcap_{U\in\mathcal{U}} U$ the \emph{foundation} of the semi-uniform structure $\mathcal{U}$. 
We say that $(X,\mathcal{U})$ is \emph{cemented} (from the Spanish \emph{cimentado}) iff 
$\mathfrak{F}(\mathcal{U})\in\mathcal{U}$. Otherwise, we say that $(X,\mathcal{U})$ is 
\emph{non-cemented}.
\end{definicion}

\begin{observacion} Note that the roof is the maximal in the semi-coarse
    structure with  order relation $A<B$ iff $A\subset B$, while the foundation
    is the minimal in the semi-uniform structure with the same order.    
\end{observacion}

\begin{proposicion} \label{prop:Unicity of foundations}

\begin{enumerate}[wide]
    \item Let $U\subset X\times X$ such that $\Delta _X\subset U$ and $U =U^{-1}$, and denote 
        by $[U]$ the filter generated by $U$. Then $(X,[U])$ is a cemented semi-uniform space 
        with foundation $U$.
    \item Let $\mathcal{U}$ and $\mathcal{T}$ be semi-uniform structures on $X$. 
        If $\mathfrak{F}(\mathcal{U})=\mathfrak{F}(\mathcal{T})$ and $(X,\mathcal{U})$ is cemented, 
        then $\mathcal{T} \subset \mathcal{U}$. It follows that, for each foundation, there exists 
        a unique cemented semi-uniform space with that foundation.
\end{enumerate}
\end{proposicion}

\begin{proof}
    \begin{enumerate}[wide]
        \item Let $U\subset X\times X$ such that $\Delta _X\subset U$ and $U =
            U^{-1}$, and denote by $[U]$ the filter generated by
            $U$. Since $\Delta_X\subset U$, then
            $\Delta_X\subset A$ for all $A\in [U]$. Furthermore, if $A\in
            [U]$, then $U\subset A$ and
            $U = U^{-1} \subset A^{-1}$, from which it follows that $A^{-1}\in
            [U]$. Thus, $(X, [U])$ is a semi-uniform space.

            More over, since $U\subset A$ for all $A\in [U]$, then
            
            \begin{align*} U \subset \bigcap\limits_{A\in [U]} A
\subset U \end{align*} 
Thus $(X,[U])$ is a cemented
semi-uniform space with foundation $[U]$.

\item Let $\mathcal{U}$ and $\mathcal{T}$ be semi-uniform structures on $X$
    such that $\mathfrak{F}(\mathcal{U})=\mathfrak{F}(\mathcal{T})$, and suppose that
    $(X,\mathcal{U})$ is cemented. Therefore, for all $B\in\mathcal{T}$, $\mathfrak{F}(\mathcal{U})
    \subset B$, and since $\mathfrak{F}(\mathcal{U})$, therefore $B\in\mathcal{U}$. Thus,
    $\mathcal{T}\subset \mathcal{U}$ and the result follows.\qedhere
\end{enumerate}
\end{proof}

\begin{proposicion}
    \label{prop:Unicity of roofs}
    Let $\mathcal{V}$ and $\mathcal{T}$ be a semi-coarse structures on a set $X$. Suppose that
    $\mathcal{V}$ is roofed, and that $\mathfrak{R}(\mathcal{V}) = \mathfrak{R}(\mathcal{T})$.
    Then $\mathcal{T} \subset \mathcal{V}$, and for each roof, there exists a unique
    roofed semi-coarse structure with that roof.
\end{proposicion}

\begin{proof} For all $B \in \mathcal{T}$, $B \subset \mathfrak{R}(\mathcal{T}) = 
    \mathfrak{R}(\mathcal{V})\in \mathcal{V}$, and therefore $B \in \mathcal{V}$.
    The result follows.
\end{proof}

\begin{proposicion}
    \label{prop:Born on roof is born}
    Let $(X,\mathcal{V})$ and $(Y,\mathcal{W})$ be semi-coarse spaces, and suppose 
    that $(X,\mathcal{V})$ is roofed. A map $f:X \to Y$ is bornologous iff
    $(f \times f)(\mathfrak{R}(\mathcal{V})) \in \mathcal{W}$.
\end{proposicion}

\begin{proof}
    Any set $V \in \mathcal{V}$ is a subset of 
    $\mathfrak{R}(\mathcal{V})$ by definition. Therefore, if $(f \times f)(V) \subset
    (f \times f)(\mathfrak{R}(V)) \in \mathcal{W}$, then $(f \times f)(V) \in\mathcal{W}$,
    so $f$ is bornologous.

    Conversely, suppose that $f$ is bornologous. Then $(f \times f)(\mathfrak{R}(\mathcal{V})) \in
    \mathcal{W}$.
\end{proof}

Similarly, for cemented semi-uniform spaces, we have

\begin{proposicion}
\label{UnifContCemented}
Let $(X,\mathcal{U})$ and $(Y,\mathcal{T})$ be semi-uniform spaces and suppose that 
$(Y,\mathcal{T})$ is cemented. A map $f:X\rightarrow Y$ is uniformly continuous iff 
$(f\times f)(\mathfrak{F}(\mathcal{U})) \subset \mathfrak{F}(\mathcal{T})$.
\end{proposicion}

\begin{proof}
Let $(X,\mathcal{U})$ and $(Y,\mathcal{T})$ be cemented semi-uniform spaces, and let
$f:(X,\mathcal{U}) \rightarrow (Y,\mathcal{T})$ be a uniformly continuous function. 
Then, for each $B\in\mathcal{T}$ there exists a set $A_B\in\mathcal{U}$ such that 
$(f\times f)(A_B)\subset B$. Therefore

\begin{align*}
    (f \times f)\left( \mathfrak{F}(\mathcal{U})\right) = & (f\times f)\left( \bigcap_{A\in\mathcal{U}} A 
        \right) 
            \subset \bigcap_{A\in\mathcal{U}}(f\times f)\left( A \right)\\
            \subset & \bigcap_{B\in\mathcal{T}}(f\times f)\left( A_B \right)
            \subset  \bigcap_{B\in\mathcal{T}} B\\
            = & \,\mathfrak{F}(\mathcal{T}).
\end{align*} 

Conversely, assume that $(f\times f)(\mathfrak{F}(\mathcal{U}))\subset \mathfrak{F}(\mathcal{T})$. 
If $B\in\mathcal{T}$, then $\mathfrak{F}(\mathcal{T})\subset B$, and therefore $(f\times f)
(\mathfrak{F}(\mathcal{U}))\subset B$ as well. However, since $(X,\mathcal{U})$ is cemented,
we have that $\mathfrak{F}(\mathcal{U}) \in \mathcal{U}$, from which it follows that $f$ is 
uniformly continuous.
\end{proof}

\begin{proposicion}
\label{Doublearrow}
Let $X$ be a set, and suppose that $\mathcal{V}$ and $\mathcal{U}$ are semi-coarse and
semi-uniform structures on $X$, respectively. Then 
\begin{enumerate}[wide,label=(\roman*)]
    \item \label{item:SU to roofed qc space}$(X,\mathcal{U}^\downarrow)$ is a 
        roofed semi-coarse space.

    \item \label{item:QC to cemented su space}$(X,\mathcal{V}^\uparrow)$ is a 
        cemented semi-uniform space.

\item If $(X,\mathcal{U})$ is a cemented semi-uniform space, then
    $(\mathcal{U}^\downarrow)^\uparrow=\mathcal{U}$.

\item If $(X,\mathcal{V})$ is a roofed semi-coarse space, then
    $(\mathcal{V}^\uparrow)^\downarrow=\mathcal{V}$.
\end{enumerate}
\end{proposicion}

\begin{proof}
    Let $X$, $\mathcal{V}$, and $\mathcal{U}$ be as in the hypothesis of the proposition.
\begin{enumerate}[wide,label=(\roman*)]

    \item First note that, by definition, $\mathfrak{F}(\mathcal{U}) \in \mathcal{U}^\downarrow$.
        Now suppose that $A\in\mathcal{U}^\downarrow$. By construction, $A\subset B$ for every 
    $B\in\mathcal{U}$, and therefore $A\subset \mathfrak{F}(\mathcal{U})$. Therefore, 
    $\mathfrak{R}(\mathcal{U}^\downarrow)=\mathfrak{F}(\mathcal{U}) \in \mathcal{U}^\downarrow$,
    and $(X,\mathcal{U}^\downarrow)$ is a roofed semi-coarse space.

\item As before, note that $\mathfrak{R}(\mathcal{V}) \in \mathcal{V}^\uparrow$ by definition. 
    Now let $B\in\mathcal{V}^\uparrow$. By construction, $A\subset B$ for every $A\in\mathcal{V}$,
    and therefore $\mathfrak{R}(\mathcal{V} \subset B$. It follows that 
    $\mathfrak{F}(\mathcal{V}^\uparrow)=\mathfrak{R}(\mathcal{V})$, so $(X,\mathcal{V}^\uparrow)$ 
    is a cemented semi-uniform space.

\item By points \ref{item:SU to roofed qc space} and \ref{item:QC to cemented su space} above, 
    $\mathfrak{F}(\mathcal{U})=\mathfrak{R}(\mathcal{U}^\downarrow)=\mathfrak{F}
    (\mathcal{U}^\downarrow)^\uparrow)$. Therefore, $\mathcal{U}=(\mathcal{U}^\downarrow)^\uparrow$
    by \autoref{prop:Unicity of foundations}.

\item By points \ref{item:SU to roofed qc space} and \ref{item:QC to cemented su space}
    above, $\mathfrak{R}(\mathcal{V})=\mathfrak{F}(\mathcal{V}^\uparrow)=
    \mathfrak{R}(\mathcal{V}^\uparrow)^\downarrow$. Therefore, 
    $\mathcal{V}=(\mathcal{V}^\uparrow)^\downarrow$ by \autoref{prop:Unicity of roofs}.
\end{enumerate}
\end{proof}

We now show that the constructions above are functorial. Let $\catname{RSCoarse}$ and
$\catname{CSUnif}$ denote the full subcategories of roofed semi-coarse spaces and cemented
semi-uniform spaces, respectively.

\begin{proposicion}
\label{prop:Functoriality QC SU}
Let $\Phi:\catname{SCoarse} \to \catname{SUnif}$ be the map 
$\Phi(X,\mathcal{V})=(X,\mathcal{V}^\uparrow)$ on objects and $\Phi(f:X\rightarrow Y)=f$
on morphisms, and let $\Psi:\catname{SUnif} \to \catname{SCoarse}$ be the map
$\Psi(X,\mathcal{U})=(X,\mathcal{U}^{\downarrow})$ on objects and $\Psi(f:X\rightarrow Y)=f$ 
on morphisms. 

Then $\Phi$ and $\Psi$ are functors, $\Phi(\catname{SCoarse})=\catname{CSUnif}$, and 
$\Psi(\catname{CSUnif})=\catname{RSCoarse}$. Moreover, we have 
$\Phi|_{\catname{RSCoarse}} \circ \Psi|_{\catname{CSUnif}}=
1_{\catname{CSUnif}}$, and $\Psi|_{\catname{CSUnif}} \circ \Phi|_{\catname{RSCoarse}}=
1_{\catname{RSCoarse}}$.
\end{proposicion}

\begin{proof}
It suffices to prove that $\Phi$ maps bornologous function to uniformly continuous functions
and that $\Psi$ maps uniformly continuous functions to bornologous functions. 
The rest of the proposition is a direct consequence of \autoref{Doublearrow}.

Let $f:(X,\mathcal{V})\rightarrow (Y,\mathcal{W})$ be a bornologous
function. Then $(f\times f)(A)\in\mathcal{W}$ for every $A\in\mathcal{V}$.
We therefore have that 
\begin{align*} (f\times f)\left(\bigcup_{A\in\mathcal{V}} A
\right)\subset \bigcup_{B\in\mathcal{W}} B. 
\end{align*} 
Thus, the
foundation of $\mathcal{W}^\uparrow$ contains the image of the foundation of
$\mathcal{V}^\uparrow$ under $f\times f$, i.e., $(f\times
f)(\mathfrak{F}(\mathcal{V}^\uparrow))\subset U$ for every $U\in
\mathcal{W}^\uparrow$. Therefore $f:(X,\mathcal{U}^\uparrow) \to (Y,\mathcal{W}^\uparrow)$
is uniformly continuous by \autoref{UnifContCemented}.

Now let $f:(X,\mathcal{U})\rightarrow (Y,\mathcal{T})$ be a uniformly continuous function.
By definition, for every $T\in\mathcal{T}$, there exists a $U_T\in\mathcal{U}$ such that 
$(f\times f)(U_T)\subset T$. We therefore have that
\begin{align*}
    (f \times f)\left(\mathfrak{F}(\mathcal{U})\right) = & 
        (f\times f)\left( \bigcap_{U\in\mathcal{U}}U \right) \subset 
        \bigcap_{U\in\mathcal{U}}(f\times f)(U)\\
    \subset & \bigcap_{T\in\mathcal{T}}(f\times f)(U_T) \subset
        \bigcap_{T\in\mathcal{T}}T\\
    = & \mathfrak{F}(\mathcal{T}).
\end{align*}
Therefore, the roof of $\mathcal{T}^\downarrow$ contains the image of the roof of $\mathcal{U}$
under $f\times f$. It follows that $(f\times f)(A)\subset \mathfrak{R}(Y,\mathcal{T}^\downarrow)$, or,
equivalently, $(f\times f)(A)\in\mathcal{T}^\downarrow$ for every $A\in\mathcal{U}^\downarrow$.
Therefore, $f:(X,\mathcal{U}^\downarrow) \to (Y,\mathcal{T}^\downarrow)$ is a bornologous
function.
\end{proof}

In addition to the above, the next result shows that, when restricted to cemented semi-uniform 
spaces and roofed semi-coarse spaces, the functors in \autoref{prop:Functoriality QC SU} 
are adjoints of each other.

\begin{proposicion}
    $\Phi$ and $\Psi$ satisfy $\Phi \dashv \Psi$, where $\Phi$ and $\Psi$ are the functors defined in \autoref{prop:Functoriality QC SU}.
\end{proposicion}

\begin{proof}
Let $(X,\mathcal{V})$ be a semi-coarse space, and let $(Y,\mathcal{U})$ be a semi-uniform space. We first show that $\Phi \dashv \Psi$. Let $f\in Hom_{\catname{SCoarse}}(X,\Psi(Y))$. Then, by definition, $(f\times f)(A)\in \mathcal{U}^\downarrow$ for all $A\in \mathcal{V}$. This equivalent to the statement that $(f\times f)(A)\subset \mathfrak{F}(\mathcal{U})$ for every $A\in\mathcal{V}$. This, in turn, is equivalent to $(f\times f)(\mathfrak{R}(\mathcal{V}))\subset \mathfrak{F}(\mathcal{U})$, since $A\subset \cup_{A'\in\mathcal{V}}A'$ for every $A\in\mathcal{V}$. Finally, we have that $(f\times f)(A)\subset \mathfrak{F}(\mathcal{U})$ for every $A\in\mathcal{V}$ implies $\cup_{A\in\mathcal{V}}(f\times f)(A)=(f\times f)\left( \cup_{A\in\mathcal{V}} A  \right)\subset \mathfrak{F}(\mathcal{U})$. Therefore, $f \in Hom_{\catname{SUnif}}(\Phi(X),Y)$ by \autoref{UnifContCemented}.

Similarly, let $f\in Hom_{\catname{SUnif}}(\Phi(X),Y)$, by definition for every $B\in\mathcal{U}$ there exists $A\in\mathcal{V}^\uparrow$ such that $(f\times f)(A)\subset B$. This equivalent to the fact that, for every $B\in\mathcal{U}$ we have $(f\times f)(\mathfrak{R}(\mathcal{V}))\subset B$, by definition of $\mathfrak{R}(\cdot)$ and the construction of $\mathcal{U}^\uparrow$. Since $\mathfrak{F}(\mathcal{U})=\cap_{B\in\mathcal{U}}B$, this identical $(f\times f)(\mathfrak{R}(\mathcal{V}))\subset \mathfrak{F}(\mathcal{V})$. Therefore, $f \in Hom_{\catname{SCoarse}}(X,\Psi(Y))$ by \autoref{prop:Born on roof is born}.

It now follows that $Hom_{\catname{SCoarse}}(X,\Psi(Y)) \cong Hom_{\catname{SUnif}}(\Phi(X),Y))$, so $\Phi \dashv \Psi$. 

\end{proof}

\subsection{Product Structure in Semi-coarse Spaces}

We discuss infinite products of semi-coarse spaces. The next two propositions give semi-coarse analogues of 
the box product and (Tychonoff) product of topological spaces.

\begin{proposicion}
    \label{prop:Box product}
Let $\Lambda$ be an index set and $\{(X_\lambda,\mathcal{V}_\lambda)\}$ be a collection of semi-coarse spaces.
Let $\prod_{\lambda\in\Lambda} \mathcal{V}_\lambda$ denote the collection of subsets of 
$\left(\prod_{\lambda \in \Lambda} X_\lambda\right)^2$ that are contained in some set of the form
\begin{equation*} \left\lbrace (x,y)\in\left( \prod_{\lambda\in\Lambda} X_\lambda \right)^2 : 
    (\pi_\lambda(x),\pi_\lambda(y))\in A_\lambda\mbox{ for each } \lambda\in\Lambda \right\rbrace 
\end{equation*}
with $A_\lambda\in\mathcal{V}_\lambda$.

Then $\prod_{\lambda \in \Lambda} \mathcal{V}_\lambda$ is a semi-coarse structure on the set
$\prod_{\lambda\in\Lambda} X_\lambda$.
\end{proposicion}

\begin{proof}
    Let $\Lambda$, $\{(X_\lambda,\mathcal{V}_\lambda)\}_{\lambda \in \Lambda}$, and 
    $\prod_{\lambda\in\Lambda} \mathcal{V}_\lambda$ be as in the statement of the proposition. We check
    that $\prod_{\lambda \in \Lambda} \mathcal{V}_\lambda$ satisfies the axioms of a semi-coarse structure.
    We have
    \begin{enumerate}[wide,label=(sc\arabic*)]
\item $\Delta_{(\prod_{\lambda\in\Lambda}X_\lambda)}=\prod_{\lambda\in\Lambda} \Delta_{X_\lambda}\in 
    \prod_{\lambda\in\Lambda} \mathcal{V}_\lambda$.
\item If $V\in \prod_{\lambda\in\Lambda} \mathcal{V}_\lambda$ and $W$ is subset of $V$, then there are 
    sets $\{A_\lambda \in \mathcal{V}_\lambda\}_{\lambda \in \Lambda}$ such that $W\subset V\subset 
    \prod_{\lambda\in\Lambda} A_\lambda$, and therefore $W\in \prod_{\lambda\in\Lambda} \mathcal{V}_\lambda$.
\item If $V,W\in\prod_{\lambda\in\Lambda} \mathcal{V}_\lambda$, then there are $\{A_\lambda, B_\lambda \in 
        \mathcal{V}_\lambda\}_{\lambda \in \Lambda}$ such that $V\subset \prod_{\lambda\in\Lambda} A_\lambda$ 
        and $W\subset\prod_{\lambda\in\Lambda} B_\lambda$. Therefore, $V\cup W \subset \prod_{\lambda\in\Lambda} 
        (A_\lambda\cup B_\lambda)$, and we conclude that $V\cup W\in \prod_{\lambda\in\Lambda}\mathcal{V}_\lambda$.
    \item If $V\in\prod_{\lambda\in\Lambda} \mathcal{V}_\lambda$, then there are 
        $\{A_\lambda \in \mathcal{V}_\lambda\}_{\lambda \in \Lambda}$ such that $V\subset \prod_{\lambda\in\Lambda} 
        A_\lambda$, and therefore $V^{-1}\subset \prod_{\lambda\in\Lambda} (A_\lambda)^{-1}$, which gives
        that $V^{-1}\in\prod_{\lambda\in\Lambda} \mathcal{V}_\lambda$.\qedhere
\end{enumerate}
\end{proof}

\begin{definicion}
    \label{def:Box product}
    The semi-coarse structure $\prod_{\lambda \in \Lambda} \mathcal{V}_\lambda$ from \autoref{prop:Box product}
    will be called the \emph{box product semi-coarse structure on 
        $\left(\prod_{\lambda \in \Lambda} X_\lambda\right)$}. We call the pair 
        $\left(\prod_{\lambda \in \Lambda} X_\lambda,\prod_{\lambda \in \Lambda} \mathcal{V}_\lambda\right)$ the
        \emph{box product of the family of semi-coarse spaces $\{(X_\lambda,\mathcal{V}_\lambda\}$.}
\end{definicion}

\begin{proposicion} 
    \label{prop:Tychonoff product} 
    Let  $\prod_{\lambda\in\Lambda}^\tau\mathcal{V}_\lambda$ be the collection of subsets of 
    $\left( \prod_{\lambda\in\Lambda} X_\lambda \right)^2$ such they are contained in a set of the form
\begin{align*}
    \left\lbrace (x,y)\in\left( \prod_{\lambda\in\Lambda} X_\lambda \right)^2 \mid \right. 
& (\pi_\lambda(x),\pi_\lambda(y))\in A_\lambda\mbox{ if } \lambda\in F,\\ 
& \left. {\vphantom{\left( \prod_{\lambda\in\Lambda} X_\lambda \right)^2}}
(\pi_\lambda(x),\pi_\lambda(y))\in \Delta_{X_\lambda}\mbox{ if } \lambda\in \Lambda \backslash F \right\}, 
\end{align*}
where $F\subset \Lambda$ is a finite set and $A_\lambda\in\mathcal{V}_\lambda$ for all $\lambda \in \Lambda$.

Then $\prod_{\lambda\in\Lambda}^\tau \mathcal{V}_\lambda$ is a semi-coarse structure on
$\prod_{\lambda\in\Lambda} X_\lambda$.

\end{proposicion}
\begin{proof}
Let $\prod_{\lambda\in\Lambda}^\tau \mathcal{V}_\lambda$ be as in the statement of the proposition. We 
check that $\prod_{\lambda \in \Lambda}^\tau \mathcal{V}_\lambda$ the axioms for semi-coarse structure.
\begin{enumerate}[wide,label=(sc\arabic*)]
\item $\Delta_{(\prod_{\lambda\in\Lambda}X_\lambda)}=\prod_{\lambda\in\Lambda}
    \Delta_{X_\lambda}\in \prod^\tau_{\lambda\in\Lambda} \mathcal{V}_\lambda$.
\item If $V\in \prod^\tau_{\lambda\in\Lambda} \mathcal{V}_\lambda$ and $W$ subset of $V$, then there are sets
    $\{V_\lambda \in \mathcal{V}_\lambda\}_{\lambda \in \Lambda}$ such that a finite number of
    $V_\lambda \neq \Delta_{X_\lambda}$ and $W\subset V\subset \prod_{\lambda\in\Lambda} V_\lambda$, 
    from which it follows that $W\in \prod^\tau_{\lambda\in\Lambda} \mathcal{V}_\lambda$.
\item If $V,W\in\prod^\tau_{\lambda\in\Lambda} \mathcal{V}_\lambda$, then there are sets 
    $\{V_\lambda \in \mathcal{V}_\lambda\}_{\lambda \in \Lambda}, 
    \{W_\lambda \in \mathcal{V}_\lambda\}_{\lambda \in \Lambda}$ such that a finite number $n \in \N$ of the 
    $V_\lambda$ satisfy $V_\lambda\neq \Delta_{X_\lambda}$ and a finite number $m$ of the $W_\lambda$ satisfy
    $W_\lambda\neq \Delta_{X_\lambda}$, and $V\subset \prod_{\lambda\in\Lambda} V_\lambda$, 
    $W\subset\prod_{\lambda\in\Lambda} W_\lambda$. Therefore, $V\cup W \subset 
    \prod_{\lambda\in\Lambda} (V_\lambda\cup W_\lambda)$ and there are at most a finite number $m+n$ 
    of unions $V_\lambda\cup W_\lambda$ which are not equal to the corresponding 
    diagonal $\Delta_{X_\lambda}$. It follows that
    $V\cup W\in \prod^\tau_{\lambda\in\Lambda}\mathcal{V}_\lambda$.
\item If $V\in\prod^\tau_{\lambda\in\Lambda} \mathcal{V}_\lambda$, then there are sets 
    $\{V_\lambda \in \mathcal{V}_\lambda\}_{\lambda\in\Lambda}$ such that for a finite number 
    of them, $V_\lambda\neq\Delta_{X_\lambda}$, and we have $V\subset \prod_{\lambda\in\Lambda} V_\lambda$.
    It follows that $V^{-1}\subset \prod_{\lambda\in\Lambda} (V_\lambda)^{-1}$, which implies
    that $V^{-1}\in\prod^\tau_{\lambda\in\Lambda} \mathcal{V}_\lambda$. \qedhere
\end{enumerate}
\end{proof}

\begin{definicion}
\label{def:Tychonoff product}
We call the semi-coarse structure $\prod^\tau_{\lambda\in\Lambda} \mathcal{V}_\lambda$ from
\autoref{prop:Tychonoff product} the \emph{Tychonoff product of the family 
$\{\mathcal{V}_\lambda\}_{\lambda\in \Lambda}$}, and we call the pair $\left(\prod_{\lambda \in \Lambda} X_\lambda,
    \prod^\tau_{\lambda \in \Lambda} \mathcal{V}_\lambda \right)$ the \emph{Tychonoff product of the family
        of semi-coarse spaces $\{(X_\lambda,\mathcal{V}_\lambda)\}_{\lambda \in \Lambda}$}.
\end{definicion}

The following proposition shows us that the box product is the preferred product for products
of infinite families of semi-coarse spaces. Naturally, when $|\Lambda| < \infty$, the box product and the Tychonoff
product coincide.

\begin{proposicion}
    Let $(Y,\mathcal{W})$ be a semi-coarse space, and let $\{(X_\lambda,\mathcal{V}_\lambda)\}_{\lambda \in \Lambda}$
    be a collection of semi-coarse spaces indexed by the set $\Lambda$.
    A map $f:(Y,\mathcal{W})\rightarrow (\prod_{\lambda\in\Lambda}X_\lambda,
    \prod_{\lambda\in\Lambda}\mathcal{V}_\lambda)$ is bornologous iff $\pi_\lambda \circ f:(Y,\mathcal{W}) \to
    (X_\lambda,\mathcal{V}_\lambda)$ is bornologous for each $\lambda\in\Lambda$.
\end{proposicion}

\begin{proof}
    First, suppose that $f$ is bornologous. Since $\pi_\lambda:\left(\prod_{\lambda \in \Lambda} X_\lambda,
    \prod_{\lambda \in \Lambda} \mathcal{V}_\lambda\right)$ is bornologous for any $\lambda \in \Lambda$, 
    it follows that $\pi_\lambda \circ f$ is bornologous for each $\lambda\in\Lambda$, since 
    the composition of bornologous functions is bornologous.

    Now suppose that $\pi_\lambda \circ f$ is bornologous for each $\lambda\in\Lambda$, and let $V\in\mathcal{W}$.
    Then $\pi_\lambda \circ f(V) \in \mathcal{V}_\lambda$. So $\prod_{\lambda\in\Lambda} (\pi_\lambda \circ 
    f(V))\in \prod_{\lambda\in \Lambda}\mathcal{V}_\lambda$. Since $f(V) \subset \prod_{\lambda \in \Lambda}
    (\pi_\lambda \circ f(V))$, $f$ is a bornologous function.
\end{proof}

Our final discussion in this section will be about the relation between products and the functors sending
semi-coarse spaces to semi-uniform spaces and vice-versa. We first recall the definition of the Cartesian 
product of semi-uniform spaces.

\begin{definicion}[Product Semi-Uniform Space; \cite{Cech_1966}, 23.D.10.]
\label{ProdInfSU}
The product of a family $\{(X_a,\mathcal{U}_a):a\in A\}$ of semi-uniform spaces, denoted by
$\prod_{a\in A}(X_a,\mathcal{U}_a)$ is defined to be the semi-uniform space $(X,\mathcal{U})$ where $X$ is 
the cartesian product of the family $\{X_a\}$, and $\mathcal{U}$, called the \emph{product semi-uniformity}, 
is the collection of subsets of $X\times X$ containing a set of the form
\begin{align*}
\{ (x,y)\in X\times X: a\in F\Rightarrow (\pi_a(x),\pi_a(y))\in U_a\}
\end{align*}
where $F$ is a finite subset of $A$ and $U_a\in\mathcal{U}_a$. Sets of the above form are called 
the \emph{canonical elements} of the product semi-uniformity.
\end{definicion}

\begin{teorema}
\label{thm:ProdInfSUPseCoar}
Let $\{(X_a,\mathcal{U}_a)\}$ be a collection of semi-uniform spaces indexed by a set $A$, $\mathcal{V}_a$ the semi-coarse structure induced by $\mathcal{U}_a$, $\mathcal{U}$ the product semi-uniformity and $\mathcal{V}$ the semi-coarse structure induced by $\mathcal{U}$. Then, $\mathcal{V}=\prod_{a\in A}\mathcal{V}_a$, the box
product semi-coarse structure. 
\end{teorema}

\begin{proof}
Let $\{(X_a,\mathcal{U}_a)\}$ be a collection of semi-uniform spaces indexed by a set $A$, $\mathcal{V}_a$ the semi-coarse structure induced by $\mathcal{U}_a$, $\mathcal{U}$ the product semi-uniformity and $\mathcal{V}$ the semi-coarse structure induced by $\mathcal{U}$.

Suppose that $V\in\mathcal{V}$. Then $V$ is contained in every canonical element of the product semi-uniformity,
which implies that $(\pi_a\times\pi_a)(V)\in U_a$ for each $U_a\in\mathcal{U}_a$, 
and therefore $(\pi_a\times\pi_a)(V)\in\mathcal{V}_a$ for each $a\in A$. We conclude that 
$V\in \prod_{a\in A}\mathcal{V}_a$.

On the other hand, if $V\in\prod_{a\in A}\mathcal{V}_a$, then $(\pi_a\times\pi_a)(V)\subset U_a$
for each $U_a\in\mathcal{U}_a$. In particular, for every finite subset $F$ of $A$ we have
$(\pi_a\times\pi_a)(V)\in U_a$ for each $U_a\in\mathcal{U}_a$, $a \in F$, and therefore $U$ is contained in 
each canonical element of the product semi-uniformity. We conclude that $V\in\mathcal{V}$, and the result follows.
\end{proof}

Making the cartesian product of two semi-pseudometric spaces, we get the following and last result.

\begin{lema}
\label{ConmutaProd}
Let $(X,d^X)$ and $(Y,d^Y)$ be metric spaces, $r$ and $s$ non-negative real numbers. Suppose that $d^X_r$ and $d^Y_s$ are
constructed as in \autoref{def:Metric to semipseudometric}, then
\begin{itemize}[wide]
\item[(i)] Let $\mathcal{U}^X_r$, $\mathcal{U}^Y_s$ and $\mathcal{U}_{r,s}$ denote the semi-uniform structures induced by $d^X_r$, $d^Y_s$ and $d_{r,s} \coloneqq\max\{d_r,\delta_s \}$, respectively. Then $\mathcal{U}^X_r\times \mathcal{U}^Y_s = \mathcal{U}_{r,s}$.
\item[(ii)] Let $\mathcal{V}^X_r$, $\mathcal{V}^Y_s$ and $\mathcal{V}_{r,s}$ be the semi-coarse structures induced by $d^X_r$, $d^Y_s$ and $d_{r,s}$, respectively. Then $\mathcal{V}^X_r\times \mathcal{V}^Y_s = \mathcal{V}_{r,s}$.
\end{itemize}
\end{lema}

\begin{proof}
Let $(X,d^X)$ and $(Y,d^Y)$ be metric spaces, let $r,s\geq 0$ be non-negative real numbers, and construct the semi-pseudometrics $d^X_r$ and $d^Y_s$ as in \autoref{def:Metric to semipseudometric}. 
\begin{itemize}[wide]
\item[(i)] Let $\mathcal{U}^X_r$, $\mathcal{U}^Y_s$ and $\mathcal{U}_{r,s}$ be the semi-uniform structures induced by $d^X_r$, $d^Y_s$ and $d_{r,s} \coloneqq \max\{d^X_r,d^Y_s \}$, respectively. First note that, if $W\in \mathcal{U}^X_{r} \times \mathcal{U}^Y_s$, then, by \autoref{ProdInfSU}, there are $U\in\mathcal{U}^X_r$ and $U'\in\mathcal{U}^Y_s$ such that $U\boxtimes U'\subset W$. By the definition of the box product $\boxtimes$, we have that the set $\{(x,x'):d^X_r(x,x')=0\}\boxtimes\{(y,y'):d^Y_s(y,y')=0\}\subset W$ is equal to $\{(x,y,x',y'):d_{r,s}((x,y),(x',y'))=0\}\subset W$. It therefore follows that $W\in\mathcal{D}_{r,s}$.

Conversely, if $W\in\mathcal{V}_{r,s}$, then $\{(x,y,x',y'):d_{r,s}((x,y),(x',y'))=0\}\subset W$, which is equivalent to the condition that $\{(x,x'):d^X_r(x,x')=0\}\boxtimes\{(y,y'):d^Y_s(y,y')=0\}\subset W$. Thus, $W\in\mathcal{U}^X_r\times\mathcal{U}^Y_s$.
\item[(ii)] Let $\mathcal{V}^X_r$, $\mathcal{V}^Y_s$ and $\mathcal{V}_{r,s}$ be the semi-coarse structures induced by $d^X_r$, $d^Y_s$ and $d_{r,s}$, respectively. Then item \emph{(ii)} in the Lemma follows from \autoref{thm:ProdInfSUPseCoar} and item \emph{(i)} above.\qedhere
\end{itemize}
\end{proof}

\subsection{Coarse Space Induced by a Semi-Coarse Space}

In this section, we discuss a method of generating a coarse space from a semi-coarse space. In 
order to proceed with the construction, we will first discuss direct limits for semi-coarse spaces.
We begin by constructing a semi-coarse structure on a disjoint union of semi-coarse spaces.

\begin{proposicion}
    Let $\{(X_\lambda,\mathcal{V}_\lambda)\}_{\lambda \in \Lambda}$ be a collection of semi-coarse 
    spaces indexed by the set $\Lambda$, and let $\sqcup_{\lambda\in\Lambda} \mathcal{V}_\lambda$
    be the collection of sets of the form 
    \begin{equation*}
        \sqcup_{\lambda \in \Lambda} A_\lambda, 
    \end{equation*}
    where each $A_\lambda \in \mathcal{V}_\lambda$. (Note that any given $A_\lambda$ may be the empty set.)
    Then $(\sqcup_{\lambda\in\Lambda} X_\lambda,\sqcup_{\lambda\in\Lambda} \mathcal{V}_\lambda)$ 
    is a semi-coarse space.
\end{proposicion}

\begin{proof}
    We show that $\sqcup_{\lambda \in \Lambda} \mathcal{V}_\lambda$ satisfies the axioms for a
    semi-coarse structure.
    \begin{enumerate}[wide,label=(sc\arabic*)]
    \item $\Delta_{\sqcup_{\lambda\in\Lambda} X_\lambda}= \sqcup_{\lambda\in\Lambda} 
        \Delta_\lambda\in \sqcup_{\lambda\in\Lambda} \mathcal{V}_\lambda$.
    \item If $A \coloneqq \sqcup_{\lambda \in \Lambda} A_\lambda \in \sqcup_{\lambda\in\Lambda}\mathcal{V}_\lambda$ 
        and $W \subset A$, then $W = \sqcup_{\lambda \in \Lambda} W_\lambda$, where 
        $W_\lambda \subset A_\lambda$ for every $\lambda \in \Lambda$. Therefore each 
        $W_\lambda \in \mathcal{V}_\lambda$, and $W\in \sqcup_{\lambda\in\Lambda}\mathcal{V}_\lambda$.
    \item If $A,W\in \sqcup_{\lambda\in\Lambda} \mathcal{V}_\lambda$, then there are
        $\sqcup_{\lambda\in\Lambda} A_\lambda$ and $\sqcup_{\lambda\in\Lambda} 
        W_\lambda$ with $A_\lambda,W_\lambda\in \mathcal{V}_\lambda$ such that $A =
        \sqcup_{\lambda\in\Lambda} A_\lambda$ and $W = \sqcup_{\lambda\in\Lambda} W_\lambda$. 
        So $A\cup W \subset \sqcup_{\lambda\in\Lambda} (A_\lambda\cup W_\lambda)$, from which we conclude
        that $A\cup W\in \sqcup_{\lambda\in\Lambda}\mathcal{V}_\lambda$.
    \item If $A\in\sqcup_{\lambda\in\Lambda} \mathcal{V}_\lambda$, then there is 
        $\sqcup_{\lambda\in\Lambda} A_\lambda$ with $A_\lambda\in\mathcal{V}_\lambda$ such that 
        $A = \sqcup_{\lambda\in\Lambda} A_\lambda$. So 
\begin{align*}
A^{-1} = \left( \sqcup_{\lambda\in\Lambda} A_\lambda \right)^{-1}= \sqcup_{\lambda\in\Lambda} (A_\lambda)^{-1},
\end{align*}
and therefore $A^{-1}\in \sqcup_{\lambda\in\Lambda} \mathcal{V}_\lambda$.
\end{enumerate}
It follows that $(\sqcup_{\lambda\in\Lambda} X_\lambda,\sqcup_{\lambda\in\Lambda} \mathcal{V}_\lambda)$ is a semi-coarse space.
\end{proof}

\begin{definicion}
\label{UnDisjPseCoar}
We call the semi-coarse structure $\sqcup_{\lambda \in \Lambda} \mathcal{V}_\lambda$ the
\emph{disjoint union semi-coarse structure}, and we call the semi-coarse space
$(\sqcup_{\lambda \in \lambda} X_\lambda,\sqcup_{\lambda \in \Lambda} \mathcal{V}_\lambda)$ the 
\emph{disjoint union of the semi-coarse spaces $\{X_\lambda,\mathcal{V}_\lambda\}_{\lambda \in \Lambda}$.}
\end{definicion}

We will now introduce the notion of a directed system in semi-coarse category, which we will use
to construct a coarse space from a semi-coarse space.

\begin{definicion}
Let $\Lambda$ be a directed set. We will call $\{(X_\alpha,\mathcal{V}_\alpha),f_\alpha^\beta,\Lambda\}$
a \emph{directed system of semi-coarse spaces} if $\{X_\alpha, f^\beta_\alpha, \Lambda\}$ is a directed 
system of sets, $(X_\alpha,\mathcal{V}_\alpha)$ are semi-coarse spaces for each $\alpha\in\Lambda$,
and each $f^\beta_\alpha$ is bornologous.
\end{definicion}

\begin{proposicion}
    \label{prop:Direct limit}
Let $\Lambda$ be a directed set, let $\{(X_\alpha,\mathcal{V}_\alpha),f_\alpha^\beta,\Lambda\}$ 
be a directed system of semi-coarse spaces, and let ${\sim}$ be the equivalence relation on 
$\sqcup_{\lambda \in \Lambda} X_\lambda$ such that 
for $x^\alpha\in X^\alpha$ and $x^\beta\in X^\beta$, $x^\alpha\sim x^\beta$ iff
there is a $\gamma\in \Lambda$ satisfying $\alpha\leq \gamma$, $\beta\leq\gamma$ and 
where $f_\alpha^\gamma x^\alpha= f_\beta^\gamma x^\beta$. Then
\begin{align*}
\lim\limits_{\rightarrow} \{(X_\alpha,\mathcal{V}_\alpha),f_\alpha^\beta,\Lambda\} = 
((\sqcup_{\lambda\in\Lambda} X_\lambda)/{\sim},(\sqcup_{\lambda\in\Lambda} \mathcal{V}_\lambda)/{\sim}),
\end{align*}
where the left hand side of the above equation is the direct limit of the directed system, 
and the right hand side is the quotient semi-coarse space from \autoref{def:Quotient qc space}
\end{proposicion}

\begin{proof}
    It is enough to show that, for any $(Y,\mathcal{W})$ and a collection of diagrams of solid arrows
    of the form
\[
\xymatrix{
(X_\lambda,\mathcal{V}_\lambda) \ar[rrrr]^{f_\lambda^\alpha} \ar[rrd]^{p_\lambda} 
\ar[rrdd]_{\psi_\lambda} & & & & 
(X_\alpha,\mathcal{V}_\alpha) \ar[lld]_{p_\alpha} \ar[lldd]^{\psi_\alpha}\\
                         & & \left(\left(\bigsqcup\limits_{\lambda\in\Lambda} X_\lambda\right)/{\sim},
                             \left(\bigsqcup\limits_{\lambda\in\Lambda} 
                         \mathcal{V}_\lambda\right)/{\sim} \right)\ar@{-->}[d]^{g}\\
& & (Y,\mathcal{W})
}
\]
where $\alpha,\lambda\in\Lambda$ satisfy $\lambda\leq \alpha$ and the solid arrows are bornologous and commute, 
there exists a bornologous map 
\begin{equation*}
    g:((\sqcup_{\lambda \in \Lambda} X_\lambda)/{\sim},(\sqcup_{\lambda \in \Lambda} \mathcal{V}_\lambda)/{\sim})
    \to (Y,\mathcal{W})
\end{equation*}
making the entire diagram commute for any choice of $\alpha,\lambda \in \Lambda$ with $\lambda \leq \alpha$. 
However, the existence of a set map
$g$ making the diagram commute is guaranteed by the fact that $\sqcup_{\lambda \in \Lambda} X_\lambda$ is the 
direct limit of the directed system viewed as sets, and $g$ is bornologous by 
\autoref{thm:Bornologous quotient maps}.
\end{proof}

We now characterize the direct limit of some special direct systems which we will use to construct a 
coarse structure from a semi-coarse structure on a set.

\begin{corolario}
\label{LimDirPseCoar}
Let $X$ be a set, suppose that $\Lambda$ is a totally ordered set, and let 
$\{\mathcal{V}_\lambda\}_{\lambda\in\Lambda}$ be a family of semi-coarse structures on
$X$ satisfying $\mathcal{V}_{\lambda}\subset \mathcal{V}_{\lambda'}$ if $\lambda\leq \lambda'$.
Denote by $\{(X_\alpha,\mathcal{V}_\alpha),i_\alpha^\beta,\Lambda\}$ the directed system such that
$X = X_\alpha$ for all $\alpha \in \Lambda$ and the maps $i_\alpha^\beta$ are all identity maps.

Then the direct limit of $\{(X_\alpha,\mathcal{V}_\alpha),i_\alpha^\beta,\Lambda\}$ is the space 
        $\left(X,\bigcup\limits_{\lambda\in\Lambda} \mathcal{V}_\lambda\right)$.
\end{corolario}

\begin{proof}
    We note that $X = (\sqcup_{\lambda \in \Lambda} X_\lambda)/{\sim}$ and 
    $\cup_{\lambda \in \Lambda} \mathcal{V}_\lambda = (\sqcup_{\lambda \in \Lambda} \mathcal{V}_\lambda)/{\sim}$.
        The result now follows from \autoref{prop:Direct limit}.
\end{proof}

\begin{corolario}
\label{LimDirPseCoarCoro}
Let $\{(X_\alpha,V_\alpha),i_\alpha^\beta,\Lambda\}$ be the directed system such that $X = X_\alpha$
for all $\alpha \in \Lambda$, the $i_\alpha^\beta$ are the identity maps, and, in addition, suppose that
the following condition holds.
\begin{equation}
    \label{eq:Coarse condition}
    \begin{split}
        \forall \lambda_1,\lambda_2\in\Lambda \; \exists \lambda_3\in\Lambda
        \text{ such that }
    V\circ W\in\mathcal{V}_{\lambda_3} \;
    \forall\, V \in \mathcal{V}_{\lambda_1}, W 
    \in \mathcal{V}_{\lambda_2}.
\end{split}
\end{equation}
Then $\lim\limits_{\rightarrow} \mathcal{V}_\lambda$ is a coarse structure on $X$.
\end{corolario}

\begin{proof}
    We must show that the semi-coarse structure $\left(\sqcup_{\lambda\in\Lambda} 
    \mathcal{V}_\lambda\right)/{\sim}$ satisfies axiom \ref{item:Coarse condition} in 
    \autoref{def:Semi-coarse}, 
    but this follows immediately from Condition (\ref{eq:Coarse condition}) in the statement of the corollary.
\end{proof}

Given a semi-coarse space $(X,\mathcal{V})$, we will proceed to construct a coarse space by repeatedly adding
sets of the form $V \circ W$ to the semi-coarse structure $\mathcal{V}$. We first
show that adding $\{ V \circ W \mid V,W \in \mathcal{V} \}$ to a semi-coarse structure $\mathcal{V}$
gives another semi-coarse structure. To do so, we will require the following lemma about properties
of the set product.

\begin{lema}
\label{PropProdInv}
Let $X$ and $Y$ be sets and let $A,B,C,D\in\mathcal{P}(X\times X)$. Then
\begin{enumerate}[label=(\roman*)]
\item $(A\circ B)^{-1}=B^{-1}\circ A^{-1}$.
\item \label{item:ii}$(A\circ B)\cup (C\circ D)\subset (A\cup C)\circ (B\cup D)$.
\item \label{item:iii}Let $f:X \to Y$ be a function of sets. Then 
    \begin{equation*}
        (f \times f)(A \circ B) \subset
(f \times f)(A) \circ (f \times f)(B). 
    \end{equation*}
\end{enumerate}
\end{lema}

\begin{proof}
Let $X$ and $Y$ be sets and let $A,B,C,D\in\mathcal{P}(X\times X)$. Then
\begin{enumerate}[wide,label=(\roman*)]
\item Let $(x,y)\in (A\circ B)^{-1}$, that is, $(y,x)\in A\circ B$. This is equivalent to 
    the existence of a point $z\in X$ such that $(y,z)\in A$ and $(z,x)\in B$.
    This, in turn, is true iff $(x,z)\in B^{-1}$ and $(z,y)\in A^{-1}$. 
    By definition, we have $(x,y)\in B^{-1}\circ A^{-1}$.
\item Let $(x,y)\in (A\circ B)\cup (C\circ D)$. Then $(x,y)\in A\circ B$ or $(x,y)\in C\circ D$, 
    and therefore there exists a point $z\in X$ such that $((x,z)\in A \text{ and }(z,y)\in B$
    or $((x,z)\in C \text{ and }(z,y)\in D)$, which implies that 
    $(x,z)\in A\cup C$ and $(z,y)\in B\cup D$. Thus, $(x,y)\in(A\cup C)\circ(B\cup D)$.
\item Suppose that $(a',b')\in (f\times f)(A\circ B)$. Then there exists an element
    pairs $(a,b)\in A\circ B$ such that $f(a)=a'$ and $f(b)=b'$, and therefore there exists
    $x\in X$ with $(a,x)\in A$ and $(x,b)\in B$. It follows that $(a',f(x))\in (f\times f)(A)$
    and $(f(x),b')\in (f\times f)(B)$, which implies that $(a',b')\in (f\times f)(A)\circ (f\times f)(B)$.
    Therefore, $(f\times f)(A\circ B)\subset(f\times f)(A)\circ (f\times f)(B)$, as desired. \qedhere 
\end{enumerate}
\end{proof}

\begin{observacion}
    In \autoref{item:ii} of the previous lemma, the other inclusion is not necessarily true: 
    Let $A$ be a non-empty set, $D=A^{-1}$ and $B=C=\varnothing$, then $A\circ B=C\circ D=\varnothing$ and $(A\cup C)\circ (B\cup D)=A\circ A^{-1}\neq \varnothing$, so that $(A\circ B)\cup (C\circ D)\nsupseteq (A\cup C)\circ(B\cup D)$.
\end{observacion}

\begin{proposicion}
    \label{prop:Set product ext}
Let $(X,\mathcal{V})$ be a semi-coarse space, and define
\begin{align*}
\mathcal{V}^{PE}\coloneqq \{C\subset X\times X \mid \exists A,B\in\mathcal{V}\mbox{ with }C\subset A\circ B\}.
\end{align*}
Then $(X,\mathcal{V}^{PE})$ is a semi-coarse space, and $\mathcal{V} \subset \mathcal{V}^{PE}$.
\end{proposicion}

\begin{proof} We first verify that $\mathcal{V}^{EP}$ satisfies the axioms of a semi-coarse structure. 
\begin{enumerate}[wide,label=(sc\arabic*)]
    \item We observe that $\Delta_X=\Delta_X\circ \Delta_X$. Thus $\Delta_X\in\mathcal{V}^{PE}$.
\item Let $B\in\mathcal{V}^{PE}$ and suppose that $A\subset B$. Then there are sets $A',B'\in\mathcal{V}$ 
    such that $A\subset B\subset A'\circ B'$, so it follows that $A\in\mathcal{V}^{PE}$.
\item Let $A,B\in\mathcal{V}^{PE}$. Then there are $A',B',A'',B''\in\mathcal{V}$
    such that $A\subset A'\circ A''$ and $B\subset B'\circ B''$. Since $A'\cup B', A''\cup B'' \in \mathcal{V}$,
    we have that $(A'\cup B')\circ (A''\cup B'')\in\mathcal{V}^{EP}$. By \autoref{PropProdInv},
    it follows that $A\cup B\subset (A'\cup B')\circ (A''\cup B'')\in \mathcal{V}^{PE}$. We conclude that 
    $A\cup B\in \mathcal{V}^{PE}$.
\item Let $C\in\mathcal{V}^{PE}$. Then there are $A,B\in\mathcal{V}$ such that $C\subset A\circ B$.
    Since $A^{-1},B^{-1}\in\mathcal{V}$, \autoref{PropProdInv} implies that $C^{-1}\subset B^{-1}\circ A^{-1}$,
    and we conclude that $C^{-1}\in \mathcal{V}^{PE}$.
\end{enumerate}
Thus $(X,\mathcal{V}^{PE})$ is a semi-coarse space.

Finally, for any $A \in \mathcal{V}$, $A\circ \Delta_X= A$, and therefore $A \in \mathcal{V}^{PE}$.
\end{proof}

\begin{definicion}
\label{ExtProdConjPseCoar}
Let $(X,\mathcal{V})$ be a semi-coarse space. We call the structure $\mathcal{V}^{EP}$ in 
\autoref{prop:Set product ext} the \emph{set product extension of $\mathcal{V}$}, and 
the ordered pair $(X,\mathcal{V}^{EP})$ will be called the \emph{set product extension of $(X,\mathcal{V})$}.
For any $k \in \N$, we recursively define $\mathcal{V}^{kPE}$ to be the set product extension of 
$\mathcal{V}^{(k-1)PE}$.
\end{definicion}

\begin{proposicion}
    \label{prop:Set product functions}
    If $f:(X,\mathcal{V}) \to (Y,\mathcal{W})$ is bornologous, then $f:(X,\mathcal{V}^{PE}) \to (Y,\mathcal{W}^{PE})$
    is bornologous.
\end{proposicion}

\begin{proof}
If $A\in\mathcal{V}^{PE}$, then there exist sets $A',A''\in\mathcal{V}$ such that 
$A\subset A'\circ A''$, so $(f\times f)(A'),(f\times f)(A'')\in\mathcal{W}$. By \autoref{PropProdInv} 
\autoref{item:iii}, we have 
\begin{equation*}
    (f \times f)(A) \subset (f \times f)(A' \circ A'') \subset (f \times f)(A') \circ (f \times f)(A'') \in
    \mathcal{V}^{PE}.
\end{equation*}
Therefore, $(f \times f)(A) \in \mathcal{V}^{PE}$ and 
$f:X\rightarrow Y$ is a $(\mathcal{V}^{EP},\mathcal{W}^{EP})$-bornologous function.
\end{proof}

\begin{observacion}
    Note that \autoref{prop:Set product functions} shows that the map 
    \begin{align*}\Psi:&\catname{SCoarse} \to \; \catname{SCoarse}\\
        &\Psi(X,\mathcal{V}) = (X,\mathcal{V}^{PE})\\
        &\Psi(f) = f
    \end{align*}
is a functor. 
\end{observacion}

As a final observation, we note that the set product of roofed semi-coarse spaces takes a special form. 

\begin{observacion}
Let $(X,\mathcal{V})$ be a roofed semi-coarse space. Then, for every $A,B\in \mathcal{V}$, $A\circ B \subset \mathfrak{R}(X,\mathcal{V})\circ \mathfrak{R}(X,\mathcal{V})$. Thus, $(X,\mathcal{V}^{EP})$ is a roofed semi-coarse space and $\mathfrak{R}(X,\mathcal{V}^{EP})=\mathfrak{R}(X,\mathcal{V})\circ \mathfrak{R}(X,\mathcal{V})$.
\end{observacion}

We are now ready to construct a coarse structure from a semi-coarse space $(X,\mathcal{V})$.

\begin{teorema}
\label{ExtProdInfinity}
Let $(X,\mathcal{V})$ be a semi-coarse space, and let $\{(X,\mathcal{V}^{kPE}),f_i^j, \mathbb{N}\}_{k=0}^\infty$
be the directed system of semi-coarse structures such that all the 
$f_i^j:(X,\mathcal{V}^{iPE})\rightarrow (X,\mathcal{V}^{jPE})$ are the identity map. 
Then $\varinjlim \{(X,\mathcal{V}^{kPE}),f_i^j, \mathbb{N}\}$ is a coarse space.
\end{teorema}

\begin{proof}
    This follows from \autoref{LimDirPseCoarCoro}, the definition of $\mathcal{V}^{kPE}$, and the fact that 
    $\mathcal{V}^{iPE} \subset \mathcal{V}^{jPE}$ for any $i < j$.
\end{proof}

\begin{definicion}
    We call the coarse structure constructed in \autoref{ExtProdInfinity} the \emph{coarse structure induced 
    by $\mathcal{V}$}, and which we denote by $\mathcal{V}^{\infty}$.
\end{definicion}

Let $(X,d)$ be a metric space, and recall that the set $\mathcal{V}_r \subset X \times X$ is defined by
\begin{equation*}
	\mathcal{V}_r = \{(x,y)\mid d(x,y)\leq r\}.
\end{equation*}
Now denote by $\mathcal{V}_\infty:= \lim_{\rightarrow}\mathcal{V}_r$, the direct limit of the directed system in \autoref{LimDirPseCoar}. 

We now have

\begin{corolario}
\label{Vinfty espcoar}
$(X,\mathcal{V}_\infty)$ is a coarse space.
\end{corolario}

\begin{proof}
Let $\Lambda=(\mathbb{R},\leq)$, let $\{ \mathcal{V}_r, f_r^{r'} ,\Lambda\}$ be the directed system given by the collection of semi-coarse structures induced by a metric $d$. For any $r\leq r'$, define the map $f_r^{r'}=i:(X,\mathcal{V}_r)\rightarrow (X,\mathcal{V}_{r'})$. By definition, $\mathcal{V}_r\subset \mathcal{V}_{r'}$ for $r\leq r'$, and $\mathcal{V}_\infty$ is a semi-coarse structure by \autoref{LimDirPseCoar}.

Now suppose that $V,W\in\mathcal{V}_\infty$. Then there are $r,r'\in [0,\infty)$ such that $V\in\mathcal{V}_r$ and $W\in\mathcal{V}_{r'}$, Without loss of generality, $r\leq r'$, and therefore $V,W\in\mathcal{V}_{r'}$.

So, if $(x,y)\in V\circ W$, then there is $z\in X$ such that $(x,z)\in V$ and $(z,y)\in W$, thereby $d(x,z)\leq r'$ and $d(z,y)\leq r'$. By the triangle inequality property, $d(x,y)\leq 2r'$, and therefor, $$V\circ W \subset \{ (x,y)\in X\times X : d(x,y)\leq 2r'\}.$$
We conclude that $V\circ W\in \mathcal{V}_{2r'}$, and by \autoref{LimDirPseCoarCoro}, $(X,\mathcal{V}_\infty)$ is a coarse space.
\end{proof}

\section*{Acknowledgements}

The authors would like to thank  Omar Antolín Camarena, Noe Barcenas, Jesús Rodríguez, and Alejandra Trujillo for helpful comments on earlier versions of the manuscript. We are also grateful to have had the opportunity to present parts of this work at the annual meeting of the Sociedad Matemática Mexicana, the Topological Data Analysis Workshop at the  Institute for Mathematical and Statistical Innovation, and the GEOTOP-A Seminar. Finally, we would like to thank the anonymous referee for a number of useful comments which helped to improve the manuscript.

\bibliographystyle{amsalpha}
\bibliography{all1,mybib}

\begin{bibdiv}
\begin{biblist}

\bib{Babson_etal_2006}{article}{
      author={Babson, Eric},
      author={Barcelo, H{\'e}l{\`e}ne},
      author={de~Longueville, Mark},
      author={Laubenbacher, Reinhard},
       title={Homotopy theory of graphs},
        date={2006},
        ISSN={0925-9899},
     journal={J. Algebraic Combin.},
      volume={24},
      number={1},
       pages={31\ndash 44},
         url={http://dx.doi.org/10.1007/s10801-006-9100-0},
      review={\MR{2245779}},
}

\bib{Boxer_Staecker_2016}{article}{
      author={Boxer, Laurence},
      author={Staecker, P.~Christopher},
       title={Fundamental groups and {E}uler characteristics of sphere-like
  digital images},
        date={2016},
        ISSN={1576-9402,1989-4147},
     journal={Appl. Gen. Topol.},
      volume={17},
      number={2},
       pages={139\ndash 158},
         url={https://doi-org.biblio.cimat.remotexs.co/10.4995/agt.2016.4624},
      review={\MR{3562604}},
}

\bib{Boxer_Staecker_2017}{article}{
      author={Boxer, Laurence},
      author={Staecker, P.~Christopher},
       title={Homotopy relations for digital images},
        date={2017},
        ISSN={1123-2536,1590-0932},
     journal={Note Mat.},
      volume={37},
      number={1},
       pages={99\ndash 126},
  url={https://doi-org.biblio.cimat.remotexs.co/10.1285/i15900932v37n1p99},
      review={\MR{3733806}},
}

\bib{Boxer_Staecker_2018}{article}{
      author={Boxer, Laurence},
      author={Staecker, P.~Christopher},
       title={Remarks on pointed digital homotopy},
        date={2018},
        ISSN={0146-4124,2331-1290},
     journal={Topology Proc.},
      volume={51},
       pages={19\ndash 37},
      review={\MR{3633236}},
}

\bib{Cech_1966}{book}{
      author={{\v{C}}ech, Eduard},
       title={Topological spaces},
      series={Revised edition by Zden\v ek Frol\' ic and Miroslav Kat\v etov.
  Scientific editor, Vlastimil Pt\'ak. Editor of the English translation,
  Charles O. Junge},
   publisher={Publishing House of the Czechoslovak Academy of Sciences, Prague;
  Interscience Publishers John Wiley \& Sons, London-New York-Sydney},
        date={1966},
      review={\MR{0211373}},
}

\bib{Drutu_Kapovich_2018}{book}{
      author={Dru\c{t}u, Cornelia},
      author={Kapovich, Michael},
       title={Geometric group theory},
      series={American Mathematical Society Colloquium Publications},
   publisher={American Mathematical Society, Providence, RI},
        date={2018},
      volume={63},
        ISBN={978-1-4704-1104-6},
        note={With an appendix by Bogdan Nica},
      review={\MR{3753580}},
}

\bib{Lupton_etal_2023}{article}{
      author={Lupton, Greg},
      author={Oprea, John},
      author={Scoville, Nicholas~A.},
       title={The digital {H}opf construction},
        date={2023},
        ISSN={0166-8641,1879-3207},
     journal={Topology Appl.},
      volume={326},
       pages={Paper No. 108405, 19},
  url={https://doi-org.biblio.cimat.remotexs.co/10.1016/j.topol.2022.108405},
      review={\MR{4533985}},
}

\bib{Lupton_etal_2021}{article}{
      author={Lupton, Gregory},
      author={Oprea, John},
      author={Scoville, Nicholas~A.},
       title={A fundamental group for digital images},
        date={2021},
        ISSN={2367-1726},
     journal={J. Appl. Comput. Topol.},
      volume={5},
      number={2},
       pages={249\ndash 311},
         url={https://doi.org/10.1007/s41468-021-00067-1},
      review={\MR{4259435}},
}

\bib{Lupton_etal_2022a}{article}{
      author={Lupton, Gregory},
      author={Oprea, John},
      author={Scoville, Nicholas~A.},
       title={Homotopy theory in digital topology},
        date={2022},
        ISSN={0179-5376,1432-0444},
     journal={Discrete Comput. Geom.},
      volume={67},
      number={1},
       pages={112\ndash 165},
  url={https://doi-org.biblio.cimat.remotexs.co/10.1007/s00454-021-00336-x},
      review={\MR{4356246}},
}

\bib{Lupton_etal_2022}{article}{
      author={Lupton, Gregory},
      author={Oprea, John},
      author={Scoville, Nicholas~A.},
       title={Subdivision of maps of digital images},
        date={2022},
        ISSN={0179-5376,1432-0444},
     journal={Discrete Comput. Geom.},
      volume={67},
      number={3},
       pages={698\ndash 742},
  url={https://doi-org.biblio.cimat.remotexs.co/10.1007/s00454-021-00350-z},
      review={\MR{4393077}},
}

\bib{Lupton_Scoville_2022}{article}{
      author={Lupton, Gregory},
      author={Scoville, Nicholas~A.},
       title={Digital fundamental groups and edge groups of clique complexes},
        date={2022},
        ISSN={2367-1726,2367-1734},
     journal={J. Appl. Comput. Topol.},
      volume={6},
      number={4},
       pages={529\ndash 558},
  url={https://doi-org.biblio.cimat.remotexs.co/10.1007/s41468-022-00095-5},
      review={\MR{4496690}},
}

\bib{Rieser_2021}{article}{
      author={Rieser, Antonio},
       title={\v{C}ech closure spaces: {A} unified framework for discrete and
  continuous homotopy},
        date={2021},
        ISSN={0166-8641},
     journal={Topology and its Applications},
      volume={296},
       pages={107613},
  url={https://www.sciencedirect.com/science/article/pii/S0166864121000249},
}

\bib{Rieser_arXiv_VR_2020}{misc}{
      author={Rieser, Antonio},
       title={Vietoris-rips homology theory for semi-uniform spaces},
        date={2021},
        note={submitted, arXiv:2008.05739},
}

\bib{Roe_2003}{book}{
      author={Roe, John},
       title={Lectures on coarse geometry},
      series={University Lecture Series},
   publisher={American Mathematical Society, Providence, RI},
        date={2003},
      volume={31},
        ISBN={0-8218-3332-4},
         url={http://dx.doi.org/10.1090/ulect/031},
      review={\MR{2007488}},
}

\bib{Staecker_2021}{article}{
      author={Staecker, P.~Christopher},
       title={Digital homotopy relations and digital homology theories},
        date={2021},
        ISSN={1576-9402,1989-4147},
     journal={Appl. Gen. Topol.},
      volume={22},
      number={2},
       pages={223\ndash 250},
         url={https://doi-org.biblio.cimat.remotexs.co/10.4995/agt.2021.13154},
      review={\MR{4359765}},
}

\bib{Trevino_2019}{thesis}{
      author={Treviño-Marroquín, Jonathan},
       title={Espacios pseuso-coarse},
        type={Master's Thesis},
        date={2019},
}

\bib{Trevino_2024_FundGroupoid_arXiv}{misc}{
      author={Treviño-Marroquín, Jonathan},
       title={Semi-coarse spaces: Fundamental groupoid and the van kampen
  theorem},
        date={2024},
         url={https://arxiv.org/abs/2404.08874},
}

\bib{Zava_2019}{article}{
      author={Zava, Nicol\`o},
       title={Generalisations of coarse spaces},
        date={2019},
        ISSN={0166-8641,1879-3207},
     journal={Topology Appl.},
      volume={263},
       pages={230\ndash 256},
  url={https://doi-org.biblio.cimat.remotexs.co/10.1016/j.topol.2019.05.028},
      review={\MR{3968132}},
}

\end{biblist}
\end{bibdiv}

\end{document}